\documentstyle[11pt]{article}
\begin{document}
\setlength{\textheight}{574pt}
\setlength{\textwidth}{432pt}
\setlength{\oddsidemargin}{18pt}
\setlength{\topmargin}{14pt}
\setlength{\evensidemargin}{18pt}
\newtheorem{theorem}{Theorem}[section]
\newtheorem{sublemma}{Sublemma}[section]
\newtheorem{lemma}{Lemma}[section]
\newtheorem{corollary}{Corollary}[section]
\newtheorem{remark}{Remark}[section]
\newtheorem{definition}{Definition}[section]
\newtheorem{problem}{Problem}
\newtheorem{example}{Example}
\newtheorem{proposition}{Proposition}[section]
\title{{\bf FINITE GENERATION  \\ OF \\ CANONICAL RINGS}}
\date{{\bf September 3, 2000}}
\author{{\bf Hajime Tsuji}}
\maketitle
\begin{abstract}We prove the finite generation of canonical rings of smooth
projective varieties of general type defined over complex numbers.
MSC 32J25
\end{abstract}
\tableofcontents
\section{Introduction}
Let $X$ be a smooth projective variety and let $K_{X}$ be the 
canonical bundle of $X$. 
The graded ring :
\[
R(X,K_{X}):= \oplus_{m=0}^{\infty}H^{0}(X,{\cal O}_{X}(mK_{X}))
\]
is called the canonical ring of $X$. 
$R(X,K_{X})$ is a birational invariant of $X$.
It is conjectured that for every  smooth projective variety $X$,
$R(X,K_{X})$  finitely  generated.
The purpose of this article is to give the following 
partial answer to the conjecture.
\begin{theorem} Let $X$ be a smooth projective varieity 
of general type defined over complex numbers. 

Then  the canonical ring $R(X,K_{X})$ is finitely generated.
\end{theorem}
This theorem has already been known in the case of
 $\dim X\leq 3$(\cite{mu,m}). 

Y. Kawamata pointed that 
the existence of a  {\bf Zariski decomposition}
of $K_{X}$ 
implies the finite generation of $R(X,K_{X})$ (\cite{ka}).
This is our starting point of the proof of Theorem 1.1.
In this case the finite generation is equivalent to 
the stable base point freeness of the nef part of the 
Zariski decomposition. 
In this case the nef part is not ample in general.
Hence to prove the stable base point freeness we need 
to use some {\bf additional positivity}.   
One of the  important observation in \cite{ka} is the fact that 
on the stable fixed component of $K_{X}$, we may use 
the positivity coming from the 
conormal bundle of it.  In \cite{ka}, using a variant of Shokurov's nonvanishing theorem
(\cite{sh}, see Theorem 4.5 below), he  proved the finite generation of canonical rings of 
smooth projective varieties of general type 
by {\bf Noetherian induction} under the assumption that 
there  exists a {\bf Zariski decomposition}
of $K_{X}$. 
The main difficulty to prove the finite generation lies
in such delicate {\bf semipositivity}.

The proof of Theorem 1.1 consists of the following 
five steps.
\begin{enumerate}
\item construct an AZD $h$ of $K_{X}$ to distinguish the positive
part of $K_{X}$,
\item construct the nontrivial numerically trivial fibration associated 
with $(K_{X},h)$ by 
{\bf the first nonvanishing theorem} (Theorem 4.1) on every stable 
fixed components,
\item find an effective {\bf  R}-divisor on a very 
general fiber of the numerically trivial fibration 
by using the structure theorem for numerically trivial 
singular hermitian line bundles (cf. Theorem 4.3).
\item construct the formal canonical model and prove 
{\bf the virtual base point freeness (cf. Definition 6.1)} of 
$R(X,K_{X})$, using {\bf the second nonvanishing theorem} (Theorem 4.4),
\item prove finite generation of $R(X,K_{X})$ by showing that 
the formal canonical model is a projective variety and is the 
canonical model of $X$. 
\end{enumerate}
Let us briefly explain  each steps. \vspace{5mm} \\
{\bf  Step 1}.
An {\bf AZD} $h$ of $K_{X}$ is a singular hermitian metric 
such that the curvature $\Theta_{h}$ is semipositive 
and 
\[
H^{0}(X,{\cal O}_{X}(mK_{X})\otimes{\cal I}(h^{m}))
\simeq H^{0}(X,{\cal O}_{X}(mK_{X}))
\]
holds for every $m\geq 0$. 
The AZD $h$ exists (more generally for any pseudoeffective
line bundle) by \cite{tu,tu2,d}. 
The singular hermitian line bundle $(K_{X},h)$ 
is considered as an analogue of the nef part of a  
Zariski decomposition. \vspace{5mm} \\
{\bf  Step 2}. 
In general $\Theta_{h}$ is not strictly positive. 
$(K_{X},h)$ has only weak positivity on every 
stable fixed components. 
This is one of the main difficulty of the proof.  
To distinguish the {\bf null direction} of $(K_{X},h)$ 
on every stable fixed component, 
we construct the nontrivial {\bf numerically trivial fibrations} 
(cf. Theorem 4.2) associated with $(K_{X},h)$. 
Here the essential ingredient is {\bf the first nonvanishing 
theorem} (Theorem 4.1) which is a generalization of 
Shokurov's nonvanishing theorem (\cite{sh}). 
\vspace{5mm}\\
{\bf  Step 3}.
On a very general fiber of the above numerically trivial fibration, 
by the structure theorem for numerically trivial 
singular hermitian line bundles (cf. Theorem 4.3),
we may distinguish a sum of at most countably many 
prime divisors with nonnegative coefficients.
In the next step, the number of the positive coefficients turns 
out to be finite.  \vspace{5mm} \\
{\bf  Step 4}. 
Taking a successive resolution of $\mbox{Bs}\mid m!K_{X}\mid (m\geq 1)$ , 
and identifying the fiber of the numerically trivial fibrations 
to a point, we costruct a formal canonical model 
$\hat{X}_{can}$. 
At this stage $\hat{X}_{can}$ is only a set. 
But we may consider  $(K_{X},h)$ as a numerically positive ``{\bf  R}-line bundle " on 
$\hat{X}_{can}$.  
Using effective base point freeness argument as in 
\cite{a-s} and {\bf the second nonvanishing theorem}(Theorem 4.4), we may prove {\bf the virtual base point freeness}
of $R(X,K_{X})$. 
The advantage of this argument is that 
{\bf we can specify a point where we want to prove the stably base point freeness of  $(K_{X},h)$  on $\hat{X}_{can}$}.
Hence we do not need to use the  Noetherian 
induction argument as in \cite{ka}. 
This is the effect of the fact that {\bf we may consider $(K_{X},h)$ to be 
 numerically positive on $\hat{X}_{can}$}.  \vspace{5mm} \\
{\bf  Step 5}. 
The last step is to prove that $\hat{X}_{can}$ is actually the canonical model 
of $X$.  To prove this we use a topological argument using 
the virtual base point freensess of $R(X,K_{X})$. \vspace{5mm} \\
We may prove Theorem 1.1 without using singular hermitian metrics. 
In this sense {\bf our proof is essentially algebraic}. 
Although for the better presentation I decided to present the proof 
 in the complex analytic language, one may easily transcript the proof 
in algebro-geometric language. 
In fact the transcription follows from the fact that one can approximate 
any plurisubharmonic functions by plurisubharmonic functions with algebraic singularities (\cite[Section 3]{d}).

In this paper ``{\bf very general}'' means outside of at most countably many union of 
proper Zarski closed subsets  and ``{\bf general}'' means in the sense of usual Zariski topology. 
\section{Multiplier ideal sheaves}

In this section, we shall review the basic definitions and properties
of multiplier ideal sheaves. 
\subsection{Multiplier ideal sheaves}
\begin{definition}
Let $L$ be a line bundle on a complex manifold $M$.
A  singular hermitian metric $h$ on $L$ is given by
\[
h = e^{-\varphi}\cdot h_{0},
\]
where $h_{0}$ is a $C^{\infty}$-hermitian metric on $L$ and 
$\varphi\in L^{1}_{loc}(M)$ is an arbitrary function on $M$.
We call $\varphi$ a  weight function of $h$.
\end{definition}
The curvature current $\Theta_{h}$ of the singular hermitian line
bundle $(L,h)$ is defined by
\[
\Theta_{h} := \Theta_{h_{0}} + \sqrt{-1}\partial\bar{\partial}\varphi ,
\]
where $\partial\bar{\partial}$ is taken in the sense of a current.
The $L^{2}$-sheaf ${\cal L}^{2}(L,h)$ of the singular hermitian
line bundle $(L,h)$ is defined by
\[
{\cal L}^{2}(L,h) := \{ \sigma\in\Gamma (U,{\cal O}_{M}(L))\mid 
\, h(\sigma ,\sigma )\in L^{1}_{loc}(U)\} ,
\]
where $U$ runs open subsets of $M$.
In this case there exists an ideal sheaf ${\cal I}(h)$ such that
\[
{\cal L}^{2}(L,h) = {\cal O}_{M}(L)\otimes {\cal I}(h)
\]
holds.  We call ${\cal I}(h)$ {\bf the multiplier ideal sheaf} of $(L,h)$.
If we write $h$ as 
\[
h = e^{-\varphi}\cdot h_{0},
\]
where $h_{0}$ is a $C^{\infty}$ hermitian metric on $L$ and 
$\varphi\in L^{1}_{loc}(M)$ is the weight function, we see that
\[
{\cal I}(h) = {\cal L}^{2}({\cal O}_{M},e^{-\varphi})
\]
holds.

If $\{\sigma_{i}\}$ are finite number of global holomorphic sections of a line bundle $L$,
for every positive rational number $\alpha$ and a $C^{\infty}$-function 
$\phi$
\[
h := e^{-\phi}\cdot\frac{1}{(\sum_{i}\mid\sigma_{i}\mid^{2})^{\alpha}}
\]
defines a singular hermitian metric  on 
the {\bf Q}-line bundle $\alpha L$. 
Here $\mid \sigma_{i}\mid^{2}$ is defined by 
\[
\mid \sigma_{i}\mid^{2}
= \frac{h_{0}(\sigma_{i},\sigma_{i})}{h_{0}},
\]
where $h_{0}$ is an arbitrary $C^{\infty}$-hermitian metric on 
$L$ (the righthandside is independent of the choice of 
$h_{0}$).
We call such a metric $h$ a singular hermitian metric 
on $\alpha L$ with  {\bf algebraic singularities}.
Singular hermitian metrics with algebraic singularities 
is particulary easy to handle, because its multiplier 
ideal sheaf or that of the multiple of the metric can 
be controlled by taking  suitable successive blowing ups 
such that the total transform of the divisor
$\sum_{i}(\sigma_{i})$ is a divisor with normal crossings. 

Let $D$ be an effective {\bf  R}-divisor on $M$ and let 
\[
\sum_{i} a_{i}D_{i}
\]
be the irreducible decomposition of $D$.
Let $\sigma_{i}$ be a global section of ${\cal O}_{M}(D_{i})$
with divisor $D_{i}$. 
Then 
\[
h = \frac{1}{\prod_{i}\mid\sigma_{i}\mid^{2a_{i}}}
\]
is a singular hermitian metric on the {\bf  R}-line bundle
${\cal O}_{M}(D)$.  
We define the multiplier sheaf ${\cal I}(D)$ associated with $D$ by
\[
{\cal I}(D) = {\cal I}(h) 
= {\cal L}^{2}({\cal O}_{X},\frac{1}{\prod_{i}\parallel \sigma_{i}\parallel^{2a_{i}}}),
\]
where $\parallel\sigma_{i}\parallel$ is the hermitian norm of
$\sigma_{i}$ with respect to a $C^{\infty}$-hermitian metric
on ${\cal O}_{M}(D_{i})$. 

If $\mbox{Supp}\,D$ is a divisor with normal crossings,
\[
{\cal I}(D) = {\cal O}_{M}(-[D])
\]
holds, where $[D] := \sum_{i}[a_{i}]D_{i}$ (for a real 
number $a$, $[a]$ 
denotes the largest integer smaller than or equal to $a$). 

\subsection{Nadel's vanishing theorem}

The following vanishing theorem plays a central role in this paper.

\begin{theorem}(Nadel's vanishing theorem \cite[p.561]{n})
Let $(L,h)$ be a singular hermitian line bundle on a compact K\"{a}hler
manifold $M$ and let $\omega$ be a K\"{a}hler form on $M$.
Suppose that $\Theta_{h}$ is strictly positive, i.e. there exists
a positive constant $\varepsilon$ such that
\[
\Theta_{h} \geq \varepsilon\omega
\]
holds.
Then ${\cal I}(h)$ is a coherent sheaf of ${\cal O}_{M}$ ideal 
and for every $q\geq 1$
\[
H^{q}(M,{\cal O}_{M}(K_{M}+L)\otimes{\cal I}(h)) = 0
\]
holds.
\end{theorem}
\begin{remark} 
The word ``a closed positive $(1,1)$ current" 
does not mean a closed  strictly positive current. 
For example the $0$-current is closed positive. 
This terminology might be misleading for algebraic geometers. 
\end{remark}
By the definition of a multiplier ideal sheaf we have the following
lemma which will be used later.
\begin{lemma}
Let $(L,h)$ be a singular hermitian line bundle on a complex 
manifold $M$.
Let $f : N\longrightarrow M$ be a modification. 
Then $(f^{*}L,f^{*}h)$ is a singular hermitian line bundle on $N$
and 
\[
f_{*}{\cal I}(f^{*}h) \subseteq {\cal I}(h)
\]
holds.
\end{lemma}

\subsection{Lelong numbers and structure of closed positive $(1,1)$-currents}
A closed positive $(1,1)$-current is considered as 
a $(1,1)$-form whose coefficients are distributions. 
Hence by the Lebesugue decomposition 
of the coefficients, every closed positive $(1,1)$-current $T$ on 
a complex manifold $M$ is uniquely decomposed as :
\[
T = T_{abc} + T_{sing},
\]
where $T_{abc}$ denotes the absolutely continuous part and 
$T_{sing}$ denotes the singular part. 
We call this decomposition the {\bf Lebesgue decomposition} of $T$.
It is important to note that $T_{abc}$ and $T_{sing}$ are 
not closed in general. 
To measure the magnitude of the singular part, 
the following definition is fundamental. 
\begin{definition}
Let $T$ be a closed positive $(1,1)$-current on a 
a unit open polydisk $\Delta^{n}$ with center $O$.
Then by $\partial\bar{\partial}$-Poincar\'{e} lemma
there exists a plurisubharmonic function  $\varphi$ 
on $\Delta^{n}$ such that
\[
T = \frac{\sqrt{-1}}{\pi}\partial\bar{\partial}\varphi .
\]
We define the Lelong number  $\nu (T,O)$ at $O$ by
\[
\nu (T,O) = \liminf_{x\rightarrow O}\frac{\varphi (x)}{\log \mid x\mid},
\]
where $\mid x\mid  = (\sum\mid x_{i}\mid^{2})^{1/2}$.
It is easy to see that $\nu (T,O)$ is independent of the choice of
$\varphi$ and local coordinates around $O$.
For an analytic subset $V$ of a complex manifold $X$, we set 
\[
\nu (T,V) = \inf_{x\in V}\nu (T,x).
\]
\end{definition}
\begin{remark} More generally 
the Lelong number is defined for a closed positive
$(k,k)$-current on a complex manifold.
\end{remark}
Let us consider  a singular hermitian metric on the trivial bundle with algebraic singularities 
\[
h = (\sum_{i=1}^{k} \mid f_{i}\mid^{2})^{-1} (f_{i}\in {\cal O}(\Delta^{n}))
\]
on $\Delta^{n}$.
Then we see that $\Theta_{h}$ is positive and for every $x\in \Delta^{n}$. 
\[
\nu(\Theta_{h},x) = 2\min_{i}\mbox{mult}_{x}(f_{i})
\]
holds, where $(f_{i})$ denotes the divisor of $f_{i}$ for every $i$.
This implies that for a singular hermitian metric with algebraic singularities, the Lelong number of the curvature is essentially the infimum of 
the vanishing order of the defining (multi)sections. 
  
In this paper we only deal with singular hermitian metrics which is 
a limit of singular hermitian metrics with algebraic singularities.
Hence in this paper we may consider that {\bf the Lelong number is nothing but the limit vanishing order of (multi)sections} in an obvious manner. 

The following theorem is fundamental. 
\begin{theorem}(\cite[p.53, Main Theorem]{s})
Let $T$ be a closed positive $(k,k)$-current on a complex manifold
$M$.
Then for every $c > 0$
\[
\{ x\in M\mid \nu (T,x)\geq c\}
\]
is a subvariety of codimension $\geq k$
in $M$.
\end{theorem}
Let $(L,h)$ be a singular hermitian line bundle on a smooth
projective variety $X$ such that $\Theta_{h}$ is a positive current.
The following lemma shows a rough relationship between 
the Lelong number of $\nu(\Theta_{h},x)$ at $x\in X$ and the stalk of the multiplier
ideal sheaf ${\cal I}(h)_{x}$ at $x$. 

\begin{lemma}(\cite[p.284, Lemma 7]{b}\cite{b2},\cite[p.85, Lemma 5.3]{s})
Let $\varphi$ be a plurisubharmonic function on 
the open unit polydisk $\Delta^{n}$ in $\mbox{\bf C}^{n}$ with center $O$.
Suppose that $e^{-\varphi}$ is not locally integrable 
around $O$.
Then we have that
\[
\nu (\sqrt{-1}\partial\bar{\partial}\varphi ,O)\geq 2
\]
holds.
And if
\[
\nu (\sqrt{-1}\partial\bar{\partial}\varphi ,O) > 2n
\]
holds, then $e^{-\varphi}$ is not locally integrable around $O$.
\end{lemma}

Let $T$ be a closed positive $(1,1)$-current on a complex manifold
$X$.  
Let ${\cal U} = \{ U_{\alpha}\}$ be an open covering of $X$ such that
for every $\alpha$ 
there exists a plurisubharmonic function $\varphi_{\alpha}$ such that
\[
T\mid U_{\alpha} = \sqrt{-1}\partial\bar{\partial}\varphi_{\alpha}
\]
holds.
We define the singular set $\mbox{Sing}\, T$ by
\[
\mbox{Sing}\, T\cap U_{\alpha}= 
\{ x\in U_{\alpha}\mid \varphi_{\alpha}(x) = -\infty\} .
\]
$\mbox{Sing}\, T$ is well defined and independent of the choice 
of $\{ U_{\alpha}\}$ and $\{ \varphi_{\alpha}\}$.
Let $Y$ be a complex manifold and let 
$f : Y\longrightarrow  X$ be a holomorphic map such that 
\[
f(Y)\not{\subseteq} \mbox{Sing}\, T.
\]
Then the pullback $f^{*}T$ is defined by 
\[
f^{*}\mid f^{-1}(U_{\alpha}) 
= \sqrt{-1}\partial\bar{\partial}(f^{*}\varphi_{\alpha}).
\]
If $Y$ is a submanifold of $X$ and $f$ is the canonical immersion,
then we denote $f^{*}T$ by $T\mid_{Y}$ and call it the restriction
of $T$ to $Y$.
To compute the Lelong number the following lemma is useful.
\begin{lemma}(\cite{s})
Let $T$ be a closed positive $(1,1)$-current on the open 
unit polydisk $\Delta^{n}$ with center $O$. 
Let us parametrize the lines passing through $O$ by 
${\bf  P}^{n-1}$ in the standard way. 
Then there exists a set $E$ of measure $0$ in 
${\bf  P}^{n-1}$ such that for every $[L]\in {\bf  P}^{n-1}-E$, 
$T\mid_{L\cap\Delta^{n}}$ is well defined and 
\[
\nu (T,O) = \nu (T\mid_{L\cap\Delta^{n}},O)
\]
holds. 
\end{lemma}
The next corollary is analogous to the corresponding fact 
about multiplicities of divisors. 
\begin{corollary}
Let $M$ be a complex manifold and let $T$ be a closed 
positive $(1,1)$-current on $M$. 
Let $f : N\longrightarrow M$ be a composition of  successive blowing ups 
with smooth centers.
Then for every $x\in M$ and $y\in f^{-1}(x)$, 
\[
\nu (f^{*}T,y) \geq \nu (T,x)
\]
holds. 
\end{corollary}
\section{Analytic Zariski decomposition}

In this section we shall introduce the notion of  analytic Zariski decomopositions which play  essential roles in this paper. 
By using analytic Zariski decompositions, we can handle a big line bundles
as if it were a nef and big line bundles.
\subsection{Definition of AZD}
\begin{definition}
Let $M$ be a compact complex manifold and let $L$ be a line bundle
on $M$.  A singular hermitian metric $h$ on $L$ is said to be 
an  analytic Zariski decomposition (AZD), if the followings hold.
\begin{enumerate}
\item $\Theta_{h}$ is a closed positive current,
\item for every $m\geq 0$, the natural inclusion
\[
H^{0}(M,{\cal O}_{M}(mL)\otimes{\cal I}(h^{m}))\rightarrow
H^{0}(M,{\cal O}_{M}(mL))
\]
is an isomorphim.
\end{enumerate}
\end{definition}
\begin{remark} If an AZD exists on a line bundle $L$ on a smooth projective
variety $M$, $L$ is pseudoeffective by the condition 1 above.
\end{remark}
As for the existence the following theorems are known. 
\begin{theorem}(\cite{tu,tu2})
 Let $L$ be a big line  bundle on a smooth projective variety
$M$.  Then $L$ has an AZD. 
\end{theorem}
More generally  the existence for general pseudoeffective line bundles, 
now we have the following theorem.
\begin{theorem}(\cite{d-p-s}, cf. \cite[Theorem 2.4]{tu4})
Let $X$ be a smooth projective variety and let $L$ be a pseudoeffective 
line bundle on $X$.  Then $L$ has an AZD.
\end{theorem}
\subsection{An explicit construction of AZD}
In the case of big line bundles, we have an explicit construction of 
an AZD (\cite{tu}). 
Here we shall review the construction.
\begin{definition} 
Let $(L,h)$ be a line bundle on a compact K\"{a}hler manifold $(X,\omega)$. 
Let $\{\phi_{0},\ldots ,\phi_{N(m)}\}$ be an orthonormal basis of 
$\Gamma (X,{\cal O}_{X}(mL))$ with respect to the $L^{2}$-inner 
product
\[
(\phi ,\phi^{\prime}):= \int_{X}h^{m}\phi\cdot\bar{\phi}^{\prime}\frac{\omega^{n}}{n!} \,\, (\phi,\phi^{\prime}\in \Gamma (X,{\cal O}_{X}(mL))).
\]
We define the $m$-th Bergman kernel $K_{m}(z,w)$ of $(L,h)$
by 
\[
K_{m}(z,w):= \sum_{i=0}^{N(m)}\phi_{i}(z)\bar{\phi}_{i}(w).
\]
Then it is trivial to see that $K_{m}(z,w)$ is independent of the 
choice of the orthonormal basis. 
For simplicity we denote the restriction of $K_{m}(z,w)$ to the 
diagonal of $X\times X$ by $K_{m}$. 
\end{definition}
\begin{theorem}(\cite{tu})
Let $L$ be a big line bundle on a compact K\"{a}hler manifold 
$(X,\omega )$.  Let $h_{0}$ be a $C^{\infty}$-hermitian metric on 
$L$.  Let $K_{m}(z,w)$ be the $m$-the Bergman kernel of $(L,h_{0})$.
Then 
\[
h := (\overline{\lim}_{m\rightarrow\infty}\sqrt[m]{K_{m}})^{-1}
\]
is an AZD of $L$. 
\end{theorem}
The reason why we presented this explicit construction 
is to show that {\bf the AZD is  a limit of singular hermitian 
metrics $\{ h_{m}\}$ with algebraic singularities} on $L$ defined by
\[
h_{m} := 1/\sqrt[m]{K_{m}}.
\] 
The Lelong number of $\Theta_{h}$ is considered as  
a limit of the Lelong number of $\Theta_{h_{m}}$ 
which is  nothing but 
\[
m^{-1}\mbox{mult}\,\mbox{Bs}\mid mL\mid .
\]
Hence for every $x\in X$ 
\[
\nu (\Theta_{h},x) = \overline{\lim}_{m\rightarrow\infty}m^{-1}\mbox{mult}_{x}\,\mbox{Bs}\mid mL\mid 
\]
holds, where $\mbox{mult}_{x}\,\mbox{Bs}\mid mL\mid$ 
denotes the multiplicity of the general member of 
$\mid mL\mid$ at $x$. 
Hence $\nu (\Theta_{h},x) (x\in X)$ is essentially an 
algebro-geometric number. 
\subsection{Some properties of AZD}
Let $(L,h)$ be a singular hermitian line bundle on a smooth projective variety $X$. 
We denote the linear system $\mid H^{0}(X,{\cal O}_{X}(mL)\otimes{\cal I}(h^{m}))\mid$ by $\mid m(L,h)\mid$.  
\begin{theorem}
Let $L$ be a big line bundle on a smooth projective 
variety $X$ and let $h$ be an 
AZD of $L$. 
Then there exists a positive constant $C$ such that 
\[
0\leq \mbox{mult}_{x}\mbox{Bs}\mid m(L,h)\mid -m\cdot\nu (\Theta_{h},x) 
\leq C
\]
holds for every $m$ and $x\in X$.
\end{theorem}
{\bf Proof}. 
The first inequality is trivial by the definition of an AZD and 
the fact that $R(X,L)$ is a ring. 
In fact by Lemma 2.2 we see that 
\[
\mbox{mult}_{x}\mbox{Bs}\mid m(L,h)\mid 
\geq m\cdot \nu (\Theta_{h},x) - n
\]
holds. 
Since $h$ is an AZD of $L$, $\mid m(L,h)\mid = \mid mL\mid$ 
holds for every $m\geq 0$. 
Hence for every $\sigma \in \Gamma (X,{\cal O}_{X}(mL)) - \{ 0\}$
and a positive integer $\ell$, 
\[
\ell \cdot \mbox{mult}_{x}(\sigma )
= \mbox{mult}_{x}(\sigma^{\ell})
\geq \ell m\cdot \nu(\Theta_{h},x)
\]
holds. 
Dividing both sides by $\ell$ and letting $\ell$ tend
to infinity, we see that 
\[
\mbox{mult}_{x} (\sigma ) \geq m \cdot \nu (\Theta_{h},x)
\]
holds. 

Nest we shall verify the second inequality.
Let $x$ be a point on $X$. 
Let $\pi : \tilde{X}\longrightarrow X$ 
be the blowing up of $X$ at $x$. 
Since $L$ is big, by Kodaira's lemma (cf. \cite[Appendix]{k-o}) there exists an 
effective {\bf  Q}-divisor $E$ such that 
$\pi^{*}L - E$ is ample.   Let $r$ be a sufficiently large 
posiive integer such that 
\[
H := r(\pi^{*}L -E)
\]
is Cartier and $H- K_{\tilde{X}}$ is  ample. 
Let $\tilde{x}$ be a very general point on the 
exceptional divisor $\pi^{-1}(x)$. 
Then by Theorem 2.2 and Lemma 2.2 we may assume that
for every $m\geq 0$  
the multiplier ideal sheaf ${\cal I}(\pi^{*}h^{m})$  
is locally free on a neighbourhood of $\tilde{x}$ 
(the neighbourhood may depend on $m$). 
Let $U$ be a small neighbourhood of $\tilde{x}$ 
and let $\rho$ be a $C^{\infty}$-function on $\tilde{X}$ 
such that 
\begin{enumerate}
\item $\mbox{Supp}\,\rho \subset\subset U$,
\item $\rho\equiv 1$ on a neighbourhood of $\tilde{x}$,
\item $0 \leq \rho \leq 1$
\end{enumerate}
hold. 
Let $d_{\tilde{x}}$ denote the distance function from $\tilde{x}$ 
with respect to a fixed K\"{a}hler metric on $\tilde{X}$. 
If we take $r$ sufficiently large we may assume that 
there exists a $C^{\infty}$ hermitian metric $\tilde{h}$ on 
$H - K_{\tilde{X}}$ such that 
\[
\Theta_{\tilde{h}} +
 2n\sqrt{-1}\partial\bar{\partial}\log (\rho\cdot d_{\tilde{x}})
\]
is strictly positive on $\tilde{X}$. 
Then by Nadel's vanishing theorem we see that 
\[
H^{0}(\tilde{X},{\cal O}_{\tilde{X}}(H + \pi^{*}(mL))
\otimes {\cal I}(\pi^{*}h^{m}))
\rightarrow {\cal O}_{\tilde{X}}(H+ \pi^{*}(mL))\otimes
{\cal I}(\pi^{*}h^{m})\otimes 
{\cal O}_{\tilde{X}}/{\cal M}_{\tilde{x}}
\]
is surjective for every $m\geq 0$, where ${\cal M}_{\tilde{x}}$ denotes 
the maximal ideal sheaf at $\tilde{x}$. 
Since $h$ is an AZD of $L$, we see that 
there exists a canonical injection
\[
H^{0}(\tilde{X},{\cal O}_{\tilde{X}}(H + \pi^{*}(mL))
\otimes {\cal I}(\pi^{*}h^{m}))
\hookrightarrow \pi^{*}H^{0}(X,{\cal O}_{X}(mL)). 
\]
Since 
\[
\pi_{*}({\cal O}_{\tilde{X}}(K_{\tilde{X}}+mL)\otimes{\cal I}(\pi^{*}h^{m}))
= {\cal O}_{X}(K_{X}+mL)\otimes{\cal I}(h^{m})
\]
holds by the definition of multiplier ideal sheaves, 
we see that 
\[
{\cal I}(h^{m})\otimes {\cal M}_{x}^{n} \subset \pi_{*}{\cal I}(\pi^{*}h^{m})
\]
holds for every $m$. 
Hence by the above argument, there exists a positive constant $C$ such that 
\[
\mbox{mult}_{x}\mbox{Bs}\mid m(L,h)\mid 
-m\cdot \nu (\Theta_{h},x) \leq C
\]
holds for every $m$. 
It is easy to see that $C$ can be taken independent of $x\in X$.
This completes the proof of Theorem 3.4. 
{\bf  Q.E.D}

\begin{corollary} 
Let $P$ be a nef and big line bundle on a smooth projective
variety $X$, then there exists a singular hermitian metric 
$h_{P}$ on $P$ such that 
$\Theta_{h_{P}}$ is a closed positive current on $X$ and 
$\nu (\Theta_{h_{P}})$ is identically $0$ on $X$.
\end{corollary}

We shall discuss about the uniqueness of the multiplier ideal sheaves  
associated with an AZD. 
First we introduce the following teminology.  
\begin{definition}
Let $h_{L}$ be a singular hermitian metric on a line bundle $L$
on a complex manifold $X$.
Suppose that the curvature of $h_{L}$ is a positive current on $X$.
We set 
\[
\bar{\cal I}(h_{L}) 
:= \lim_{\varepsilon\downarrow 0}{\cal I}(h_{L}^{1+\varepsilon})
\]
and call it  the closure of ${\cal I}(h_{L})$. 
\end{definition} 
Let us explain the reason why we take the closure.
Let $h_{L}$ be a singular hermitian metric on a line bundle $L$ 
on a complex manifold $X$ with positive curvature current. 
Then $\bar{\cal I}(h_{L})$ is coherent ideal sheaf on $X$ by Theorem 2.1.
Let $f : Y \longrightarrow X$ be a modification such that 
$f^{*}\bar{\cal I}(h_{L})$ is locally free. 
If we take $f$ properly, we may assume that there exists a divisor $F = \sum F_{i}$ with normal crossings on $Y$ such that 
\[
K_{Y} = f^{*}K_{X} + \sum a_{i}F_{i}
\]
and 
\[
f^{*}\bar{\cal I}(h_{L}) = {\cal O}_{Y}(-\sum b_{i}F_{i})
\]
hold on $Y$ for some nonnegative integers $\{ a_{i}\}$ and $\{ b_{i}\}$. 
Then by Lemma 2.2, 
\[
b_{i} = [\nu (f^{*}\Theta_{h_{L}},F_{i}) - a_{i}]
\]
holds for every $i$. 
In this way $\bar{\cal I}(h_{L})$ is determined by the Lelong numbers of 
the curvature current on some modification. 
This is not the case, unless we take the closure as in the following example. 
\begin{example}
Let $h_{P}$ be a singular hemritian metric on the trivial line bundle on {\bf  C} 
defined by
\[
h_{P} = \frac{1}{\mid z\mid^{2}(\log \mid z\mid  )^{2}}.
\]
Then $\nu (\Theta_{h_{P}},0) = 1$ holds. 
But ${\cal I}(h_{P}) = {\cal O}_{\bf  C}$ holds. 
On the other hand $\bar{\cal I}(h_{P}) = {\cal M}_{0}$ holds. 
\end{example}
\begin{remark}
Theorem 2.1 still holds, even if we replace the multiplier ideal by its closure. 
This can be verified as follows. 
Let $X$ be a compact K\"{a}hler manifold and let $(L,h_{L})$ be a 
singular hermitian line bundle on $X$ with strictly positive curvature. 
Let $h_{\infty}$ be a $C^{\infty}$-hermitian metric on $L$. 
Then $h_{L}^{1+\varepsilon}\cdot h_{\infty}^{-\varepsilon}$ has strictly positive curvature for every sufficiently small positive number $\varepsilon$. 
By Theorem 2.1, we see that 
\[
H^{q}(X,{\cal O}_{X}(K_{X}+L)\otimes {\cal I}(h_{L}^{1+\varepsilon})) = 0
\]
holds for every $q\geq 1$ and every sufficiently small positive number $\varepsilon$.
Letting $\varepsilon\downarrow 0$ we obtain that 
\[
H^{q}(X,{\cal O}_{X}(K_{X}+L)\otimes\bar{\cal I}(h_{L})) = 0
\]
holds for every $q\geq 1$. 
\end{remark}
Now we shall prove the following uniqueness theorem for the multiplier 
ideal sheaves associated with an AZD. 
\begin{proposition}
Let $L$ be a big line bundle on a smooth projective variety $X$.
Let $h$ be an AZD of $L$. 
For any positive integer $m$,  
$\bar{\cal I}(h^{m})$ is independent of the choice of the AZD $h$.
\end{proposition}
{\em Proof of Proposition 3.1}
Let $L$,$h$ be as above. 
By Theorem 3.4,  for any modification 
\[
f : Y \longrightarrow X
\]
and $y\in Y$, we see that 
\[
\nu (f^{*}\Theta_{h},y) = \lim_{m\rightarrow\infty}m^{-1}
\mbox{mult}_{y}f^{*}\mbox{Bs}\mid mL\mid
\]
holds.  
This implies that for any positive integer $m$
$\bar{\cal I}(h^{m})$ is independent of the choice of $h$.
{\bf Q.E.D.}

\begin{definition}
Let $L$ be a pseudoeffective line bundle on a smooth projective 
variety $X$ and let $h$ be a singular hermitian metric on $L$ 
with positive curvature current.  
Let $x$ be a point on $X$. 
$\mid m(L,h)\mid$ is said to be base point free at $x$, if 
\[
\mbox{mult}_{x}\mbox{Bs}\mid m(L,h)\mid = m\cdot \nu (\Theta_{h},x)
\]
holds. 
\end{definition}
Let $h$ be an AZD on a big line bundle $L$ on a smooth projective 
variety $X$. 
Then by Theorem 3.4, we see that 
$(L,h)$ is {\bf asymptotically base point free} 
in the context of Definition 3.4. 
\subsection{Volume of subvarieties}
Let $L$ be a big line bundle on a smooth projective variety
$X$. To measure the total positivity of $L$ on 
a subvariety of $X$. We define the following notion.
\begin{definition}(\cite{tu3})
Let $L$ be a big line bundle on a smooth projective variety
$X$ and let $h$ be an AZD of $L$.
Let $Y$ be a subvariety of $X$ of dimension $r$.
We define the  volume $\mu (Y,L)$ of $Y$ with respect to 
$L$ by 
\[
\mu (Y,L) := r!\cdot\overline{\lim}_{m\rightarrow\infty}m^{-r}
\dim H^{0}(Y,{\cal O}_{Y}(mL)\otimes{\cal I}(h^{m})/tor),
\]
where $tor$ denotes the torsion part of ${\cal O}_{Y}(mL)\otimes{\cal I}(h^{m})$.
\end{definition}
\begin{remark}
If we define $\mu (Y,L)$ by 
\[
\mu (Y,L) := r!\cdot\overline{\lim}_{m\rightarrow\infty}m^{-r}
\dim H^{0}(Y,{\cal O}_{Y}(mL))
\]
then it is totally different unless $Y = X$. 
The above definition is meaningful when ${\cal O}_{Y}(mL)\otimes{\cal I}(h^{m})$ is generically rank one for $m >> 1$. 
Otherwise $\mu (Y,L)$ may be infinity.
If $\mu (Y,L)$ is finite, by Proposition 3.1, it is easy to see that 
$\mu (Y,L)$ is independent of 
the choice of the AZD $h$. 
In fact if $\mu (Y,L) > 0$, then there exists a positive integer $m_{0}$ 
and an effective divisor $E$ on $Y$ so that for every nonnegative integer $m$,the natural injection
\[
{\cal O}_{Y}((m+m_{0})L -E)\otimes {\cal I}(h^{m+m_{0}})
\rightarrow {\cal O}_{Y}(mL)\otimes {\cal I}(h^{m})
\]
exists.   
\end{remark}
If $L$ is a nef and big line bundle on a smooth projective
variety, then by Corollary 3.1 and 
Lemma 2.2, for every subvariety $Y$ in $X$, 
\[
\mu (Y,L) = L^{\dim Y}\cdot Y
\]
holds. 
For a general singular hermitian line bundle with positive curvature, 
we define the volume as follows. 
\begin{definition}
$L$ be a pseudoeffective line bundle on a smooth projective variety $X$
and let $h$ be a singular hermitian metric on $X$ such that 
$\Theta_{h}$ is a closed positive current. 
Let $Y$ be a subvariety of $X$.  
We define the volume $\mu (Y,(L,h))$ of $Y$ with respect to $(L,h)$ by 
\[
\mu (Y,(L,h)) := (\dim Y)!\cdot\overline{\lim}_{m\rightarrow\infty}
m^{-\dim Y}\dim H^{0}(Y,{\cal O}_{Y}(mL)\otimes{\cal I}(h^{m})/tor). 
\]
\end{definition}
\subsection{Intersection theory for singular hermitian line bundles}
In this subsection we review the definition an intersection number for 
a singular hermitian line bundle with positive curvature current 
on a smooth projective variety and an irreducible curve 
on it (cf. \cite{tu4}).  
This intersection number is different from the usual intersection 
number of the underlying line bundle  and the curve. 
The new intersection number measures 
the intersection of the {\bf positive part} of the singular hermitian 
line bundle and the curve. 
Next we shall consider the restriction of singular hermitan 
line bundles to subvarieties. 
\begin{definition}
Let $L$ be a line bundle on a complex manifold $M$.
Let $h$ be a singular hermitian metric on $L$ given by 
\[
h = e^{-\varphi}\cdot h_{0},
\]
where $h_{0}$ is a $C^{\infty}$-hermitian metric on $L$ and 
$\varphi\in L^{1}_{loc}(M)$.
Suppose that the curvature current $\Theta_{h}$ is bounded 
from below by some $C^{\infty}$-(1,1)-form. 
For a subvariety $V$ of $M$, we say that the restriction 
$h\mid_{V}$ is well defined, if 
$\varphi$ is not identically $-\infty$ on $V$. 
\end{definition}
Let $(L,h)$,$h_{0}$,$V$, $\varphi$ be as in Definition 3.7.
Then $\varphi$ is an almost plurisubharmonic 
function i.e. locally a sum of a plurisubharmonic function and $C^{\infty}$-function.
Let $\pi :\tilde{V} \longrightarrow V$ be an arbitrary 
resolution of $V$. 
Then  $\pi^{*}(\varphi\mid_{V})$ is locally integrable on $\tilde{V}$, 
since $\varphi$ is almost plurisubharmonic. 
Hence 
\[
\pi^{*}(\Theta_{h}\mid_{V}) := \Theta_{\pi^{*}h_{0}\mid_{V}}
+ \sqrt{-1}\partial\bar{\partial}\pi^{*}(\varphi\mid_{V})
\]
is well defined.
We shall define the intersection number 
for a singular hermitian metric with positive curvature current 
and an irreducible curve such that the restriction of 
the singular hermitian metric is well defined. 
\begin{definition}
Let $(L,h)$ be a singular hermitian lien bundle on a smooth 
projective variety $X$ such that the curvature current $\Theta_{h}$
is closed positive. 
Let $C$ be an irreducible curve on $X$ such that 
$h\mid_{C}$ is well defined. 
The  intersection number $(L,h)\cdot C$ is defined by
\[
(L,h)\cdot C := 
\overline{\lim}_{m\rightarrow\infty}m^{-1}\dim H^{0}(C,{\cal O}_{C}(mL)\otimes
{\cal I}(h^{m})/tor),
\]
where $tor$ denotes the torsion part of 
${\cal O}_{C}(mL)\otimes {\cal I}(h^{m})$. 
\end{definition}
Let $(L,h)$,$C$ be as above.
Let 
\[
\pi : \tilde{C}\longrightarrow C
\]
be the normalization of $C$. 
We define the multiplier ideal sheaf \\
${\cal I}(h^{m}\mid_{C}) (m\geq 0)$ 
on $C$ by 
\[
{\cal I}(h^{m}\mid_{C}) := \pi_{*}{\cal I}(\pi^{*}h^{m}\mid_{C}). 
\]
And the Lelong number $\nu (\Theta_{h}\mid_{C},x) (x\in C)$ by 
\[
\nu (\Theta_{h},x) = \sum_{\tilde{x}\in \pi^{-1}(x)}
\nu (\pi^{*}\Theta_{h}\mid_{C},\tilde{x}). 
\]
\begin{lemma}(\cite[Lemma 2.4]{tu4})
Let $(L,h)$ be a singular hermitian line bundle on 
a smooth projective variety $X$ such that 
$\Theta_{h}$ is closed positive.  
Let $C$ be an irreducible curve on $X$ such that $h\mid_{C}$ 
is well defined. 
Suppose that $(L,h)\cdot C = 0$ holds. 
Then 
\[
\Theta_{h}\mid_{C} = \sum_{x\in C}\nu (\Theta_{h}\mid_{C},x)x
\]
holds in the sense that 
\[
\pi^{*}(\Theta_{h}\mid_{C}) = 
\sum_{\tilde{x}\in C}\nu (\pi^{*}\Theta_{h}\mid_{C},\tilde{x})\tilde{x}
\]
holds. 
\end{lemma}
\begin{definition} 
Let $(L,h)$ be a singular hermitian line bundle on 
a smooth projective variety $X$ such that $\Theta_{h}$ 
is positive. 
$(L,h)$ is said to be  numerically trivial, if for every 
irreducible curve $C$ on $X$ such that $h\mid_{C}$ is 
well defined, 
\[
(L,h)\cdot C = 0
\]
holds. 
\end{definition}

\subsection{Restriction of the intersection theory to 
divisors}
In the previous subsection we define an 
intersection number of a singular hemitian line bundle 
with positive curvature and an irreducible curve on 
which the restriction of the singular hermitian metric 
is well defined.  In this subsection, we shall 
extend the definition of the intersection number. 

Let $(L,h)$ be a singular hermitian line bundle on 
a smooth projective variety $X$ such that $\Theta_{h}$ 
is positive. 

Let $D$ be a {\bf smooth} divisor on $X$. 
We set  
\[
v_{m}(D) := \mbox{mult}_{D}\mbox{Spec}({\cal O}_{X}/{\cal I}(h^{m}))
\]
and 
\[
\tilde{\cal I}_{D}(h^{m}) := {\cal O}_{D}(v_{m}(D)D)\otimes {\cal I}(h^{m}).
\] 
Then $\tilde{\cal I}_{D}(h^{m})$ is an ideal sheaf on $D$, since 
${\cal O}_{D}(mL)\otimes{\cal I}(h^{m})$ is a subsheaf of 
the locally free sheaf ${\cal O}_{D}(mL-v_{m}(D)D)$  on 
the smooth variety $D$.
We define the ideal sheaf $\sqrt[m]{\tilde{\cal I}_{D}(h^{m})}$ 
on $D$ by 
\[
\sqrt[m]{\tilde{\cal I}_{D}(h^{m})}_{x}
:= \cup {\cal I}(\frac{1}{m}(\sigma))_{x} (x\in D),
\]
where $\sigma$ runs all the germs of $\tilde{\cal I}_{D}(h^{m})_{x}$.
And we set 
\[
{\cal I}_{D}(h):= \cap_{m\geq 1}\sqrt[m]{\tilde{\cal I}_{D}(h^{m})}
\]
and call it {\bf the multipler ideal sheaf of $h$ on $D$}.
Also we set 
\[
\bar{\cal I}_{D}(h):= \lim_{\varepsilon\downarrow 0}
{\cal I}_{D}(h^{1+\varepsilon}).
\]
If $h\mid_{D}$ is well defined, then 
\[
\bar{\cal I}_{D}(h) = \bar{\cal I}(h\mid_{D})
\]
holds (\cite[Theorem 2.8]{tu4}). 
For every irreducible curve $C$ on $D$,
we say that the intersection number $(L,h)\cdot C$ is well defined, if 
$\nu (\Theta_{h},x) = \nu (\Theta_{h},D)$ holds for 
a very general point $x$ on $C$. 
In this case  
${\cal I}_{D}(h^{m})\mid_{C}$ is an ideal sheaf 
on $C$. 

 We define the 
{\bf intersection number} $(L,h)\cdot C$ by 
\[
(L,h)\cdot C := 
\overline{\lim}_{m\rightarrow\infty} 
m^{-1}\dim H^{0}(C,{\cal O}_{C}(mL-v_{m}(D)D)
\otimes\tilde{\cal I}_{D}(h^{m})/tor).
\]
Then we see that 
\[
(L,h)\cdot C =  (L- \nu (\Theta_{h},D)D)\cdot C + 
\overline{\lim}_{m\rightarrow\infty}
m^{-1}\deg \tilde{\cal I}_{D}(h^{m})
\]
holds. 

We can  also define the {\bf volume} of $r$-dimensional subvariety $Y$ of $D$ 
with respect to $(L,h)$ by using $\tilde{\cal I}_{D}(h^{m})$ as 
\[
\mu (Y,(L,h)) := r!\cdot\overline{\lim}_{m\rightarrow\infty}
m^{-r}\dim H^{0}(Y,{\cal O}_{Y}(mL-v_{m}(D)D)\otimes \tilde{\cal I}_{D}(h^{m})/tor). 
\]
But this coincides the definition before as is easily be seen.

We may define the {\bf Lelong number} $\nu_{D}(\Theta_{h},x) (x\in D)$ 
by 
\[
\nu_{D}(\Theta_{h},x) := \overline{\lim}_{m\rightarrow\infty}
m^{-1}\mbox{mult}_{x}\mbox{Spec}({\cal O}_{D}/\tilde{\cal I}_{D}(h^{m})),
\]
where $\mbox{mult}_{x}$ denotes the multiplicity on $D$. 
Then we see that the set 
\[
S_{D}:= \{ x\in D\mid \nu_{D} (\Theta_{h},x) > 0\}
\]
is at most countable union of subvarieties on $D$. 
This follows from the approximation theorem 
\cite[p.380,Proposition 3.7]{d}. 
Also this is obvious, 
if $h$ is an AZD constructed as in Section 3.2.
\subsection{Another definition of the intersection number}
Let $(L,h)$ be a singular hermitian line bundle on 
a smooth projective variety $X$ such that $\Theta_{h}$ 
is positive. 
And let $C$ be an irreducible curve on $X$ such that 
the restriction $h\mid_{C}$ is well defined. 
The another candidate of the intersection number of $(L,h)$ 
and $C$ is :
\[
(L,h)*C := L\cdot C - \sum_{x\in C}\nu (\Theta_{h}\mid_{C},x).
\]
But we have the following theorem.
\begin{theorem}(\cite[Theorem 2.7]{tu4}) 
\[
(L,h)\cdot C = (L,h)*C
\]
holds. 
\end{theorem}
\subsection{Limit multiplicities}

Let $(L,h)$ be a singular hermitian line bundle on a smooth projective 
variety $X$. 
Suppose that $\Theta_{h}$ is strictly positive. 
In this subsection, we shall consider the behavior of 
\[
\mbox{mult}_{x}\mbox{Bs}\mid m(L,h)\mid
\]
as $m$ goes to infinity. We shall prove the following theorem. 
\begin{theorem}
Let $(L,h)$ be a singular hermitian line bundle on 
a smooth projective variety $X$. 
Suppose that $\Theta_{h}$ is strictly positive. 
Let $x_{0}\in X$ be a point such that $\nu (\Theta_{h},x_{0}) = 0$. 
Let $c$ be a positive number such that 
\[
c < \mu (X,(L,h)).
\]
Then for every $x\in X$
\[
\nu (x) := \lim_{m\rightarrow\infty}
m^{-1}\mbox{mult}_{x}\mbox{Bs}\mid m(L,h)\otimes{\cal M}_{x_{0}}^{[cm]}\mid 
\]
exists, where 
$\mid m(L,h)\otimes{\cal M}_{x_{0}}^{[cm]}\mid$ 
denotes 
\[
\mid H^{0}(X,{\cal O}_{X}(mL)\otimes{\cal I}(h^{m})\otimes{\cal M}_{x_{0}}^{[cm]})\mid .
\]
Moreover for any modification 
\[
f : Y\longrightarrow X
\]
and $y\in Y$, 
\[
\nu (y):=  \lim_{m\rightarrow\infty}
m^{-1}\mbox{mult}_{y}f^{*}\mbox{Bs}\mid m(L,h)\otimes{\cal M}_{x_{0}}^{[cm]}\mid 
\]
exists.
\end{theorem}
{\bf Proof of Theorem 3.6.}
For $x\in X$ we set  
\[
\bar{\nu}(x): = \overline{\lim}_{m\rightarrow\infty}
m^{-1}\mbox{mult}_{x}\mbox{Bs}\mid m(L,h)\otimes{\cal M}_{x_{0}}^{[cm]}\mid.
\]
We claim that for any $\epsilon > 0$ and $x\in X$, there exists a positive integer 
$m(\epsilon )$  such that for every $m\geq m(\epsilon )$
\[
\mbox{mult}_{x}\mbox{Bs}\mid m(L,h)\otimes{\cal M}_{x_{0}}^{[cm]}\mid 
\geq (1-\epsilon )\bar{\nu}(x)m
\]
holds. 

Let $\delta$ be a small poitive number such that 
\[
\mu (X,(L,h)) > c+\delta 
\]
and $c+\delta$ is a rational number. 
Let us fix $x\in X$. 
Let 
\[
\pi : \tilde{X}\longrightarrow X
\] 
be the blowing up at $\{ x,x_{0}\}$.
We set 
\[
E := \pi^{-1}(x) 
\]
and
\[
E_{0} := \pi^{-1}(x_{0}). 
\]
We shall prove the following lemma. 
\begin{lemma}
There exists a singular hermitian metric $\tilde{h}_{\delta}$ on 
$\pi^{*}L$ such that 
\begin{enumerate}
\item 
$\pi_{*}\bar{\cal I}(h_{\delta}^{m}) 
\subseteq {\cal I}(h^{m})\otimes {\cal M}_{x_{0}}^{[(c+\frac{1}{2}\delta )m]}$ holds
for every sufficiently large $m$,
\item $\Theta_{h_{\delta}}$ is strictly positive on $\tilde{X}$.
\end{enumerate}
\end{lemma}
{\em Proof}.
Let $H$ be a very ample divisor on $\tilde{X}$ and let 
$h_{H}$ be a $C^{\infty}$-hermtian metric on $H$ with strictly positive 
curvature. 
Let $\varepsilon_{H}$ be a sufficiently small positive rational number 
such that 
\[
\mu (\tilde{X},(L-\varepsilon_{H}H,h\cdot h_{H}^{-\varepsilon_{H}})) > 0 
\]
holds. 
For every sufficiently large $\ell$, let  
\[
\tilde{\sigma}_{\ell}\in H^{0}(\tilde{X},
{\cal O}_{\tilde{X}}(\ell !(\pi^{*}L-\varepsilon_{H}H -(c+\delta )E_{0})\otimes
\pi^{*}{\cal I}(h^{\ell !}))
\]
be a nontrivial section and set 
\[
\tilde{h}_{\ell} = \frac{1}{\mid\tilde{\sigma}_{\ell}\mid^{2/\ell !}}.
\]
Since for every sufficiently large $\ell$, $\Theta_{h_{\ell}}$ is a closed positive current which represents
$2\pi c_{1}(\pi^{*}L-\varepsilon_{H}H)$, we may assume that there exists a subsequence $\{ \Theta_{h_{\ell_{j}}}\}$ 
of $\{ \Theta_{h_{\ell}}\}$ such that 
\[
\Theta_{\infty}:= \lim_{j\rightarrow\infty}\Theta_{\tilde{h}_{\ell_{j}}}
\]
exists as a closed positive current. 
Let $\tilde{h}$ be a singular hermitian metric on $\pi^{*}L-\varepsilon_{H}H$ such that 
\[
\Theta_{\tilde{h}} = \Theta_{\infty}
\]
holds. 
Then 
\[
\tilde{h}_{\delta}:= \tilde{h}\cdot h_{H}^{\varepsilon_{H}}
\]
is a singular hermitian metric on $\pi^{*}L$ with strictly positive curvature 
current. 
By the construction we see that 
for every modification 
\[
f : Y\longrightarrow X
\]
and $y\in Y$, 
\[
\nu (f^{*}\Theta_{h_{\delta}},y)
\geq \nu (f^{*}\Theta_{h},y) + (c+\delta)\cdot\mbox{mult}_{y}f^{*}E_{0}
\]
holds. 
Hence by Lemma 2.1 and Lemma 2.2 we see that 
\[
\pi_{*}\bar{\cal I}(h_{\delta}^{m}) 
\subseteq {\cal I}(h^{m})\otimes {\cal M}_{x_{0}}^{[(c+\frac{1}{2}\delta )m]}
\]
holds
for every sufficiently large $m$.
This completes the proof of Lemma 3.2.
{\bf Q.E.D.} \vspace{5mm} \\

Suppose there exists a point $x\in X$ such that 
for some $\epsilon > 0$, there exists an 
increasing  sequence of positive integers $\{ m_{j}\}$ 
such that 
\[
\mbox{mult}_{x}\mbox{Bs}\mid m_{j}(L,h)\otimes{\cal M}_{x_{0}}^{[cm_{j}]}\mid 
< (1-\epsilon )\bar{\nu}(x)m_{j}
\]
holds. 
Let 
\[
\sigma_{j}\in H^{0}(X,{\cal O}_{X}(m_{j}L)\otimes{\cal I}(h_{j}^{m})\otimes
{\cal M}^{\lceil cm_{j}\rceil})
\]
(here for a real number $a$, $\lceil a \rceil$ denotes the smallest integer 
larger or equal to $a$) be a nonzero element such that 
\[
\mbox{mult}_{x}(\sigma_{j} )\leq (1-\epsilon )\bar{\nu}(x)m_{j}
\]
holds.
We define the singular hermitian metric $h_{j}$ of $L$ by
\[
h_{j} :=\frac{1}{\mid\sigma_{j}\mid^{2/m_{j}}}. 
\]
Let $\tilde{x}\in \tilde{X}$ be a point on $E$ such that 
for every $m$, ${\cal I}(\pi^{*}h^{m})$ is locally free on 
a neighbourhood of $\tilde{x}$ (the neighbourhood may depend on $m$). 

Let $U$ be a small neighbourhood of $\tilde{x}$ 
and let $\rho$ be a $C^{\infty}$-function on $\tilde{X}$ 
such that 
\begin{enumerate}
\item $\mbox{Supp}\,\rho \subset\subset U$,
\item $\rho\equiv 1$ on a neighbourhood of $\tilde{x}$,
\item $0 \leq \rho \leq 1$
\end{enumerate}
hold. 
Let $d_{\tilde{x}}$ denote the distance function from $\tilde{x}$ 
with respect to a fixed K\"{a}hler form  $\omega$ on $\tilde{X}$. 

Let $\nu_{0}$ be a sufficiently large positive integer 
such that 
\[
\nu_{0}\Theta_{\tilde{h}_{\delta}}+\mbox{Ric}_{\omega} 
+ 2n\sqrt{-1}\partial\bar{\partial}(\rho \log d_{\tilde{x}})
\]
is strictly positive and 
\[
\pi_{*}{\cal I}(h_{\delta}^{\nu_{0}}) 
\subseteq {\cal I}(h^{\nu_{0}})\otimes {\cal M}_{x_{0}}^{[(c+\frac{1}{2}\delta )\nu_{0}]}
\]
holds. 
Let $h_{0}$ be a $C^{\infty}$-hermitian metric on $L$.  
Let $\epsilon_{0}$ be a sufficiently small positive number 
such that 
\[
\nu_{0}\Theta_{\tilde{h}_{\delta}}+\mbox{Ric}_{\omega} 
+ 2n\sqrt{-1}\partial\bar{\partial}(\rho \log d_{\tilde{x}})
- \epsilon_{0}\pi^{*}\Theta_{h_{0}}
\]
is strictly positive. 
Then by Nadel's vanishing theorem (Theorem 2.1), 
\[
H^{1}(\tilde{X},{\cal O}_{\tilde{X}}(\pi^{*}L)
\otimes {\cal I}(\pi^{*}h_{j}^{m+\epsilon_{0}}\cdot
\tilde{h}_{\delta}^{\nu_{0}}e^{-2n\rho\log d_{\tilde{x}}}))
= 0
\]
holds for every $j$ and $m\geq 0$. 
This implies that there exists
\[
\tilde{\sigma}\in H^{0}(\tilde{X},{\cal O}_{\tilde{X}}(\pi^{*}(m+\nu_{0})L)
\otimes {\cal I}(\pi^{*}h_{j}^{m+\epsilon_{0}}\cdot\tilde{h}_{\delta}^{\nu_{0}}))
\]
which generates 
\[
{\cal O}_{\tilde{X}}(\pi^{*}(m+\nu_{0})L)
\otimes {\cal I}(\pi^{*}h_{j}^{m+\epsilon_{0}}\cdot\tilde{h}_{\delta}^{\nu_{0}})
\]
at $\tilde{x}$. 
Let us fix an arbitrary positive integer $m$. 
We note that for every sufficiently large $j$, 
\[
\pi_{*}({\cal I}(\pi^{*}h_{j}^{m+\epsilon_{0}}\cdot\tilde{h}_{\delta}^{\nu_{0}}))
\subseteq {\cal I}(h^{m+\nu_{0}})\otimes{\cal M}_{x_{0}}^{[c(m+\nu_{0})]}
\]
holds by the construction of $\{ h_{j}\}$. 
Hence we see that 
\[
\mbox{mult}_{\tilde{x}}\mbox{Bs}\mid (m+\nu_{0})(L,h)\otimes {\cal M}_{x_{0}}^{[cm]}\mid
<  (1-\epsilon )(m+\nu_{0})\bar{\nu}(x)
\]
holds. 
Since $m$ is an arbitrary positive integer, this is the contradiction.
Hence we conclude that for every $\epsilon > 0$, 
there exists a positive integer $\nu_{0}$ such that
for every $m\geq \nu_{0}$
\[
\mbox{mult}_{\tilde{x}}\mbox{Bs}\mid m(L,h)\otimes {\cal M}_{x_{0}}^{[cm]}\mid
\geq   (1-\epsilon )m\bar{\nu}(x)
\]
holds. 
This implies that 
\[
\underline{\lim}_{m\rightarrow\infty}
m^{-1}\mbox{mult}_{x}\mbox{Bs}\mid m(L,h)\otimes{\cal M}^{[cm]}\mid
= \bar{\nu}(x)
\]
holds. 
By the definition of $\bar{\nu}(x)$ we see that 
\[
\nu (x) := \lim_{m\rightarrow\infty}
m^{-1}\mbox{mult}_{x}\mbox{Bs}\mid m(L,h)\otimes{\cal M}^{[cm]}\mid
\]
exists.
Since $x$ is arbitrary, this completes the proof of Theorem 3.6
except the last statement.
The proof of the last statement is similar.
{\bf Q.E.D.}
\section{Fibration theorem}
\subsection{The first nonvanishing theorem}
Let $X$ be a smooth projective variety. 
Let $D$ be a divisor on $X$ and let $A = \sum a_{i}A_{i}$ be a {\bf  Q}-divisor 
on $X$. 
Assume the following conditions :
\begin{enumerate}
\item $D$ is nef,
\item $\lceil A\rceil:= \sum_{i}\lceil a_{i}\rceil A_{i}$ is effective, 
\item $\mbox{Supp}\{ A\}$ is a divisor with normal crossings, 
where $\{ A\}$ denotes the fractional part of $A$, i.e., 
$\{ A\} := A - [A]$. 
\item there exists a positive integer $a$ such that 
$aD+A-K_{X}$ is nef and big. 
\end{enumerate}
In \cite{sh}, Shokurov proved that under these 
conditions, 
for every sufficiently large positive integer $b$,
\[
H^{0}(X,{\cal O}_{X}(bD+\lceil A\rceil))\neq 0
\]
holds.
In this section, we shall prove a similar nonvanishing theorem  which 
plays  essential roles in this paper.  
\begin{theorem}(The first nonvanishing theorem)
Let $X$ be a smooth projective variety and let $(L,h_{L})$ be a 
singular hermitian line bundle on $X$ such that the curvature current
$\Theta_{L}$ is positive.
Let $(A,h_{A})$ be a singular hermitian line bundle on $X$
with strictly positive curvature current $\Theta_{A}$.
Then one of the followings holds.
\begin{enumerate}
\item 
$H^{0}(X,{\cal O}_{X}(K_{X}+A+mL)\otimes {\cal I}(h_{A}h_{L}^{m}))\neq 0
$ 
holds for every sufficiently large $m$,
\item there exists
a nontrivial numerically trivial fiber space structure  
\[
f : X-\cdots\rightarrow Y,
\] 
i.e.,
\begin{enumerate}
\item $f$ is regular over the generic point of $Y$, 
\item for a very  general fiber $F$ the restriction 
$(L,h_{L})\mid_{F}$ is numerically trivial, 
\item for a very general point $x\in X$ and  every 
irreducible horizontal curve $C$ containing $x$, $(L,h)\cdot C > 0$
holds,
\item $\dim Y < \dim X$ 
is minimal among such fibrations.
\end{enumerate}
\end{enumerate}  
\end{theorem}
Let us compare Theorem 4.1 with Shokurov's nonvanishing 
theorem above. 
The positivity of $\Theta_{L}$ in Theorem 4.1 corresponds 
to the nefness of $D$ in Shokurov's theorem. 
The strict positivity of $(A,h_{A})$ corresponds to 
the the third condition in Shokurov's theorem. 
But the second condition in Shokurov's theorem does not 
have a counterpart in Theorem 4.1. 
That is why we have two cases. 
Roughly speaking Theorem 4.1 tells us what happens, if 
we drop the second condition in Shokurov's theorem. 
To construct a nontrivial holomorphic  section 
of ${\cal O}_{X}(K_{X}+A+mL)$ on $X$ in the second case,
we need to construct a section on a general fiber $F$
of the numerically trivial fibration 
$f : X-\cdots\rightarrow  Y$. 
This problem will be treated in the second nonvanishing
theorem (Theorem 4.4) later.

The following (more algebraic) corollary follows 
from the proof of Theorem 4.1(Corollary 4.1 is a corollary of the proof of Theorem 4.1.  See Remark 4.1 below.). 
\begin{corollary}(\cite[Corollary 8.1]{tu3})
Let $X$ be a smooth projective variety and let $L$ be a nef 
line bundle on $X$.
Let $A$ be a big line bundle on $X$.
Then one of the followings holds.
\begin{enumerate}
\item \[
H^{0}(X,{\cal O}_{X}(K_{X}+A+mL))\neq 0
\]
holds for every sufficiently large $m$,
\item there exists a rational fibration 
\[
f : X-\cdots\rightarrow Y
\]
such that 
\begin{enumerate} 
\item $f$ is regular over the generic point of $Y$,
\item  for a very general fiber $F$ the restriction
$L_{\mid F}$ is numerically trivial,
\item for every general point $x$ on $X$ and every 
irreducible horizontal (with respect to $f$) curve 
$C$ containing $x$, $L\cdot C > 0$ holds,
\item $\dim Y < \dim X$ is minimal among such fibrations. 
\end{enumerate}    
\end{enumerate}  
\end{corollary}
\begin{remark}
In Corollary 4.1 $L$ may not admit a singular hermitian metric $h$ such that 
$\Theta_{h}$ is positive and ${\cal I}(h^{m}) = {\cal O}_{X}$ for 
every $m\geq 0$ on $X$. 
But the proof is parallel to that of Theorem 4.1, if we change the volume
$\mu (X,(A+mL,h_{A}h_{L}^{m}))$ of 
a subvariety $V$ in $X$ with respect to $A+mL$ (see Lemma 4.1 below) by the intersection number 
$(A+mL)^{\dim V}\cdot V$.
\end{remark}

\begin{example}
To illustrate our method let us consider the following example.
Let $X$ be an irreducible quotient of the open unit bidisk $\Delta^{2}$ 
in $\mbox{\bf  C}^{2}$, i.e.
\[
X = \Delta^{2}/\Gamma ,
\]
where $\Gamma$ is an irreducible cocompact torsion free lattice. 
Let $L$ denotes the line bundle such that whose curvature form 
comes from the Poincar\'{e} metric on the first factor. 
Then one see that $L$ is nef and $L^{2} = 0$ holds.
In particular $L$ is not big.
Let $c_{1}(L)$ be the first Chern form of $L$ induced by 
the Poincar\'{e} metric on the first factor. 
On the other hand for every ample line bundle $A$, 
$K_{X}+ mL + A$ is very ample for $m>> 1$.
Moreover for every singular hermitian line bundle 
$(A,h_{A})$ with strictly positive curvature (in the sense of 
current),
\[
H^{0}(X,{\cal O}_{X}(K_{X}+mL+ A)\otimes{\cal I}(h_{A})) 
\]
gives a birational rational map from $X$ into a projective space
and even it separates jets of any fixed order $k$ at very general points 
on $X$ for every sufficiently large $m$ (of course such  $m$ depends
on $k$).
In a sense $L$ behaves more or less like an ample line bundle.

\end{example}

\subsection{Numerically trivial fibrations}
The following theorems are key ingredients for our proof 
of Theorem 4.1 and Theorem 1.1.  
\begin{theorem}(\cite[Theorem 1.1]{tu4})
Let $(L,h)$ be a singular hermitian line bundle on 
a smooth projective variety $X$. 
Suppose that the curvature current $\Theta_{h}$ is positive. 
Then there exists a unique (up to birational equivalence) rational fibration 
\[
f : X -\cdots \rightarrow Y
\]
such that 
\begin{enumerate}
\item $f$ is regular over the generic point of $Y$, 
\item for every very general fiber $F$, 
 $(L,h)\mid_{F}$ is well defined and is numerically trivial (cf. Definition 3.9),
\item $\dim Y$ is minimal among such fibrations,
\item for a very general point $x\in X$ and any irreducible  
 horizontal  curve (with respect to $f$) $C$ containing $x$,
$(L,h)\cdot C > 0$ holds.
\end{enumerate}
\end{theorem}
We call the above fibration {\bf the numerically trivial fibration}
associated with $(L,h)$.
\begin{remark}
Let $X$,$(L,h)$ be as above. Then for any smooth divisor $D$ on $X$,
there exists a numerically trivial fibration 
\[
f_{D} : D -\cdots \rightarrow W.
\]
This is simply because the restriction of the intersection theory 
on $D$ exists and the proof of the above theorem essentially
does not require the existence of the restriction of $\Theta_{h}$
on $D$. 
\end{remark} 

The structure of numerically trivial singular hermitian line bundles 
with positive curvature current is given  as follows.
\begin{theorem}(\cite[Theorem 1.2]{tu4})
Let $(L,h)$ be a singular hermitian line bundle on a 
smooth projective variety $X$. 
Suppose that $(L,h)$ is numerically trivial on $X$. 
Then there exist at most countably many prime divisors $\{ D_{i}\}$
and nonnegative numbers $\{ a_{i}\}$ such that 
\[
\Theta_{h} = 2\pi\sum_{i}a_{i}D_{i}
\]
holds. 
More generally  let $Y$  be a subvariety of $X$ such that the restriction  
$h\mid_{Y}$ is well defined. 
Suppose that  $(L,h)$ is numerically trivial 
on $Y$. 
Then the restriction $\Theta_{h}\mid_{Y}$ is  a sum of
at most  countably many 
prime divisors with nonnegative coefficients on $Y$. 
\end{theorem}
Theorem 4.3 gives an information on the restricion of 
the singular hermitian metric on a very general fiber 
of a numerically trivial fibration.  
\begin{corollary}(\cite[Corollary 3.2]{tu4})
Let $X$ be a smooth projective variety and let 
$(L,h)$ be a singular hermitian line bundle on $X$
such that $\Theta_{h}$ is positive. 
Let $D$ be a smooth divisor on $X$. 
Suppose that $(L,h)$ is numerically trivial on $D$. 
Then 
\[
S_{D}:= \{ x\in D\mid \nu_{D}(\Theta_{h},x) > 0\}
\]
is a sum of at most countably many prime divisors on $D$, 
where $\nu_{D}(\Theta_{h},x)$ is the Lelong number 
defined as in Section 3.6.
Also 
\[
(L - \nu (\Theta_{h},D)\cdot D)_{D} - \sum_{E}\nu_{D}(\Theta_{h},E)\cdot E
\]
is numerically trivial on $D$, where $E$ runs all the prime divisors on $D$. 
\end{corollary}
\begin{remark}
Corollary 4.2 still holds for a subvariety $V$ on 
$D$, if there exists a curve on $V$ such that 
$(L,h)\cdot C$ is well defined (cf. \cite[Remark 3.1]{tu4}). 
\end{remark}
  
\subsection{Proof of Theorem 4.1}

Let $X,(L,h_{L}), (A,h_{A})$ be as in Theorem 4.1. 
Let 
\[
f : X -\cdots\rightarrow Y
\]
be the numerically fibration associated with $(L,h_{L})$. 
If $\dim Y < \dim X$ holds, then this is the desired fibration.
Hence we shall assume that $f$ is the identity morphism.
In other words, {\bf for a very general point $x$ on $X$ 
and  
any irreducible curve $C$ containing $x$, 
$h_{L}\mid_{C}$ is well defined and  
\[
(L,h_{L})\cdot C > 0
\] 
holds.} 
We say that $(L,h_{L})$ is {\bf very generically numerically 
positive}. 
\begin{lemma}(\cite[Lemma 4.1]{tu4})
Suppose that $(L,h_{L})$ is not numerically trivial.
Then for every ample line bundle $H$ on $X$ 
\[
\overline{\lim}_{m\rightarrow\infty}m^{-1}\mu (X,(H+mL,h_{H}h^{m}_{L})) > 0
\]
holds, where $h_{H}$ is any $C^{\infty}$ hermitian metric 
with strictly positive curvature on $H$. 
\end{lemma}
{\bf Proof}.
Let $n$ be the dimension of $X$.
We prove this lemma by induction on $n$.
If $n =1$,  Lemma 4.1 is trivial.
Let $\pi :\tilde{X} \longrightarrow {\bf  P}^{1}$
be a Lefschetz pencil associated with a very ample linear 
system  say $\mid H\mid$ on $X$.
If we take the pencil very general, we may assume that 
${\cal I}(h_{L}^{\ell})$ is an ideal sheaf on 
all fibers of $\pi$ for every $\ell \geq 1$.
And let 
\[ 
b : \tilde{X}\longrightarrow X
\]
be the modification associated with the pencil and let
$E$ be the exceptional divisor of $b$.
Then by the inductive assumption for a very general fiber $F$ of
$\pi$, we see that 
\[
\overline{\lim}_{m\rightarrow\infty}m^{-1}\mu (F,b^{*}(H+mL,h_{H}h_{L}^{m})) > 0
\]
holds.
Let us consider the direct image
\[
{\cal E}_{m,\ell} := \pi_{*}{\cal O}_{\tilde{X}}(\ell b^{*}(H+mL))\otimes
{\cal I}(b^{*}h^{m\ell}_{L})).
\]
By Grothendiek's theorem, we see that
\[
{\cal E}_{m,\ell} \simeq \oplus_{i=1}^{r}{\cal O}_{{\bf  P}^{1}}(a_{i})
\]
for some $a_{i} = a_{i}(m,\ell )$ and $r = r(m,\ell )$.
By the inductive assumption, we see that 
\[
\overline{\lim}_{m\rightarrow\infty}m^{-1}(\overline{\lim}_{\ell\rightarrow\infty}
\ell^{-(n-1)}r(m,\ell)) > 0
\]
holds.
We note that $\ell_{0}b^{*}H - E$ is ample
 for some positive integer $\ell_{0}$.
Hence we see that 
\[
{\cal O}_{\tilde{X}}(\ell_{0}b^{*}H- E)
\]
admits a $C^{\infty}$-hemitian metric with  strictly positive curvature.
Hence by Nadel's vanishing theorem \cite[p.561]{n} there is a positive constant 
$c$ such that
\[
H^{1}(\tilde{X},{\cal O}_{\tilde{X}}(\ell b^{*}(H+mL-\frac{1}{\ell_{0}}E))
\otimes{\cal I}(b^{*}h^{m\ell}_{L})
\otimes \pi^{*}{\cal O}_{{\bf  P}^{1}}(-\lceil c\ell\rceil )) = 0
\]
holds for every sufficiently large $\ell$ divisible by $\ell_{0}$.
This implies that 
\[
\overline{\lim}_{\ell\rightarrow\infty}\ell^{-1}(\min_{i}a_{i})\geq  c
\]
holds  and 
\[
\overline{\lim}_{\ell\rightarrow\infty}
\ell^{-n}\dim H^{0}(\tilde{X},{\cal O}_{\tilde{X}}(\ell p^{*}(H+mL))
\otimes{\cal I}(p^{*}(h^{m\ell}_{L}))) \geq 
\]
\[
\,\,\,\,\,\,\,\,\,\,\,\,\,\geq c\cdot\overline{\lim}_{\ell\rightarrow\infty}
\ell^{-(n-1)}r(m,\ell)
\] 
holds.
Hence we see that 
\[
\overline{\lim}_{m\rightarrow\infty}m^{-1}(\overline{\lim}_{\ell\rightarrow\infty}
\ell^{-n}\dim H^{0}(\tilde{X},{\cal O}_{\tilde{X}}(\ell b^{*}(H+mL))
\otimes{\cal I}(b^{*}(h^{m\ell}_{L}))) > 0
\]
holds.
Since 
\[
b_{*}{\cal I}(b^{*}h^{m\ell}_{L}) \subseteq
{\cal I}((h^{m\ell}_{L}))
\]
holds by Lemma 2.1,
we see that
\[
\overline{\lim}_{m\rightarrow\infty}m^{-1}\mu (X,(H+mL,h_{H}h^{m}_{L})) > 0
\]
holds.
Here we have assumed that $H$ to be sufficiently very ample. 
To prove the general case of Lemma 4.1, we argue as follows.
Let $H$ be any ample line bundle on $X$.
Then 
\[
\mu (X,(a(H+mL),h_{H}^{a}h^{am}_{L})) = a^{n}\cdot \mu (X,(H+mL,h_{H}h^{m}_{L}))
\]
holds for every positive integer $a$. 
Now it is clear that Lemma 4.1 holds for any ample line bundle $H$. 
This completes the proof of Lemma 4.1. {\bf  Q.E.D.} \vspace{5mm} \\

Let $H$ be a very ample line bundle on $X$ and let $h_{H}$ be 
a $C^{\infty}$-hermitian metric on $H$. 
Since the curvature $\Theta_{A}$ of $h_{A}$ is strictly positive, we see that 
for every sufficiently small positive number $\epsilon$. 
\[
\mu (X,(A-\epsilon H,h_{A}h_{H}^{-\epsilon})) > 0
\]
holds by Theorem 2.1. 

By the assumption $(L,h_{L})$ is very generically numerically positive on $X$. 
Let $\nu_{0}$ be a positive integer and let us consider the singular hermitian line bundle
$(\nu_{0}L+A, h_{L}^{\nu_{0}}h_{A})$ on $X$.
By Lemma 4.1  we can take a sufficiently
large $\nu_{0}$ so that 
\[
\mu (X,(\nu_{0}L+A,h_{L}^{\nu_{0}}h_{A}))   > 
\epsilon^{n}\mu (X,(\nu_{0}L+\epsilon H,h_{L}^{\nu_{0}}h_{H}^{\epsilon})) >  2^{n}(n+1)^{2n} 
\]
hold. 
\begin{lemma}
Let $x\in X$ be a very general point such that $\nu(\Theta_{L},x) = 
\nu (\Theta_{A},x) = 0$ hold. 

Then for every sufficiently large positive integer $m$,
\[
H^{0}(X,{\cal O}_{X}(m(\nu_{0}L+A)\otimes{\cal I}((h_{A}h_{L}^{\nu_{0}})^
{m})\otimes {\cal M}_{x}^{\otimes 2(n+1)^{2}m})\neq 0
\]
holds, 
where ${\cal M}_{x}$ denotes the maximal ideal sheaf at $x$. 
\end{lemma}
{\bf Proof}. 
Let us consider the following morphism
\begin{eqnarray*}
H^{0}(X,{\cal O}_{X}(m(\nu_{0}L+A)\otimes{\cal I}((h_{A}h_{L}^{\nu_{0}})^
{m}))
\rightarrow  \\
H^{0}(X,{\cal O}_{X}(m(\nu_{0}L+A)\otimes{\cal I}((h_{A}h_{L}^{\nu_{0}})^
{m})/{\cal M}_{x}^{\otimes 2(n+1)^{2}m})
\end{eqnarray*}
The kernel of this morphism is exactly
\[
H^{0}(X,{\cal O}_{X}(m(\nu_{0}L+A)\otimes{\cal I}((h_{A}h_{L}^{\nu_{0}})^
{m})\otimes {\cal M}_{x}^{\otimes 2(n+1)^{2}m}).
\]
If we take $x$ very general we may assume that $\nu (\Theta_{L},x) = 0$ 
holds. 
Hence ${\cal I}(h_{L}^{m})_{x} = 
 {\cal O}_{X,x}$ holds for every $m\geq 0$
by Lemma 2.2.
Since 
\[
\dim H^{0}(X,{\cal O}_{X}(m(\nu_{0}L+A)\otimes{\cal I}((h_{A}h_{L}^{\nu_{0}})^{m})/{\cal M}_{x}^{\otimes 2(n+1)^{2}m})
=  \frac{2^{n}(n+1)^{2n}}{n!}m^{n} + O(m^{n-1})
\]
and 
\[
\mu (X,(\nu_{0}L+A,h_{L}^{\nu_{0}}h_{A})) > 2^{n}(n+1)^{2n}
\]
hold, for every sufficiently large $m$,
\begin{eqnarray*}
\dim H^{0}(X,{\cal O}_{X}(m(\nu_{0}L+A)\otimes{\cal I}((h_{A}h_{L}^{\nu_{0}})^
{m}))
>   \\
\dim H^{0}(X,{\cal O}_{X}(m(\nu_{0}L+A)\otimes{\cal I}((h_{A}h_{L}^{\nu_{0}})^
{m})/{\cal M}_{x}^{\otimes 2(n+1)^{2}m})
\end{eqnarray*}
holds.
This completes the proof of Lemma 4.2.\vspace{5mm} {\bf  Q.E.D.}

Let $x$ be a very general point of $X$ such that 
$\nu (\Theta_{L},x) = \nu (\Theta_{A},x) = 0$ hold and 
for every irreducible curve $C$ containing $x$ 
$h_{L}\mid_{C}$ is well defined and  satisfies 
\[
(L,h_{L})\cdot C > 0.
\]
Let  $\varepsilon$ be a sufficiently small positive number.
Let $m_{0}$ be a sufficiently large positive number and let
\[
\sigma_{0}\in
H^{0}(X,{\cal O}_{X}(m_{0}(\nu_{0}L+A)\otimes{\cal I}((h_{A}h_{L}^{\nu_{0}})^
{m_{0}})\otimes {\cal M}_{x}^{\otimes 2(n+1)^{2}m_{0}})
\]
be a general nonzero section.

Let $h_{0}$ be the singular hermiaitn metric on $\nu_{0}L+A$ 
defined by
\[
h_{0} = \frac{1}{\mid\sigma_{0}\mid^{2/m_{0}}}.
\]
Let $\alpha_{0}$ be the positive number defined by
\[
\alpha_{0} := \inf \{\alpha \mid ({\cal O}_{X}/{\cal I}(h_{0}^{\alpha}))_{x}\neq 0\} .
\]
We set  
\[
V_{1} := \lim_{\delta\downarrow 0}\mbox{Spec} (\cal O_{X}/{\cal I}(\alpha_{0}+\delta )).
\]
And let $X_{1}$ be a branch of $V_{1}$ containing $x$.
Then since $\sigma_{0}$ is an element of 
\[
H^{0}(X,{\cal O}_{X}(m_{0}(\nu_{0}L+A)\otimes{\cal I}((h_{A}h_{L}^{\nu_{0}})^
{m_{0}})\otimes {\cal M}_{x}^{\otimes 2(n+1)^{2}m_{0}}),
\]
we see that
\[
\alpha_{0} < \frac{1}{2n}
\]
holds.
Let us take $m_{0}$ sufficiently large and 
$\sigma_{0}$ very general. 
\begin{remark}
By Theorem 3.6 for every $m_{0}^{\prime} \geq m_{0}$ 
and very general 
\[
\sigma_{0}^{\prime}\in
H^{0}(X,{\cal O}_{X}(m_{0}^{\prime}(\nu_{0}L+A)\otimes{\cal I}((h_{A}h_{L}^{\nu_{0}})^
{m_{0}^{\prime}})\otimes {\cal M}_{x}^{\otimes 2(n+1)^{2}m_{0}^{\prime}})
\]
we see that 
\[
X_{1} \subseteq \mbox{Spec}({\cal O}_{X}/{\cal I}((h_{0}^{\prime})^{\alpha_{0}+\delta_{0}}))
\]
holds, where $\delta_{0}$ is a positive number which 
tends to $0$ as $m_{0}$ tends to infinity. 
Hence even if we move $m_{0}$, only finitely many subvarieties
appears as  $X_{1}$ as far as we take $\sigma_{0}$ very general. 
\end{remark}
We set 
\[
n_{1} := \dim X_{1}.
\]
We note that $h_{L}\mid_{X_{1}}$ is well defined and 
$(L,h_{L})$ is not numerically trivial on $X_{1}$ 
by the choice of $x$.  
In this case by Lemma 4.1 we take a sufficiently large positive integer $\nu_{1}$ so that
\[
\mu (X_{1},(\nu_{1}L+A,h_{L}^{\nu_{1}}h_{A})) > 2^{n_{1}}(n_{1}+ 1)^{n_{1}}n^{n_{1}}.
\]
Then we have the following lemma.
\begin{lemma}
Let $x_{1}$ be a very general point on $X_{1,reg}$. 
Then for  every sufficiently large positive integer $m$,
\[
H^{0}(X_{1},({\cal O}_{X_{1}}(m(\nu_{0}L+A)\otimes{\cal I}((h_{A}h_{L}^{\nu_{1}})^{m}))/tor\otimes {\cal M}_{x_{1}}^{\otimes 2(n+1)^{2}m})\neq 0
\]
holds.
\end{lemma}
{\bf Proof}. 
Let $x_{1}$ be a very general point 
 on $X_{1,reg}$ such that
 $\nu (\Theta_{h_{A}h_{L}^{\nu_{1}}},x_{1}) = 0$ holds. 
Then for every $m$, ${\cal I}((h_{A}h_{L}^{\nu_{1}})^{m})_{x_{1}} = {\cal O}_{X,x_{1}}$ holds. 
Then the proof of Lemma 4.3 is parallel to that of Lemma 4.2. {\bf  Q.E.D.} 
\vspace{10mm} \\ 
Let $E$ be an effective {\bf  Q}-divisor such that 
$A+\nu_{1}L-E$ is ample (such a divisor exists by 
Kodaira's lemma \cite[Appendix]{k-o}). 
We set 
\[
H_{1} = r(A+\nu_{1}L- E),
\]
where $r$ is a positive integer such that $H$ is an integral divisor 
on $X$.
Then by Nadel's vanishing theorem, we have the following lemma.
\begin{lemma}
If we take $r$ sufficiently large, then 
\[
\phi_{m} : H^{0}(X,{\cal O}_{X}(m(\nu_{1}L+A)+H)\otimes{\cal I}((h_{A}h_{L}^{\nu_{0}})^{m}))\rightarrow 
\]
\[
H^{0}(X_{1},{\cal O}_{X_{1}}(m(\nu_{1}L+A)+H_{1})\otimes{\cal I}((h_{A}h_{L}^{\nu_{0}})^{m}))
\]
is surjective for every $m\geq 0$.
\end{lemma}
{\bf Proof}.
Let us take a locally free resolution of the ideal sheaf ${\cal I}_{X_{1}}$
of $X_{1}$.
\[
0\leftarrow {\cal I}_{X_{1}}\leftarrow {\cal E}_{1}\leftarrow {\cal E}_{2}
\leftarrow \cdots \leftarrow {\cal E}_{\ell}\leftarrow 0.
\]
Then by the trivial extension of Nadel's vanishing theorem to 
the case of vector bundles, if $r$ is sufficiently large, we have : 
\begin{sublemma}
\[
H^{q}(X,{\cal O}_{X}(m(\nu_{1}L+A)+H_{1})
\otimes{\cal I}((h_{A}h_{L}^{\nu_{0}})^{m})
\otimes{\cal E}_{j}) = 0
\]
holds for every $m\geq 1$, $q\geq 1$ and  $1\leq j\leq k$.
\end{sublemma}
In fact if we take $r$ sufficiently large, we see that for every $j$, 
${\cal O}_{X}(H - K_{X})\otimes {\cal E}_{j}$ admits a $C^{\infty}$-hermitian metric $g_{j}$ such that
\[
\Theta_{g_{j}} \geq \mbox{Id}_{E_{j}}\otimes \omega
\]
holds, where $\omega$ is a K\"{a}hler form on $X$.
By   \cite[Theorem 4.1.2 and Lemma 4.2.2]{ca}, we completes the proof of 
Sublemma 4.1.  {\bf  Q.E.D.}\vspace{5mm} \\ 
Let 
\[ 
p_{m} : Y_{m}\longrightarrow X
\]
be  a composition of successive blowing ups with smooth centers such that 
$p_{m}^{*}{\cal I}((h_{A}h_{L}^{\nu_{0}})^{m})$ is locally free on $Y_{m}$. 
\begin{sublemma}
\[
R^{p}p_{m\,*}p_{m}^{*}({\cal O}_{Y_{m}}(K_{Y_{m}})\otimes {\cal I}(p_{m}^{*}(h_{A}h_{L}^{m})) = 0
\]
holds for every $p\geq 1$ and $m\geq 1$.
\end{sublemma}
{\bf Proof}. 
This sublemma follows from Theorem 2.1.  
\vspace{10mm}
{\bf  Q.E.D.}  \\
We note that by the definition of the multiplier ideal sheaves 
\[
p_{m\, *}({\cal O}_{Y_{m}}(K_{Y_{m}})\otimes{\cal I}(p_{m}^{*}(h_{A}h_{L}^{\nu_{0}})^{m}))= {\cal O}_{X}(K_{X})\otimes {\cal I}((h_{A}h_{L}^{\nu_{0}})^{m})
\]
holds. 
Hence by Sublemma 4.1, Sublemma 4.2 and the Leray spectral sequence,  we see that 
\[
H^{q}(Y_{m},{\cal O}_{Y_{m}}(K_{Y_{m}}+p_{m}^{*}(m(\nu_{1}L+A)+H_{1} - K_{X}))
\otimes {\cal I}(p_{m}^{*}(h_{A}h_{L}^{\nu_{0}})^{m})
\otimes p_{m}^{*}{\cal E}_{j}) = 0
\]
holds for every $q\geq 1$ and $m\geq 1$.
Hence 
\[
H^{1}(Y_{m},{\cal O}_{Y_{m}}(K_{Y_{m}}+ p_{m}^{*}(m(\nu_{1}L+A)+H_{1}-K_{X})\otimes p_{m}^{*}{\cal I}((h_{A}h_{L}^{\nu_{0}})^{m}))\otimes p_{m}^{*}{\cal I}_{X_{1}}) = 0
\]
holds. 
Hence every element of 
\[
H^{0}(Y_{m},{\cal O}_{Y_{m}}(K_{Y_{m}}+ p_{m}^{*}(m(\nu_{1}L+A)+H_{1}-K_{X})\otimes {\cal I}(p_{m}^{*}(h_{A}h_{L}^{\nu_{0}})^{m}))\otimes {\cal O}_{Y_{m}}/p_{m}^{*}{\cal I}_{X_{1}}) 
\]
extends to an element of 
\[
H^{0}(Y_{m},{\cal O}_{Y_{m}}(K_{Y_{m}}+ p_{m}^{*}(m(\nu_{1}L+A)+H_{1}-K_{X})\otimes {\cal I}(p_{m}^{*}(h_{A}h_{L}^{\nu_{0}})^{m}))) 
\]
Also there exists a natural map 
\[
H^{0}(X_{1},{\cal O}_{X_{1}}(m(\nu_{1}L+A)+H_{1})\otimes{\cal I}((h_{A}h_{L}^{\nu_{0}})^{m}))
\rightarrow
\]
\[
H^{0}(Y_{m},{\cal O}_{Y_{m}}(K_{Y_{m}}+ p_{m}^{*}(m(\nu_{1}L+A)+H_{1}-K_{X}))\otimes {\cal I}(p_{m}^{*}(h_{A}h_{L}^{\nu_{0}})^{m}))\otimes {\cal O}_{Y_{m}}/p_{m}^{*}{\cal I}_{X_{1}}). 
\]
Hence we can extend every element of 
\[
p_{m}^{*}H^{0}(X_{1},{\cal O}_{X_{1}}(m(\nu_{1}L+A)+H_{1})\otimes{\cal I}((h_{A}h_{L}^{\nu_{0}})^{m}))
\]
to an element of 
\[
H^{0}(Y_{m},{\cal O}_{Y_{m}}(K_{Y_{m}}+ p_{m}^{*}(m(\nu_{1}L+A)+H_{1}-K_{X}))\otimes {\cal I}(p_{m}^{*}(h_{A}h_{L}^{\nu_{0}})^{m}))) 
\]
Since 
\[
H^{0}(Y_{m},{\cal O}_{Y_{m}}(K_{Y_{m}}+ p_{m}^{*}(m(\nu_{1}L+A)+H_{1}-K_{X})\otimes {\cal I}(p_{m}^{*}(h_{A}h_{L}^{\nu_{0}})^{m})))
\simeq 
\]
\[
 H^{0}(X,{\cal O}_{X}(m(\nu_{1}L+A)+H_{1})\otimes{\cal I}((h_{A}h_{L}^{\nu_{0}})^{m}))
\]
holds by the isomorphism 
\[
p_{m\, *}({\cal O}_{Y_{m}}(K_{Y_{m}})\otimes{\cal I}(p_{m}^{*}(h_{A}h_{L}^{\nu_{0}})^{m})))= {\cal O}_{X}(K_{X})\otimes {\cal I}((h_{A}h_{L}^{\nu_{0}})^{m}), 
\]
this completes the proof of Lemma 4.4.
\vspace{10mm}{\bf  Q.E.D.} \\ 
Let $\tau_{1}$ be a nonzero element of $H^{0}(X,{\cal O}_{X}(H_{1}))$. 
Let $x_{1}$ be a very general point 
 on $X_{1,reg}$ such that
 $\nu (\Theta_{h_{A}h_{L}^{\nu_{1}}},x_{1}) = 0$ holds. 
Let $m_{1}$ be a sufficiently large positive integer and let
\[
\sigma_{1}^{\prime} \in 
H^{0}(X_{1},{\cal O}_{X_{1}}(m_{1}(\nu_{1}L+A)\otimes{\cal I}((h_{A}h_{L}^{\nu_{0}})^{m_{1}})/tor \otimes {\cal M}_{x_{1}}^{\otimes 2(n+1)^{2}m_{1}})
\]
be a nonzero element.
We note that if $X_{1}$ is smooth (and if we take $x$ very general),   
\[
{\cal O}_{X_{1}}(m_{1}(\nu_{1}L+A)\otimes{\cal I}((h_{A}h_{L}^{\nu_{0}})^{m_{1}})
\]
is torsion free, since it is a subsheaf of a locally free sheaf
on a smooth variety.

Let 
\[
p : \tilde{X} \longrightarrow X
\]
be an embedded resolution and let $X_{1}^{\prime}$ be the strict transform of $X_{1}$. 
We may consider $\sigma_{1}^{\prime}$ as an element of 
\[
H^{0}(X_{1}^{\prime},{\cal O}_{X_{1}^{\prime}}(p^{*}(m_{1}(\nu_{1}L+A))\otimes p^{*}({\cal I}((h_{A}h_{L}^{\nu_{0}})^{m_{1}}) \otimes {\cal M}_{x_{1}}^{\otimes 2(n+1)^{2}m_{1}})).
\]
Hence $\sigma_{1}^{\prime}$ can be lifted to an element of 
\[
H^{0}(X_{1},{\cal O}_{X_{1}}(m_{1}(\nu_{1}L+A)\otimes{\cal I}((h_{A}h_{L}^{\nu_{0}})^{m_{1}}) \otimes {\cal M}_{x_{1}}^{\otimes 2(n+1)^{2}m_{1}}),
\]
if it vanishes on $(p^{*}X_{1} - X_{1}^{\prime})\cap X_{1}^{\prime}$. 
Such a nonzero element $\sigma_{1}^{\prime}$ certainly exists, if $m_{1}$ is sufficiently large. 
Hence we may assume that 
$\sigma_{1}^{\prime}$ is  an element of 
\[
H^{0}(X_{1}^{\prime},{\cal O}_{X_{1}^{\prime}}(p^{*}(m_{1}(\nu_{1}L+A))\otimes p^{*}({\cal I}((h_{A}h_{L}^{\nu_{0}})^{m_{1}}) \otimes {\cal M}_{x_{1}}^{\otimes 2(n+1)^{2}m_{1}})).
\]

Let 
$\sigma_{1}$ be an extension of 
\[
\sigma_{1}^{\prime}\otimes\tau_{1}
\in
H^{0}(X_{1},{\cal O}_{X_{1}}(m_{1}(\nu_{1}L+A)+H_{1})\otimes{\cal I}((h_{A}h_{L}^{\nu_{0}})^{m_{1}}))
\]
to $X$. 
This extension is possible by Lemma 4.4.
Then we set 
\[
h_{1} := \frac{1}{\mid \sigma_{1}\mid^{\frac{2}{r+m_{1}}}}.
\]
Then $h_{1}$ is a singular hermitian metric of $A + \nu_{1}L$
with positive curvature current.

Suppose that $x$ is a regular point of $X_{1}$. 
In this case we shall take $x_{1} = x$. 
Let $\varepsilon_{0}$ be a sufficiently small positive number. 
We define a positive number $\alpha_{1}$ by
\[
\alpha_{1} = \inf\{\alpha > 0 \mid ({\cal O}_{X}/{\cal I}(h_{0}^{\alpha_{0}-\varepsilon_{0}}\cdot h_{1}^{\alpha})_{x}\neq 0\} .
\]
Let us recall the following lemma.
\begin{lemma}(\cite[p.12, Lemma 6]{t})
Let $a,b$ be  positive numbers. Then
\[
\int_{0}^{1}\frac{r_{2}^{2n_{1}-1}}{(r_{1}^{2}+r_{2}^{2a})^{b}}
dr_{2}
=
r_{1}^{\frac{2n_{1}}{a}-2b}\int_{0}^{r_{1}^{-{2}{a}}}
\frac{r_{3}^{2n_{1}-1}}{(1 + r_{3}^{2a})^{b}}dr_{3}
\]
holds, where 
\[
r_{3} = r_{2}/r_{1}^{1/a}.
\]
\end{lemma}

By Lemma 4.5 (if we take $\varepsilon_{0}$ sufficiently small),
we see that 
\[
\alpha_{1}\leq \frac{1}{2n}
\]
holds.
 
Suppose that $x$ is a singular point of $X_{1}$.
In this case letting $x_{1}$ tend to $x$,  
we define the singular hermitian metric
$h_{1}$.    
To estimate $\alpha_{1}$, we use the following lemma (\cite{a-s}). 

\begin{lemma}
Let $\varphi$ be a plurisubharmonic function on $\Delta^{n}\times{\Delta}$.
Let $\varphi_{t}(t\in\Delta )$ be the restriction of $\varphi$ on
$\Delta^{n}\times\{ t\}$.
Assume that $e^{-\varphi_{t}}$ does not belong to $L^{1}_{loc}(\Delta^{n},O)$
for every $t\in \Delta^{*}$.

Then $e^{-\varphi_{0}}$ is not locally integrable at $O\in\Delta^{n}$.
\end{lemma}
Lemma 4.6 is an immediate consequence of the 
$L^{2}$-extension theorem (\cite[p.200, Theorem]{o-t}).
By Lemma 4.6 we have the same estimate
\[
\alpha_{1}\leq \frac{1}{2n}
\]
also in the case that $x$ is a singular point of $X_{1}$.
We define 
\[
V_{2} = \lim_{\delta\downarrow 0}\mbox{Spec}({\cal O}_{X}/{\cal I}(h_{0}^{\alpha_{0}-\varepsilon_{0}}\cdot h_{1}^{\alpha_{1}+\delta}))
\]
and let $X_{2}$ be a branch of $V_{2}$ containing $x$.
By the choice of $x$, 
$h_{L}\mid_{X_{2}}$ is well defined and $(L,h_{L})$ 
is not numerically trivial on $X_{2}$. 

By the above argument, inductively we obtain the strictly decreasing 
sequence of subvarieties:
\[
X = X_{0}\supset X_{1}\supset \cdots X_{r}\supset X_{r+1} = \{ x\}
\]
(the last subvariety $X_{r+1}$ is a point by the 
choice of $x$, i.e. by the {\bf numerical positivity of $(L,h_{L})$ }
at $x$) and the positive numbers $\{\alpha_{i}\}_{i=0}^{r}$ depending 
on small positive numbers 
$\{ \varepsilon_{i}\}_{i=0}^{r-1}$.
Since
\[
\sum_{i=0}^{r}\alpha_{i}\leq \frac{1}{2} 
\] 
holds by the construction, 
we can define a singular hermitian metric $\tilde{h}_{x}$ on 
$mL + A$ for every $m >  \sum_{i=0}^{r}\alpha_{i}\nu_{i}$ by 
\[
\tilde{h}_{x} = (\prod_{i=0}^{r-1}h_{i}^{(\alpha_{i}-\varepsilon_{i})})\cdot h_{r}^{\alpha_{r}+\varepsilon_{r}}\cdot h_{A}^{(1-(\sum_{i=0}^{r-1}(\alpha_{i}-\varepsilon_{i}))-(\alpha_{r}+\varepsilon_{r}))}\cdot
h_{L}^{m-\sum_{i=0}^{r}\alpha_{i}\nu_{i}},
\]
where $\varepsilon_{0},\ldots ,\varepsilon_{r}$ are sufficiently
small positive numbers.
Then the curvature $\Theta_{\tilde{h}_{x}}$ is a closed strictly positive 
$(1,1)$-current on $X$ since 
\[
\Theta_{\tilde{h}_{x}} =  \sum_{i=0}^{r-1}(\alpha_{i}-\varepsilon_{i})\Theta_{h_{i}}
+ (\alpha_{r}+\varepsilon_{r})\Theta_{h_{r}} + (1-(\sum_{i=0}^{r-1}(\alpha_{i}-\varepsilon_{i}))-(\alpha_{r}+\varepsilon_{r}))\Theta_{h_{A}} +
(m-\sum_{i=0}^{r}\alpha_{i}\nu_{i})\Theta_{h_{L}},
\]
\[
\Theta_{h_{i}} (0\leq i\leq r),\Theta_{L} (= \Theta_{h_{L}}) \geq 0, 
\Theta_{A}( = \Theta_{h_{A}}) > 0,
\]
\[
1-(\sum_{i=0}^{r-1}(\alpha_{i}-\varepsilon_{i}))-(\alpha_{r}+\varepsilon_{r})) > 0, 
m-\sum_{i=0}^{r}\alpha_{i}\nu_{i} > 0
\]
hold.

And moreover ${\cal I}(\tilde{h}_{x})$ defines a subscheme of isolated 
support at $x$, if we have taken $x$  to be a very general point on $X$. 

If
\[
{\cal I}(\tilde{h}_{x}) \subseteq {\cal I}(h_{A}h_{L}^{m})
\] 
holds, applying Nadel's vanishing theorem (Theorem 2.1), we see that 
for every $m > \sum \alpha_{i}\nu_{i}$ 
there exists a section 
\[
\sigma \in 
H^{0}(X,{\cal O}_{X}(K_{X}+A+mL)\otimes{\cal I}(h_{A}h_{L}^{m}))
\]
such that $\sigma (x)\neq 0$ holds. 
But at this stage it is not clear that the above inclusion  
holds.
 The reason is that 
$h_{0},\ldots ,h_{r}$ may have weaker singularities than
$h_{A}h_{L}^{\nu_{0}},\ldots ,h_{A}h_{L}^{\nu_{r}}$ 
at some points on $X$ respectively, i.e.
$h_{A}h_{L}^{\nu_{i}}/h_{i} (i= 0,\ldots ,r)$
 may not be bounded on $X$.

Suppose that 
\[
{\cal I}(\tilde{h}_{x}) \not\subseteq {\cal I}(h_{A}h_{L}^{m})
\]
holds. 
Let us fix $m$ such that 
\[
m > \sum_{i=0}^{r}\alpha_{i}\nu_{i}.
\]
and  
\[
\rho_{m} : X^{(m)}\longrightarrow X
\]
be a modification such that $\rho_{m}^{*}{\cal I}(h_{L}^{m}h_{A})$
is locally free on $X^{(m)}$. 
Let $D_{m} = \sum_{k} a_{m,k}D_{m,k}$ be the integral divisor such that 
\[
{\cal O}_{X^{(m)}}(-D_{m}) = \rho_{m}^{*}{\cal I}(h_{L}^{m}h_{A}).
\]
As is seen in Remark 4.4, we have only finitely many choices of  $X_{1},\ldots X_{r}$. 
Hence we may assume that $X_{1},\ldots ,X_{r}$ are 
independent of $m_{0},\ldots , m_{r-1}$.
\begin{sublemma}
There exists a positive constant $C$ such that 
\[
\nu (\rho_{m}^{*}\Theta_{h_{A}h_{L}^{\nu_{i}}},D_{m,k})
-\nu (\rho_{m}^{*}\Theta_{h_{i}},D_{m,k})
\leq \frac{C}{m_{i}}
\]
holds for every $y\in X^{(m)}$,
$0\leq i\leq r$ and $k$.
\end{sublemma}
{\bf Proof}. 
Since 
\[
\sigma_{i} \in 
H^{0}(X,{\cal O}_{X}((m_{i}+r_{i})(\nu_{i}L+A))
\otimes {\cal I}((h_{A}h_{L}^{\nu_{i}})^{m_{i}}))
\]
holds, 
by Lemma 2.2, by the definition of $h_{i}$ we see that 
\[
\nu (\rho_{m}^{*}\Theta_{h_{i}},D_{m,k})
\geq \frac{m_{i}}{m_{i}+r_{i}}
\nu (\rho_{m}^{*}(\Theta_{A}+\nu_{i}\Theta_{L}),D_{m,k})
- \frac{1}{m_{i}+r_{i}}(a_{k}+1)
\]
holds, where $a_{k} \geq 0$ is the coefficient of $D_{m,k}$ in  
the discrepancy $K_{X^{(m)}} - \rho_{m}^{*}K_{X}$, i.e.  
\[
K_{X^{(m)}}- \rho_{m}^{*}K_{X} = 
\sum_{k}a_{k}D_{m,k}
\]
holds. 
Hence there exists a positive constant $C$ such that 
\[
\nu (\rho_{m}^{*}\Theta_{h_{A}h_{L}^{\nu_{i}}},D_{m,k})
-\nu (\rho_{m}^{*}\Theta_{h_{i}},D_{m,k})
\leq \frac{C}{m_{i}}
\]
holds for every $y\in X^{(m)}$,
$0\leq i\leq r$ and $k$. {\bf  Q.E.D.} \vspace{5mm} \\
Sublemma 4.3 means that if we take $m_{i}$ very large, then 
$h_{A}h_{L}^{\nu_{i}}/h_{i}$ has {\bf very small singularities} 
on $X$, even if it is not bounded on $X$. 

To assure the  inclusion
\[
{\cal I}(h_{x})\subseteq {\cal I}(h_{A}h_{L}^{m}),
\]
we modify the argument as follows. 
Let us fix $m > \sum_{i=0}^{r}\alpha_{i}\nu_{i}$.
We set 
\[
S (=S_{m}) := \mbox{Spec}({\cal O}_{X}/{\cal I}(h_{L}^{m}h_{A}))_{red}.
\] 
Let $\varphi$ be an almost plurisubharmonic function which is expressed locally:
\[
\varphi =  \log \sum_{j}\mid f_{j}\mid^{2} + \mbox{$C^{\infty}$-function}, 
\]
where $\{ f_{j}\}$ is a finite set of local generators of the ideal of $S$.  
We set 
\[
h_{x} := \tilde{h}_{x}\cdot e^{-\delta \varphi}, 
\]
where $\delta$ is a small positive number so that 
$\Theta_{h_{x}}$ is strictly positive. 

By this modification and Sublemma 4.3, we see that 
\[
{\cal I}(h_{x}) \subseteq {\cal I}(h_{L}^{m}h_{A})
\]
holds, if $m_{0},\ldots ,m_{r}$ are sufficiently large. 

In fact this  can be verified as follows. 
If we take $\delta$ so that 
\[
\delta >>\sum_{i=0}^{r}\frac{C}{m_{i}}
\]
holds (this does not violate the fact that we need to take 
$\delta$ sufficiently small, because $\{ m_{i}\}$'s 
can be arbitrary large),
\[
\nu (\rho_{m}^{*}\Theta_{h_{x}},D_{m,i})
> \nu (\rho_{m}^{*}(\Theta_{A}+m\Theta_{L}),D_{m,i})
\]
holds for every $i$, if we take $m_{0},\ldots ,m_{r}$ sufficiently 
large. 

Hene by Lemma 2.2, we have the inclusion:
\[
{\cal I}(h_{x}) \hookrightarrow  {\cal I}(h_{L}^{m}h_{A}).
\]
Then by Nadel's vanishing theorem (Theorem 2.1) 
we see that
\[
H^{0}(X,{\cal O}_{X}(K_{X}+mL+A)\otimes{\cal I}(h_{L}^{m}h_{A}))
\neq 0
\]
for every $m > \sum_{i=0}^{r}\alpha_{i}\nu_{i}$.
This completes the proof of Theorem 4.1. {\bf  Q.E.D.} 

\begin{remark}
By modifying the above proof in the first case of Theorem 4.1, 
it is not hard to show that 
\[
H^{0}(X,{\cal O}_{X}(K_{X}+A+mL)\otimes {\cal I}(h_{A}h_{L}^{m}))
\]
gives a birational rational map from $X$ into a projective space
for every sufficiently large $m$.
\end{remark}

\subsection{The second nonvanishing theorem}

In this subsection we shall consider the existence of sections 
of a numerically trivial singular hermitian line bundles 
twisted by some line bundle. 

Let $(L,h)$ be a singular hermitian  line bundle on a smooth projective 
variety $Y$ such that $\Theta_{h}$ is closed positive on $Y$. 
Let $\sum_{i=1}^{r}Z_{i}$ be a divisor with normal crossings
on $Y$. 
Let $X$ be a smooth subvariety defined by
\[
X = Z_{1}\cap\cdots \cap Z_{r}.
\] 
We say such a subvariety $X$ a {\bf transverse complete intersection}
in $Y$.  

Suppose that $(L,h)$ is {\bf numerically trivial} on $Y$. 
Then by Theorem 4.3, we see that 
there exists at most coutably many prime divisors 
$\{ F_{k}\}$ and nonnegative real numbers
$\{ a_{k}\}$ such that
\[
\Theta_{h} = 2\pi\sum_{k}a_{k}F_{k}
\]
holds. 
Let $\xi_{k}$ be a nonzero global section of 
${\cal O}_{Y}(F_{k})$ with divisor $F_{k}$. 
Then we see that there exists a positive constant $C$ 
such that 
\[
h = C \cdot\prod_{k}\frac{1}{\mid\xi_{k}\mid^{2a_{k}}}
\]
holds. 

We shall assume that $\{ Z_{1},\ldots Z_{r}\}$
contains all the divisorial components of 
$\{ y\in Y \mid \nu (\Theta_{h},y) > 0\}$
containing $X$. 
In this case we say that $X$ is {\bf a transverse complete intersection
with respect to $(L,h)$}. 
Then for every $m\geq 1$, we see that 
\[
{\cal O}_{X}(mL)\otimes{\cal I}(h^{m})
= {\cal O}_{X}(mL - \sum_{i=1}^{r}[m\cdot\nu (\Theta_{h},Z_{i})]Z_{i})
\otimes \tilde{\cal I}_{X}(h^{m}),
\]
holds for some ideal sheaf $\tilde{\cal I}_{X}(h^{m})$ on 
$X$, since  the lefthandside is a subsheaf of 
the locally free sheaf ${\cal O}_{X}(mL - \sum_{i=1}^{r}[m\cdot\nu (\Theta_{h},Z_{i})]Z_{i})$ on the smooth variety $X$.
We define the Lelong number $\nu_{X}(\Theta_{h},x) (x\in X)$ 
by 
\[
\nu_{X}(\Theta_{h},x) := 
\overline{\lim}_{m\rightarrow\infty}
m^{-1}\mbox{mult}_{x}\mbox{Spec}({\cal O}_{X}/\tilde{I}_{X}(h^{m})).
\]
And define 
\[
{\cal I}_{X}(h^{m}):= \cap_{m\geq 1}\sqrt[m]{\tilde{\cal I}_{X}(h^{m})},
\]
where 
$\sqrt[m]{\tilde{\cal I}_{X}(h^{m})}$ is defined by
\[
\sqrt[m]{\tilde{\cal I}_{X}(h^{m})}_{x}
:= \cup {\cal I}(\frac{1}{m}(\sigma ))_{x}\hspace{10mm} (x\in X),
\]
where $\sigma$ runs all the germs of $\tilde{\cal I}_{X}(h^{m})_{x}$.

Then by successive use of Corollary 4.2 we see that 
\[
S : = \{ x\in X \mid \nu_{X}(\Theta_{h},x) > 0\}
\]
consists of countably many prime divisors on $X$. 
We set 
\[
S = \sum_{j}D_{j},
\] 
\[
d_{j} := \nu_{X}(\Theta_{h},D_{j})
\]
and 
\[
D: = \sum_{j}d_{j}D_{j}.
\]
Then since $\Theta_{h}$ is a sum of 
at most countably prime divisors with nonnegative 
real coefficients by Theorem 4.3, we see that 
\[
{\cal I}_{X}(h^{m}) = {\cal I}(mD)
\]
holds for every $m\geq 0$, since 
$X$ is a transverse complete intersection 
with respect to $(L,h)$.

The following theorem is as important as Theorem 4.1. 
\begin{theorem}(The second nonvanishing theorem)
Let $Y$ be a smooth projective variety and let 
$(L,h)$ be a numerically trivial singular hermitian line 
bundle on $Y$.
Let $X = Z_{1}\cap\cdots\cap Z_{r}$ be a 
transverse complete intersection subvariety
with respect to $(L,h)$. 
We set 
\[
\nu_{i} : = \nu (\Theta_{h},Z_{i})
\]
and let  $\zeta_{i}$ be nonzero global section in  $\Gamma(Y,{\cal O}_{Y}(Z_{i}))$ with divisor $Z_{i}$ for $1\leq i\leq r$.
For $1\leq i\leq r$ let $h_{Z_{i}}$ be a $C^{\infty}$-hermitian metric on ${\cal O}_{Y}(Z_{i})$ respectively.
Let $\omega$ be a $C^{\infty}$-K\"{a}hler form on 
$X$. 
Let $A = \sum_{k}a_{k}A_{k}$ be a ${\bf R}$-divisor on $X$ and 
let $\tau_{k}$ denotes a nonzero global section of ${\cal O}_{X}(A_{k})$
with divisor $A_{k}$ for every $k$. 
Suppose that the following conditions are satisfied.  
\begin{enumerate}
\item $\mbox{Supp}\, A = \sum_{k}A_{k}$ is a divisor with normal crossings,
\item $\lceil A\rceil$ is effective,
\item
there exists a $C^{\infty}$-hermitian metric
$h_{A - K_{X}}$ with strictly positive curvature on the  
{\bf R}-line bundle $A - K_{X}$ such that
there exists a positive number $\delta$ such that 
for every $m\geq 1$ satisfying 
\[
\{ m\nu_{i}\} \leq \delta
\]
the curvature current of
the singular hermitian metric
\[
h_{m}:= h^{m}\cdot h_{A-K_{X}}\cdot
(\prod_{i=1}^{r}h_{Z_{i}}^{\{ m\nu_{i}\}}\cdot\mid\zeta_{i}\mid^{2m\nu_{i}})\cdot (\prod_{k}\mid\tau_{k}\mid^{-2(\lceil a_{i}\rceil -a_{i})})
\]
on 
${\cal O}_{X}(\lceil A\rceil -K_{X}+ mL -\sum_{i}[m\nu_{i}]Z_{i})$,
satisfies the inequality
\[
\Theta_{h_{m}} > c\cdot\omega,
\]
where $c$ is a positive number independent of such $m$.
\end{enumerate}
Then  there are finitely  many  positive numbers $t_{0}$, 
$\alpha_{1},\ldots ,\alpha_{\ell}$, $\beta_{1},\ldots ,\beta_{p}$ and a small positive number
$\varepsilon$  such that 
for every   $m \geq t_{0}$ satisfying 
the inequalities
\[
\mid \langle m\alpha_{s}\rangle -m\alpha_{s}\mid 
< \varepsilon  \hspace{15mm}  (1\leq s\leq \ell),
\]
\[
0 < \mid\lceil m\beta_{q}\rceil- m\beta_{q}\mid < \varepsilon  \hspace{15mm} (1\leq q \leq p)
\]
and 
$\{ m\nu_{i}\} \leq \delta \hspace{5mm} (1\leq i\leq r)$, 
\[
H^{0}(X,{\cal O}_{X}(\lceil A\rceil + mL)\otimes {\cal I}_{X}(h_{m}))
\neq 0
\]
holds.
In other words for such $m\geq t_{0}$   
\[
H^{0}(X,{\cal O}_{X}(\lceil A\rceil + mL)\otimes {\cal I}_{X}((\prod_{i=1}^{r}h_{Z_{i}}^{\{ m\nu_{i}\}}\cdot\mid\zeta_{i}\mid^{2m\nu_{i}})
\cdot h_{\lceil A\rceil}\cdot h^{m})) 
\neq 0
\]
holds, where $h_{\lceil A\rceil}$ denotes the singular hermitian metric on 
${\cal O}_{X}(\lceil A\rceil )$ defined by 
\[
h_{\lceil A\rceil} =  \prod_{k}h_{k}^{a_{i}}\mid\tau_{k}\mid^{-2(\lceil a_{i}\rceil -a_{i})}.
\]
Moreover the set of such $m$ is nonempty and infinite. 
\end{theorem}
{\bf Proof of Theorem 4.4}.
Let $n$ denote $\dim X$. 
We prove this theorem by induction on $n$. 
If $n=1$, by the assumption  
for every $m\geq 1$ such that 
\[
\{ m\nu_{i}\} < \delta \hspace{15mm} (1\leq i\leq r) 
\]
hold,
\[
\deg_{X}{\cal O}_{X}(\lceil A\rceil + mL)\otimes{\cal I}(h_{m})
> \deg_{X}K_{X}
\]
holds (such a positive integer $m$ certainly exists by Lemma 4.7 below). 
Since $\dim X = 1$, 
${\cal O}_{X}(\lceil A\rceil + mL)\otimes{\cal I}(h_{m})$ 
is an invertible sheaf on $X$.
Hence for such $m$ 
\[
H^{0}(X,{\cal O}_{X}(\lceil A\rceil + mL)\otimes{\cal I}(h_{m}))
\neq 0
\]
holds  by the Kodaira vanishing theorem and the Riemann-Roch 
theorem. 

Suppose that the theorem holds for every $X$ 
with $\dim X < n$.
Let us consider the case that $\dim X = n$. 
We set 
\[
S : = \{ x \in X\mid \nu_{X}(\Theta_{h},x) > 0\} .
\]
Then $S$ consists of at most countably many prime divisors on $X$. 
Let  
\[
S : = \sum_{j\in J} D_{j}
\]
be the irreducible decomposition. 
We set  
\[
d_{j} := \nu_{X}(\Theta_{h},D_{j}).
\]

We shall consider the following two cases. \vspace{5mm} \\
{\bf Case 1}: $\sharp J = \infty$, \\
{\bf Case 2}: $\sharp J < \infty$. \vspace{5mm} \\
First let us consider  Case 1. 
Let $H$ be a very ample smooth  divisor on $X$.
We may assume that $2\pi\omega$ is a 1-st Chern form of 
${\cal O}_{X}(H)$.  
Now we have the following sublemma.
\begin{sublemma}
There exists a positive number $\delta_{H}$ such that 
for every effective {\bf R}-divisor $E$ on $X$ such that 
\[
H^{n-1}\cdot E < \delta_{H}
\]
$H - E$ is ample. 
\end{sublemma}
{\bf  Proof of Sublemma 4.4.} 
Since for every positive number $a$
\[
\{ c_{1}(E) \in H^{2}(X,\mbox{\bf R})\mid E \, \mbox{: effective {\bf R}-divisor}, H^{n-1}\cdot E < a \}
\]
 is relatively compact in $H^{2}(X,\mbox{\bf R})$,
This follows from Kleinman's criterion for ampleness. 
{\bf Q.E.D.} \vspace{5mm} \\
Then there exists some $D_{j}$, say $D_{0}$ such that
\begin{enumerate}
\item $D_{0}\not\subset\mbox{Supp}\, A$,
\item $d_{0}\cdot H^{n-1}\cdot D_{0} << \frac{1}{2}c\cdot\delta_{H}$ 
\end{enumerate}
hold. 
Since $H^{n-1}\cdot (\sum_{j\in J}d_{j}D_{j})$ is finite and 
$\sharp J$ is infinite,
there exists such $D_{0}$.
We may assume that $D_{0}$ is a smooth divisor. 
In fact let 
\[
\pi_{Y} : \tilde{Y}\longrightarrow Y
\]
be an embedded resolution of $D_{0}$ obtained 
by successive blow ups with smooth centers. 
Let $\tilde{D}_{0}$ denote the strict transform of $D_{0}$
in $\tilde{X}$. 
Let $\tilde{X}$ be the strict transform of $X$ in $\tilde{Y}$ and
let 
\[
\pi : \tilde{X}\longrightarrow X
\]
be the restriction of $\pi_{Y}$ to $\tilde{X}$.
Let $E$ be the effective divisor defined by 
\[
E := K_{\tilde{X}} - \pi^{*}K_{X}. 
\]
We define the divisor $\tilde{A}$ on $\tilde{X}$ by
\[
\tilde{A} := \pi^{*}A + E.
\]
If we take $\pi_{Y}$ properly we may assume that 
$\mbox{Supp}\,\tilde{A}$ is a divisor with normal crossings.
We note that $\lceil\tilde{A}\rceil$ is effective, if and only 
if ${\cal I}(-\tilde{A})\simeq {\cal O}_{\tilde{X}}$ holds.
We note that that 
\[
K_{\tilde{X}}-\tilde{A} = \pi^{*}(K_{X}- A)
\]
holds by the definition of $\tilde{A}$. 
Hence $\lceil \tilde{A}\rceil$ is also effective.
Also the above formula implies that 
for every $m\geq 1$ such that 
\[
\{ m\nu_{i}\} < \delta \hspace{15mm} (1\leq i\leq r), 
\]
there exists a $C^{\infty}$-hermitian metric
$h_{\tilde{A} - K_{\tilde{X}}}$ on the  
{\bf R}-line bundle $A - K_{X}$ such that 
\[
\tilde{h}_{m}:= \pi^{*}h^{m}\cdot h_{\tilde{A}-K_{\tilde{X}}}\cdot
\pi^{*}(\prod_{i=1}^{r}h_{Z_{i}}^{\{ m\nu_{i}\}}\cdot\mid\zeta_{i}\mid^{2m\nu_{i}})\cdot
\pi^{*}(\prod_{k}\mid\tau_{k}\mid^{-2(\lceil a_{i}\rceil -a_{i})})
\]
is a singular hermitian metric on 
${\cal O}_{X}(\lceil \tilde{A}\rceil -K_{\tilde{X}}+ m\pi^{*}L -
\sum_{i}[m\nu_{i}]Z_{i})$ such that the curvature current 
$\Theta_{\tilde{h}_{m}}$ 
satisfies the inequality 
\[
\Theta_{\tilde{h}_{m}} \geq c\cdot\pi^{*}\omega .
\]
We note that there exists an effective divisor $E^{\prime}$ supported 
on $E_{red}$ and a $C^{\infty}$-hermitian metric $h_{E^{\prime}}$ on ${\cal O}_{\tilde{X}}(E^{\prime})$ such that 
\[
\tilde{\omega} (= \tilde{\omega}(\lambda )) = c\cdot\pi^{*}\omega + 
\lambda\cdot\Theta_{h_{E^{\prime}}}
\]
is a $C^{\infty}$-K\"{a}hler form on $\tilde{X}$ for every sufficiently small 
positive number $\lambda$. 
Let $\lambda_{0}$ be a sufficiently small positive number such that 
$\lceil \tilde{A} - \lambda_{0}E^{\prime}\rceil $
is effective and $\tilde{\omega}(\lambda_{0})$ is 
a K\"{a}hler form on $\tilde{X}$. 
Then  if we replace $(Y,X,(L,h),D_{0},A)$ by 
$(\tilde{Y},\tilde{X},\pi^{*}_{Y}(L,h),\tilde{D}_{0},\tilde{A}-\lambda_{0}E^{\prime})$, 
we may assume that $D_{0}$ is smooth, since by the definition of $\tilde{A}$
\[
\pi^{*}H^{0}(\tilde{X},{\cal O}_{\tilde{X}}(\lceil \tilde{A}-\lambda_{0}E^{\prime}\rceil + mL)\otimes {\cal I}_{\tilde{X}}(\pi^{*}((\prod_{i=1}^{r}h_{Z_{i}}^{\{ m\nu_{i}\}}\cdot\mid\zeta_{i}\mid^{2m\nu_{i}})
\cdot h_{\lceil A\rceil}\cdot h^{m})))  
\]
is contained in 
\[
H^{0}(X,{\cal O}_{X}(\lceil A\rceil + m\pi^{*}L)\otimes {\cal I}_{X}((\prod_{i=1}^{r}h_{Z_{i}}^{\{ m\nu_{i}\}}\cdot\mid\zeta_{i}\mid^{2m\nu_{i}})
\cdot h_{\lceil A\rceil}\cdot h^{m})) 
\] 
holds. 
Now we assume that $D_{0}$ is smooth.
 
We need the following lemma. 
\begin{lemma}
Let $a_{1},\ldots a_{\ell}$ are positive numbers. 
Let us consider the sequence 
\[
\{([ma_{1}],\ldots ,[ma_{\ell}])\mid m\in \mbox{\bf N}\} 
\]
in $(\mbox{\bf R}/\mbox{\bf Z})^{\ell}$.
Then there exist a connected subgroup $T$   of 
$(\mbox{\bf R}/\mbox{\bf Z})^{\ell}$
such that 
\[
T\cap \{([ma_{1}],\ldots ,[ma_{\ell}])\mid m\in \mbox{\bf N}\} 
\]
is dense in $T$. 
\end{lemma} 
{\bf Proof of Lemma 4.7.} 
Let $T^{\prime}$ denote the closure of $\{([ma_{1}],\ldots ,[ma_{\ell}])\mid m\in \mbox{\bf N}\}$.
Then $T^{\prime}$ contains the origin.
Let $T$ be the connected component containing the identity.
Then by the definition of $T$, for every 
positive integer $m$,  $([ma_{1}],\ldots ,[ma_{\ell}])$ 
has an inverse in $T$. 
This means that $T$ is a connected subgroup of 
$(\mbox{\bf R}/\mbox{\bf Z})^{\ell}$.
It is clear that $T$ satisfies the desired property. 
{\bf Q.E.D.} \vspace{5mm}\\
Suppose that $m$ is a positive integer such that 
 $\{ -md_{0}\}$ is a {\bf nonzero number} satisfying  
\[
\{ -md_{0}\}\cdot (H^{n-1}\cdot D_{0}) < \frac{1}{2}c\cdot\delta_{H}
\]
and
\[
\{ m\nu_{i}\} \leq \delta \hspace{15mm} (1\leq i\leq r)
\]
hold. 
By Lemma 4.7 such a positive integer $m$ exists, if we take $d_{0}$ sufficiently small.  
For such $m$, by the assumption  
there exists a $C^{\infty}$ hermitian metric 
$h_{A - K_{X}}$ on the  
{\bf R}-line bundle $A - K_{X}$ such that 
\[
h_{m}:= h^{m}\cdot h_{A-K_{X}}\cdot
(\prod_{i=1}^{r}h_{Z_{i}}^{\{ m\nu_{i}\}}\cdot\mid\zeta_{i}\mid^{2m\nu_{i}})
 \cdot (\prod_{k}\mid\tau_{k}\mid^{-2(\lceil a_{i}\rceil -a_{i})})
\]
is a singular hermitian metric on 
${\cal O}_{X}(\lceil A\rceil -K_{X}+ mL -\sum_{i=1}^{r}[m\nu_{i}]Z_{i})$ such that $\Theta_{h_{m}}$ satisfies the inequality
\[
\Theta_{h_{m}} \geq c\cdot \omega.
\]
Let $h_{0}$ be a $C^{\infty}$-hermitian metric on 
${\cal O}_{X}(D_{0})$ and let $\sigma_{0}\in 
\Gamma (X,{\cal O}_{X}(D_{0}))$ be a nonzero global section 
with divisor $D_{0}$.
Then if we take $h_{0}$ properly, by 
the choice of $D_{0}$ and  Sublemma 4.4 
\[
\hat{h}_{m}:=  h_{m}
\cdot (\frac{1}{h_{0}(\sigma_{0},\sigma_{0})^{\{-md_{0}\}}})
\]
is a singular hermitian metric on 
${\cal O}_{X}(\lceil A\rceil -K_{X}+ mL -\sum_{i=1}^{r}[m\nu_{i}]Z_{i})$ 
such that the curvature current $\Theta_{\hat{h}_{m}}$ 
satisfies the inequality 
\[
\Theta_{\hat{h}_{m}}\geq \frac{1}{2}c\cdot \omega
\]
on $X$.
This implies that 
\[
H^{1}(X,{\cal O}_{X}(\lceil A\rceil+ mL)
\otimes{\cal I}_{X}(\hat{h}_{m})) = 0
\]
holds. 
Hence 
\[
H^{0}(X,{\cal O}_{X}(\lceil A\rceil+ D_{0}+ mL)
\otimes{\cal I}_{X}(\hat{h}_{m})) 
\rightarrow 
H^{0}(D_{0},{\cal O}_{D_{0}}(\lceil A\rceil+ D_{0}+ mL)
\otimes{\cal I}_{X}(\hat{h}_{m}))
\]
is surjective for such $m$. 
By the construction of $\hat{h}_{m}$ we see that 
\[
{\cal I}_{X}(\hat{h}_{m})\subseteq 
{\cal I}_{X}(h_{\lceil A\rceil}\cdot h^{m})\otimes{\cal O}_{X}(-D_{0})
\]
holds.
In fact this can be verified as follows.
First by the definition of $\hat{h}_{m}$, it is clear that 
\[
{\cal I}_{X}(\hat{h}_{m})\subseteq 
{\cal I}_{X}(h_{\lceil A\rceil}\cdot h^{m})
\]
holds. 
Let $\varpi : X^{\prime} \longrightarrow X$ be 
a modification such that 
$\varpi^{*}{\cal I}_{X}(\hat{h}_{m})$ and 
$\varpi^{*}{\cal I}_{X}(h_{\lceil A\rceil}\cdot h^{m})$
are locally free.
Take a sufficiently small open set $U$ in $X$ and  any local section $\sigma$ of 
${\cal I}_{X}(\hat{h}_{m})(U)$.
Then $\varpi^{*}\sigma$ is 
an element of 
$\varpi^{*}{\cal I}_{X}(h_{\lceil A\rceil}\cdot h^{m})(\varpi^{-1}(U))$ 
which is identically $0$ on the strict transform of $D_{0}$
in $\varpi^{-1}(U)$.
Hence we have that 
$\sigma\in ({\cal I}_{X}(h_{\lceil A\rceil}\cdot h^{m})\otimes{\cal O}_{X}(-D_{0}))(U)$ holds.

Since $D_{0}$ is not contained in the support of 
$A$, we see that 
$\lceil A\mid_{D_{0}}\rceil$ is effective. 
Since $K_{D_{0}}= (K_{X}+D_{0})\mid_{D_{0}}$ holds
by the adjunction formula, we see that for such $m$, 
there exists a $C^{\infty}$-hermitian metric
$h_{A\mid_{D_{0}} - K_{D_{0}}}$ on the  
{\bf R}-line bundle $A\mid_{D_{0}} - K_{D_{0}}$ such that 
\[
h_{m,D_{0}}:= \hat{h}_{m}\mid_{D_{0}}= 
 h^{m}\cdot h_{A\mid_{D_{0}}-K_{D_{0}}}\cdot
(\prod_{i}h_{Z_{i}}^{\{ m\nu_{i}\}}\cdot\mid\zeta_{i}\mid^{2m\nu_{i}})\cdot 
(\prod_{k}\mid\tau_{k}\mid^{-2(\lceil a_{i}\rceil -a_{i})})
\cdot (\frac{1}{h_{0}(\sigma_{0},\sigma_{0})^{\{-md_{0}\}}})
\]
is a singular hermitian metric on 
\[
{\cal O}_{D_{0}}(\lceil A\rceil -K_{X}+ mL 
-\sum_{i=1}^{r}[m\nu_{i}]Z_{i}
-[md_{0}]\cdot D_{0})
\]
such that
the curvature current $\Theta_{h_{m},D_{0}}$ satisfies the inequality
\[
\Theta_{h_{m,D_{0}}}\geq \frac{1}{2}c\cdot\omega_{\mid D_{0}}.
\]
Since 
\[
{\cal I}(h_{m,D_{0}})\subseteq {\cal I}(\hat{h}_{m})
\otimes{\cal O}_{D_{0}}
\]
holds by the $L^{2}$-extension theorem (\cite{o-t}, cf. Lemma 4.6 above),
\[
H^{0}(D_{0},{\cal O}_{D_{0}}(\lceil A\mid_{D_{0}}\rceil
+ mL)\otimes {\cal I}(h_{m,D_{0}}))
\subseteq 
H^{0}(D_{0},{\cal O}_{D_{0}}(\lceil A\mid_{D_{0}}\rceil
+ mL)\otimes {\cal I}(\hat{h}_{m}))
\]
holds.
Hence if 
\[
H^{0}(D_{0},{\cal O}_{D_{0}}(\lceil A\mid_{D_{0}}\rceil
+ mL)\otimes {\cal I}(h_{m,D_{0}}))\neq 0
\]
holds, then by the above argument, we see that 
\[
H^{0}(X,{\cal O}_{X}(\lceil A\rceil +mL)\otimes{\cal I}_{X}(h_{m}))\neq 0
\]
holds.  

By the definition of $D_{0}$, $\lceil A\mid_{D_{0}}\rceil$ is effective.  
Hence repeating the same procedure, we may continue the argument and 
reduce the problem to the (1-dimension) lower dimensional case.
We note that $D_{0}$ may not be a transverse complete intersection with respect to $(L,h)$ in $Y$, but it is a smooth divisor on $X$. 
Hence essentially we may apply the inducton.
In fact only difference is that we should consider the {\bf R}-line bundle
$m(L -\sum_{i=1}^{r} \nu_{i}D_{i})$ on $X$ instead of $mL$. But as above this does not matter, if the residual divisor $\sum_{i=1}^{r}\{ m\nu_{i}\}D_{i}$ is sufficiently small.
By the inductive argument (see also the second case below), setting $d_{0}$ to be one of 
$\{ \beta_{1},\ldots ,\beta_{p}\}$ in the statement of Theorem 4.4,  we completes the proof of 
Theorem 4.4 in this case.  \vspace{5mm} \\

Next let us consider the second case, i.e. the case that 
$\sharp J < \infty$.    
By taking a suitable modification of $X$ in $Y$, 
we may assume that 
\[
A + \sum_{j\in J}d_{j}D_{j}
\]
is a divisor with normal crossings on $X$. 

In this case we quote the following theorem. 
\begin{theorem}(\cite[p. 427,Theorem 3]{ka})
Let $M$ be a smooth projective variety and let 
$A$ be an divisor on $X$ with real coefficients such that
\begin{enumerate}
\item $\mbox{Supp}\{ A\}$ is a divisor with normal crossings, 
\item $\lceil A\rceil$ is effective,
\item $A - K_{M}$ is ample.
\end{enumerate}
Let $L$ be a line bundle and let $D = \sum d_{j}D_{j}$ be an effective divisor
with real coefficients on $M$
 such that $L - D$ is nef and $\mbox{Supp}\, D$ is 
a divisor with normal crossings.
Then there exist positive numbers $t_{0}$ and $\varepsilon$  such that 
for every integer $m$ satisfying 
\[
m \geq t_{0}
\,\,\,\, \mbox{and}\,\,\,\, 
\mid\langle md_{j}\rangle - md_{j}\mid < \varepsilon ,
\]
\[
H^{0}(M,{\cal O}_{M}(\lceil A\rceil- \langle mD\rangle +mL))\neq 0
\]
holds, where for a real number $d$, $\langle d \rangle$ denotes
the integer such that 
\[
d - \frac{1}{2} \leq  \langle d\rangle < d + \frac{1}{2}  
\]
and 
\[
\langle mD\rangle := \sum_{j}\langle md_{j}\rangle D_{j}.
\]
\end{theorem}
\begin{remark}
By Lemma 4.7, the set of $m$ satisfying the inequalities 
in Theorem 4.5 is nonempty and infinite. 
\end{remark}
Let us continue the proof of Theorem 4.4.
We note that  
\[
L - \sum_{i=1}^{r}\nu_{i}Z_{i}- \sum_{j\in J}d_{j}D_{j}
\]
is numerically trivial on $X$.
Also by the assumption for every positive number $\varepsilon$
for every $1\leq i\leq r$ and $j\in J$
there exists a positive integer $m$
such that 
\[
\mid\langle md_{j}\rangle - md_{j}\mid < \varepsilon \hspace{20mm} (j\in J)
\]
and 
\[
\mid\lceil m\nu_{i}\rceil\mid < \delta  \hspace{20mm} (1\leq i\leq r)
\]
hold. 
The existence of such $m$ follows from Lemma 4.7. 

We cannot apply Theorem 4.5 directly in our situation,
since ${\cal O}_{X}(Z_{i})\mid_{X}$ is not effective. 
To use the Cartier divisor $mL - \sum_{i=1}^{r}[m\nu_{i}]Z_{i}$ 
instead of the {\bf R}-divisor $mL -\sum_{i=1}^{r}m\nu_{i}Z_{i}$,
we need to dispose of  the residual divisor
$\sum_{i=1}^{r}\{ m\nu_{i}\} Z_{i}$. 
This residual divisor can be absorbed in $A$ in the 
following manner, if 
$\{ m\nu_{i}\} ,1\leq i\leq r$ are sufficiently small. 
Let us take a very ample divisor $H$ as above. 
We may assume that $A + H$ is a divisor with normal crossings. 
Let us take a positive rational number $\varepsilon_{0}$ so that 
 $A - \varepsilon_{0}H - K_{X}$ is ample. 
Then there exists a positive rational number $\delta_{0}$ 
such that if 
\[
\{ m\nu_{i}\} < \delta_{0} (1\leq i\leq r),
\]
then 
\[
\varepsilon_{0}H - \sum_{i=1}^{r}\{ m\nu_{i}\} Z_{i} 
\]
is {\bf Q}-linearly equivalent to an ample effective
{\bf Q}-divisor $B$. 
We may assume that $A + H +B$ is a divisor with normal
crossings. 
Then we see that 
\[
 \sum_{i=1}^{r}\{ m\nu_{i}\} Z_{i} 
\sim_{\mbox{\bf Q}} \varepsilon_{0} H - B
\]
and 
\[
 A - \varepsilon_{0}H + \sum_{i=1}^{r}\{ m\nu_{i}\}
Z_{i} 
\sim_{\mbox{\bf Q}} A - B 
\]
hold, where $\sim_{\mbox{\bf Q}}$ denotes the 
{\bf Q}-linear equivalence relation. 
Also we note that we may assume that 
$\lceil A - B\rceil = \lceil A\rceil$ holds. 
Thus if $\{ m\nu_{i}\} < \delta_{0}$ holds for every 
$1\leq i\leq r$, we may neglect the residual divisor 
$ \sum_{i=1}^{r}\{ m\nu_{i}\}Z_{i}$ by the 
perturbation of the divisor $A$.  
This argument has already been used in \cite{ka}
to prove Theorem 4.5. 
The essential part of the proof of Theorem 4.5 is 
this argument and the rest of the proof is parallel 
to the proof of the Shokurov's nonvanishing theorem 
(\cite{sh}). 
In this way we can dispose of {\bf R}-divisors, 
if the residual part is sufficiently small. 

Then by Theorem 4.5 if we replace $\delta$ by $\min (\delta ,\delta_{0})$, we see that 
there exists a positive numbers $t_{0}$ and  $\varepsilon$ such that 
\[
H^{0}(X,{\cal O}_{X}(\lceil A\rceil + mL - [m\cdot\sum \nu_{i}Z_{i}]
-[m\sum_{j} d_{j}D_{j}]))
\neq 0  
\]
holds, if 
\[
\mid\langle md_{j}\rangle-md_{j}\mid < \varepsilon \hspace{20mm}
(j\in J)
\]
and 
\[
\{ m\nu_{i}\} < \delta 
\hspace{20mm}(1\leq i\leq r)
\]
and $m\geq t_{0}$.
This completes the proof of Theorem 4.4.

 {\bf Q.E.D.}
\begin{remark}
We note that in  Theorem 4.5
for every flat line bundle $F$ on $M$ and 
a positive integer $m$ such that 
\[
m \geq t_{0}
\,\,\,\, \mbox{and}\,\,\,\, 
\mid\langle md_{j}\rangle - md_{j}\mid < \varepsilon ,
\]
\[
H^{q}(M,{\cal O}_{M}(\lceil A\rceil- \langle mD\rangle +mL+F)) = 0
\]
holds for every $q \geq 1$ by the assumption.
And we note that the curvature of the singular hermitian 
line bundle is stable under tensoring a flat line bundle.
Hence we see that 
Theorem 4.4 holds also for the nonvanishing of 
\[
H^{0}(X,{\cal O}_{X}(\lceil A\rceil + mL+F_{X})\otimes {\cal I}_{X}(h_{m})),
\]
where $F_{X}$ is an arbitrary flat line bundle on $X$. 
\end{remark}

\subsection{Volume of the stable fixed component}

We call the set 
\[
\mbox{SBs}(K_{X}) := \cap_{\ell\geq 1}\mbox{Supp}\,\mbox{Bs}\mid mK_{X}\mid
\]
the stable base locus of $K_{X}$. 
\begin{theorem}
Let $V$ be a divisorial component of $\mbox{SBs}(K_{X})$.
Then 
\[
\mu (V,K_{X}) (= \mu (V,(K_{X},h))) = 0
\]
holds (for the definition of $\mu (V,K_{X})$ see 
Definition 3.5).
\end{theorem}
{\bf Proof}.
Let $V$ be a divisorial component of $\mbox{SBs}(K_{X})$.
Taking an embedded resolution of $V$, we may assume that 
$V$ is smooth. 
Suppose that $\mu (V,K_{X}) > 0$ holds.

Let $m_{0}$ be a  positive integer such that 
$\Phi_{\mid m_{0}K_{X}\mid}$ gives a birational rational map
onto its image.
By taking a suitablve modification, we may assume that 
$\mbox{Bs}\mid m_{0}K_{X}\mid$ is a divisor.
Let 
\[
\mid m_{0}K_{X}\mid = \mid P\mid + F
\]
be the decomposition into the free part $\mid P\mid$
and the fixed component $F$.
Let $h_{P}$ be a $C^{\infty}$-hermitian metric on
${\cal O}_{X}(P)$ with semipositive curvature 
defined by a pull back of the Fubini-Study metric 
on ${\cal O}(1)$ by $\Phi_{\mid P\mid}$.
We set 
\[
r_{V}^{\prime} = \mbox{mult}_{V}\mbox{Bs}\mid m_{0}K_{X}\mid 
- m_{0}\nu (\Theta_{h},V)
\]
and 
\[
r_{V} = \left\{ \begin{array}{cl}
r_{V}^{\prime} & \mbox{if $r^{\prime}_{V} > 0$} \\
1   &  \mbox{if $r^{\prime}_{V} = 0$}
\end{array}
\right.
\]
We note that if $\nu (\Theta_{h},V)$ is $0$,
$r^{\prime}_{V} \geq 1$ holds by the assumption.

By Kodaira's lemma (\cite[Appendix]{k-o}) there exists an effective {\bf  Q}-divisor 
$E$ such that $P - E$ is positive.
Let $h_{P,E}$ be a $C^{\infty}$-hermitian metric on the 
{\bf  Q}-line bundle ${\cal O}_{X}(P - E)$
with strictly positive curvature.
Then $h_{P,E}$ is considered as a singular hermitian metric 
on $m_{0}K_{X}$ as follows. 
Let $a$ be a positive integer such that $aE$ is an integral divisor. 
Let $\sigma$ be a global section of ${\cal O}_{X}(aE)$ with divisor $aE$.
Then 
\[
\hat{h}_{P,E}:= \frac{h_{P,E}}{\mid\sigma\mid^{\frac{2}{a}}}
\]
is a singular hermitian metric on $m_{0}K_{X}$ with strictly positive
curvature. 
Let $\varepsilon$ be a sufficiently small positive number.
For $m > m_{0}/r_{V}$, we set
\[
h_{m,\varepsilon} =  h_{P}^{\frac{1}{r_{V}}}\cdot
h^{(m-\frac{m_{0}}{r_{V}}-m_{0}\varepsilon )}\cdot \hat{h}_{P,E}^{\varepsilon}.
\]
Then $h_{m,\varepsilon}$ is a singular hermitian metric
on $mK_{X}$ with strictly positive curvature.
By Nadel's vanishing theorem (Theorem 2.1) we have that 
\[
H^{1}(X,{\cal O}_{X}((m+1)K_{X})\otimes{\cal I}(h_{m,\varepsilon}))
= 0
\]
holds.
Hence
\[
H^{0}(X,{\cal O}_{X}((m+1)K_{X}+V)\otimes{\cal I}(h_{m,\varepsilon})) 
\rightarrow H^{0}(V,
{\cal O}_{V}((m+1)K_{X}+V)\otimes{\cal I}(h_{m,\varepsilon}))
\]
is surjective. 
Let $G$ be any effective divisor on $V$.
Since $\mu (V,K_{X})$ is positive, 
\[
H^{0}(V,
{\cal O}_{V}(mK_{X}-G)\otimes{\cal I}(h^{m}))\neq 0
\]
holds for every sufficiently large $m$.
Since $G$ is an arbitrary effective divisor on $V$, this implies that 
for any fixed effective divisor $G^{\prime}$ on $V$,
\[
H^{0}(V,
{\cal O}_{V}(mK_{X}-G^{\prime})\otimes{\cal I}(h_{m,\varepsilon}))\neq 0
\]
holds for every sufficiently large $m$. 
Hence we see that 
\[
H^{0}(V,
{\cal O}_{V}((m+1)K_{X}+V)\otimes{\cal I}(h_{m,\varepsilon}))\neq 0
\]
holds for every sufficiently large $m$. 
We note that by the definition of $h_{m,\varepsilon}$
 on the generic point of $V$
\[
{\cal O}_{V}((m+1)K_{X}+V)\otimes{\cal I}(h_{m,\varepsilon})
\]
is a subsheaf of  
\[
{\cal O}_{V}((m+1)K_{X}-[(m-1)\cdot\nu (\Theta_{h},V)V]),
\]
if $r^{\prime}_{V} > 0$ and 
is a subsheaf of 
\[
{\cal O}_{V}((m+1)K_{X}-[((m-1)\cdot\nu (\Theta_{h},V)-1)V]),
\]
if $r^{\prime}_{V} = 0$.

If $\nu (\Theta_{h},V)$ is positive,
taking $\varepsilon$ sufficiently small, we see that 
\[
\mbox{mult}_{V}\mbox{Bs}\mid (m+1)K_{X}\mid
< m\cdot\nu (\Theta_{h},V)
\]
holds. 
This is the contradiciton.
If $\nu (\Theta_{h},V)$ is $0$, then we see that 
\[
\mbox{mult}_{V}\mbox{Bs}\mid (m+1)K_{X}\mid = 0
\]
holds.
This also contradicts the assumption that 
$V$ is in $\mbox{SBs}(K_{X})$. {\bf  Q.E.D.} 

\subsection{Fibration theorem}

Using Theorem 4.1, we have the following theorem.

\begin{theorem} 
Let $X$ be a smooth projective variety of general type 
and let $h$ be an AZD of $K_{X}$. 
Let $F$ be a divisorial irreducible component
of the stable base locus of $K_{X}$.
Then $(K_{X},h)$ defines a nontrivial numerically trivial 
fiber space structure on $F$, i.e. there exists a unique (up to birational equivalence) rational fibration 
\[
f : F -\cdots\rightarrow W
\]
 such that
\begin{enumerate}
\item for a very general fiber $V$, $(K_{X},h)$ is 
numerically trivial on $V$, where $h$ is an AZD of $K_{X}$,
\item $\dim W$ is minimal among such fibrations 
and is less than $\dim F$,
\item for a very general point $x$ on $F$ and any irreducible 
horizontal curve 
$C\subset F$ containing $x$, $(K_{X},h)\cdot C > 0$  holds. 
\end{enumerate}
Moreover if $F$ is smooth, then $f$ is regular over the 
generic point of $W$. 
\end{theorem}  
{\bf Proof.}
By taking an embedded resolution, we may assume that 
$F$ is smooth. 
Let $m_{0}$ be a positive integer such that 
$\Phi_{\mid m_{0}K_{X}\mid}$ is a birational rational map onto 
its image.
Taking a suitable modification, if necessary, we may assume that
$\mbox{Bs}\mid m_{0}K_{X}\mid$ is a divisor 
with normal crossings.
Let 
\[
\mid m_{0}K_{X}\mid = \mid P\mid + \sum_{i}a_{i}D_{i}
\]
be the decomposition of $\mid m_{0}K_{X}\mid$ into the 
free part $\mid P\mid $ and the fixed component 
$\sum_{i} a_{i}D_{i}$.

Taking a suitable modification, if necessary,
there exists a divisor $E = \sum_{j}E_{j}$ and positive numbers $\{ \delta_{i}\} ,\{ \delta_{j}\}$ such that 
\[
P^{*}= P - \sum_{i} \delta_{i}D_{i} -\sum_{j} \delta_{j}E_{j}
\]
is ample and $\sum_{i} \delta_{i}D_{i} + \sum_{j} \delta_{j}E_{j}$
is a divisor with normal crossings.
By Kleiman's criterion for ampleness we may assume that 
$\{\delta_{i}\} ,\{ \delta_{j}\}$ are sifficietnly small positive numbers. 
Let $h_{P^{*}}$ be a $C^{\infty}$hermitian metric 
on the {\bf  R}-line bundle ${\cal O}(P^{*})$ with strictly positive
curvature.

Let $h$ be an AZD of $K_{X}$. 
If $F$ is not contained in $\mbox{Bs}\mid m_{0}(K_{X},h)\mid$ (cf. Definition 3.4),  by Theorem 4.6, $F$ is blown down 
by $\Phi_{\mid m_{0}K_{X}\mid}$ and this defines a numerically trivial fiber space sturcture 
of $F$.
Hence $(K_{X},h)$ defines a nontrivial numerically trivial 
fiber space structure on $F$. 

Next suppose that $F$ is contained in 
$\mbox{Bs}\mid m_{0}(K_{X},h)\mid$.
By changing the indices we may assume that $F = D_{0}$ holds.
We  set  
\[
r = \mbox{mult}_{F}\mbox{Bs}\mid m_{0}K_{X}\mid -m_{0}\cdot\nu (\Theta_{h},F) > 0,
\]
and
\[
c_{F} = \frac{1 +\nu (\Theta_{h},F)}{r +\delta_{0}}+\delta,
\]
where $\delta$ is a sufficiently small positive number.
Then for $b > c_{F}m_{0}+1$
\[
c_{F}(P - \sum_{i}\delta_{i}D_{i}-\sum_{j}\delta_{j}E_{j}) 
+ (b - c_{F}m_{0}-1)(K_{X} -\sum_{i}\nu(\Theta_{h},D_{i})D_{i})
\]
has a singular hermitian metric 
\[
h_{P^{*}}^{c_{F}}\cdot h^{b-c_{F}m_{0}-1}
\]
with strictly positive curvature.
By Nadel's vanishing theorem (Theorem 2.1), we have see that 
the homomorphism
\[
H^{0}(X,{\cal O}_{X}(K_{X}+ F+ (b-1)K_{X})\otimes{\cal I}(h_{P^{*}}^{c_{F}}\cdot h^{b-c_{F}m_{0}-1}))
\rightarrow
\]
\[
H^{0}(F,{\cal O}_{F}(K_{F}+(b-1)K_{X})
\otimes{\cal I}(h_{P^{*}}^{c_{F}}\cdot h^{b-c_{F}m_{0}-1}))
\]
is surjective.
By Theorem 4.1 and its proof, we see that one 
of the followings holds.
\begin{enumerate}
\item 
$H^{0}(F,{\cal O}_{F}(K_{F}+(b-1)K_{X})
\otimes{\cal I}(h_{P^{*}}^{c_{F}}\cdot h^{b-c_{F}m_{0}-1}))\neq 0$
holds for infinitely many  positive integers $b$.
\item 
there exists a nontrivial rational fiber
 space structure 
\[
f : F\longrightarrow W
\] 
such that for a very general fiber $V$, 
$(K_{X},h)$ is numerically trivial. 
\end{enumerate}
We note that {\bf we may not apply Theorem 4.1 directly}, since $h\mid_{F}$ 
is not well defined in general. 
But we shall modify the proof of Theorem 4.1 as follows.

Let us consider the case that $F$ does not admit a nontrivial numerically trivial fibration associated with $(K_{X},h)$. 
 Let $x$ be a very general point on $F$. 
Then adding one more strata, i.e. 
constructing the stratification 
\[
X \supset F \supset F_{1}\supset \cdots \supset F_{r+1} = \{ x\} ,
\]
starting from $X$ (where the strata $F$ is  associated with $h_{P^{*}}^{c_{F}}$) as in the proof of Theorem 4.1,we directly prove the nonvanishing 
\[
H^{0}(X,{\cal O}_{X}(bK_{X}+F)\otimes {\cal I}(h_{P^{*}}^{c_{F}}\cdot h^{b-c_{F}m_{0}-1}))\neq 0
\]
for every sufficiently large $b$.
In this case we construct a singular hermitian metric
$h_{x}$ with strictly positive curvature on $(b-1)K_{X}$ such that 
\begin{enumerate}
\item $\mbox{Spec}({\cal O}_{X}(-[b\nu (\Theta_{h},F)] F)/{\cal I}(h_{x}))$
has isolated support at $x$, 
\item $h_{x}$ is of the form :
\[
h_{x} = h_{P^{*}}^{(1-\varepsilon_{0})c_{F}}\cdot (\prod_{i=0}^{r-1}h_{i}^{(\alpha_{i}-\varepsilon_{i})})
\cdot h_{r}^{\alpha_{r}+\varepsilon_{r}}\cdot 
h_{P^{*}}^{\varepsilon_{0}-(\sum_{i=0}^{r-1}(\alpha_{i}-\varepsilon_{i}))-(\alpha_{r}+\varepsilon_{r}))}\cdot h^{b-(c_{F}m_{0}+1)-\sum_{i=0}^{r}\alpha_{i}\nu_{i}},
\]
where $\{ \nu_{i}\}$ are sufficiently large
positive integers,  $\{ h_{i}\}$ are singular hermitian metrics on 
$\{ \nu_{i}K_{X}+P^{*}\}$ constructed as in the proof of 
Theorem 4.1, $\{ \alpha_{i}\}, \alpha_{i} > 0$ are invariants 
defined as in the proof of Theorem 4.1 such that 
\[
\delta <<\sum_{i=0}^{r}\alpha_{i} < \varepsilon_{0}
\]
holds. 
and $\varepsilon_{i} (0\leq i\leq r)$ are sufficiently small 
positive numbers.
\end{enumerate}  
As in the proof of Theorem 4.1, we may also  need to 
consider the correction term ($e^{-\delta\varphi}$ in the
proof of Theorem 4.1, but  since we are constructing 
a section of $bK_{X}$ with the desired properties, 
this is not essential). 
We note that by the definition of $c_{F}$ 
\[
{\cal O}_{X}(bK_{X}+F)\otimes {\cal I}(h_{P^{*}}^{c_{F}}\cdot h^{b-c_{F}m_{0}-1})
\]
is a subsheaf of 
\[
{\cal O}_{X}(bK_{X}+\lceil A\rceil )\otimes {\cal I}(h^{b}),
\]
where 
\[
A := \sum_{i\neq 0}(-c_{F}r_{i}+\nu (\Theta_{h},D_{i})-\delta_{i})D_{i}.
\] 
Hence by Lemma 2.2,  we see that 
\[
{\cal O}_{X}(bK_{X}+F)\otimes {\cal I}(h_{P^{*}}^{c_{F}}\cdot h^{b-c_{F}m_{0}-1})
\]
is a subsheaf of
${\cal O}_{X}(bK_{X})$ for every sufficiently large $b$ and 
is isomorphic to 
\[
{\cal O}_{F}(bK_{X})\otimes{\cal I}(h^{b})
\,\,\,\,\,\,\mbox{or} \,\,\,\,\,\,\, 
{\cal O}_{F}(bK_{X}+F)\otimes{\cal I}(h^{b})
\]
(in the latter case $\nu (\Theta_{h},F) > 0$ holds) 
on the generic point of $F$. 
By Nadel's vanishing theorem 
\[
H^{1}(X,{\cal O}_{X}(K_{X}+ F+ (b-1)K_{X})\otimes{\cal I}(h_{x})) = 0
\]
holds, we see that for every sufficiently large $b$ there exists a section 
\[
\sigma\in H^{0}(X,{\cal O}_{X}(K_{X}+ F+ (b-1)K_{X})\otimes{\cal I}(h_{P^{*}}^{c_{F}}\cdot h^{b-c_{F}m_{0}-1}))
\]
such that 
\[
\mbox{mult}_{F}(\sigma ) = [b\cdot \nu (\Theta_{h},F)]
\]
holds. By Theorem 3.4, we see that $b\cdot \nu (\Theta_{h},F)$ is 
an integer and by Theorem 4.6,  
$F$ is blown down by $\Phi_{\mid bK_{X}\mid}$ for some $b$.
This contradicts he very generic numerical positivity of $(K_{X},h)$ on $F$.
Hence this case cannot occur and we 
see that $(K_{X},h)$ defines a nontrivial numerically trivial fiber
space structure on $F$.
The last statement follows from Theorem 4.2 and Remark 4.2. 
\vspace{5mm}{\bf  Q.E.D.} \\

The following theorem follows from the proof of Theorem 4.6.

\begin{theorem} 
Let $X,F$ be as in Theorem 4.6.
If for a positive integer $b$ 
\[
\mbox{mult}_{F}\mbox{Bs}\mid bK_{X}\mid = b\cdot\nu (\Theta_{h},F)
\]
holds. 
Then $\Phi_{\mid bK_{X}\mid}$ blows down $F$. 
\end{theorem}
\section{Local base point freeness}
The goal of this section is to prove the following 
proposition. 
\begin{proposition}
Let $X$ be a smooth projective variety of general type. 
Let $u : \tilde{X}\longrightarrow X$ be an arbitrary  composition of successive blowing ups with smooth centers. 
Then for every prime divisor $D$ on $\tilde{X}$,  there exists a positive integer $m(D)$ depending on $D$ such that 
\[
\mbox{mult}_{D}\mbox{Bs}\mid m(D)u^{*}K_{X}\mid  = m(D)\cdot
\nu (u^{*}\Theta_{h},D)
\]
holds. 
\end{proposition} 
Proposition 5.1 and Theorem 4.6 imply that  
every irreducible stable fixed component of $K_{X}$ 
is contracted by $\mid mK_{X}\mid$ for some $m > 0$. 
\subsection{Projective limit of projective varieties}
In this section we deal with a projective limit of 
projective varieties. 
Usually it is not easy to handle such spaces because 
the usual algebro-geometric tools break down on such spaces.
{\bf But in this paper, we only need such spaces to state 
the results. 
The actual proofs are carried out on  usual algebraic 
varieties.} 

First we shall define usuful objects. 
Let 
\[
\cdots \rightarrow M_{m}\stackrel{f_{m}}{\rightarrow} M_{m-1}\stackrel{f_{m-1}}{\rightarrow} \cdots\stackrel{f_{2}}{\rightarrow} M_{1}\stackrel{f_{1}}{\rightarrow} M_{0}:= M
\]
be  successive modifications of a projective variety $M$ such that 
every $M_{m} (m\geq 0)$ is smooth. 
Then we consider the projecive limit 
\[
\hat{M} := \lim_{\leftarrow}M_{m}.
\]
In general $\hat{M}$ is not a projective variety.  

We define the topology on $\hat{M}$ as a minimal topology 
such that  the natural map 
$\hat{M}\longrightarrow M_{m}$ is continuous for every 
$m$. 

Let $m_{0}$ be a nonnegative integer and let $\{ D_{m}\}_{m=m_{0}}^{\infty}$ be a system of divisors such that 
\begin{enumerate}
\item $D_{m}$ is a divisor on $X_{m}$, 
\item $D_{m} = (f_{m+1})_{*}D_{m+1}$ holds for every $m\geq m_{0}$.
\end{enumerate}
We note that $m_{0}$ can be taken to be $0$ by setting
\[
D_{m} := (f_{m_{0},m})_{*}D_{m_{0}}
\]
for every $m\leq m_{0}$, where 
\[
f_{m_{0},m} : M_{m_{0}}\longrightarrow M_{m}
\]
be the composition of $f_{m+1}\circ\cdots \circ f_{m_{0}}$.
In this case we may consider the projective limit 
\[
\hat{D} : = \lim_{\leftarrow} D_{m}
\]
and call it a {\bf divisor} on $\hat{M}$. 
A {\bf prime  divisor} on $\hat{M}$ is a 
projective system of divisor $\{ D_{m}\}$ such that 
every $D_{m}$ is a prime divisor (we consider $0$ is a 
prime divisor).  

Let $D$ be a  divisor of some $M_{m}$. 
For $\ell \geq m$, let $D_{\ell}$ be the strict transform of 
$D$ in $M_{\ell}$. 
Then we see that 
\[
g_{\ell ,*}D_{\ell} = D_{\ell -1}
\]
holds for every $\ell > m$. 
Hence $\{ D_{\ell}\}_{\ell\geq m}$ defines a divisor 
$\hat{D}$ in $\hat{M}$. 
We  call  $\hat{D}$ the {\bf strict transform} of $D$ in 
$\hat{M}$.
We note that every strict transform of a prime divisor 
on some $M_{m}$ is always a prime divisor on $\hat{M}$.

A {\bf sheaf} $\hat{\cal F}$ on $\hat{M}$ is a system of 
sheaves $\{ {\cal F}_{m}\}_{m\geq 0}$ such that 
\[
(f_{m+1})_{*}{\cal F}_{m+1} = {\cal F}_{m}
\]
holds for every $m\geq 0$. 
In particular we can define the structure sheaf ${\cal O}_{\hat{M}}$
is defined as 
\[
{\cal O}_{\hat{M}} : = \lim_{\leftarrow}{\cal O}_{M_{m}}.
\]

Let $h$ be a singular hermitian metric on a line bundle $L$ on $M$ 
such that $\Theta_{h}$ is bounded from below by a $C^{\infty}$-form on 
$M$. 

Let
\[
g_{m} : M_{m}\longrightarrow M
\]
be the natural morphism.
Then we see that 
\[
{\cal O}_{M_{m}}(K_{M_{m}})\otimes{\cal I}(g_{m}^{*}h) 
= (f_{m+1})_{*}({\cal O}_{M_{m+1}}(K_{M_{m+1}})\otimes{\cal I}(g_{m+1}^{*}h))
\]
holds for every $m\geq 0$ by the definition of multiplier ideal sheaves. 
Hence we may define 
\[
{\cal O}_{\hat{M}}(K_{\hat{M}})\otimes {\cal I}(\hat{g}^{*}h)
\]
as the projective limit 
\[
\lim_{\leftarrow}{\cal O}_{M_{m}}(K_{M_{m}})\otimes{\cal I}(g_{m}^{*}h),
\]
where 
\[
\hat{g} : \hat{M}\longrightarrow M
\]
is the natural morphism, i.e, 
the projective limit
\[
\hat{g} : = \lim_{\leftarrow}g_{m}.
\]
We note that ${\cal I}(\hat{g}^{*}h)$ is {\bf not} a well defined sheaf
on $\hat{M}$.
In this paper we always consider 
${\cal O}_{\hat{M}}(K_{\hat{M}})\otimes {\cal I}(\hat{g}^{*}h)$
instead of ${\cal I}(\hat{g}^{*}h)$. 
We call ${\cal O}_{\hat{M}}(K_{\hat{M}})\otimes {\cal I}(\hat{g}^{*}h)$ 
the {\bf multiplier canonical sheaf} of 
$\hat{\pi}^{*}(L,h)$ on $\hat{M}$. 
Also for a prime divisor $\hat{D} = \lim_{\leftarrow}D_{m}$ on $\hat{M}$,
we define the {\bf Lelong number} 
$\nu (\hat{\pi}^{*}\Theta_{h},\hat{D})$ by 
\[
\nu (\hat{g}^{*}\Theta_{h},\hat{D}) 
= \nu (g_{\ell_{0}}^{*}\Theta_{h},D_{\ell_{0}}),
\]
where $\ell_{0}$ is a sufficiently large positive integer.

\subsection{Formal canonical model}
For $m\geq 1$ let 
\[
\pi_{m} : X_{m}\longrightarrow X
\]
be a resolution of $\mbox{Bs}\mid m!K_{X}\mid$ (we set 
$X_{0} = X$).
We may assume 
\begin{enumerate}
\item for $m\geq 2$ there exists a morphism 
\[
\phi_{m} : X_{m}\longrightarrow X_{m-1}
\]
such that 
\[
\pi_{m} = \phi_{m}\circ\pi_{m-1}
\]
holds,
\item $\phi_{m}$ is a sequence of blowing ups with smooth centers contained 
in the indeterminancy locus of the rational map $\Phi_{\pi_{m-1}^{*}\mid m!K_{X}\mid}$.
\item the exceptional divisor of $\pi_{m}$ is a divisor with 
normal crossings,
\item $\pi_{m}^{*}({\cal O}_{X}(m!K_{X})\otimes{\cal I}(h^{m!}))$ 
is locally free on $X_{m}$
\end{enumerate} 
Let $F^{(m)}$ denote the exceptional divisor of $\pi_{m}$.
Let us consider the projective limit: 
\[
\hat{X} = \lim_{\leftarrow} X_{m}.
\]
{\bf  $\hat{X}$ is not a projective variety and depends on the choice of $\{ \pi_{m}\}$ (but we note that we are not considering all such choices at the same time)}.
Let 
\[
\hat{\pi} : \hat{X}\longrightarrow X
\]
be the natural morphism.
We decompose $\phi_{m+1}$ as a sequence of 
blowing ups:
\[
X_{m+1}\stackrel{p^{\ell (m)}_{m}}{\rightarrow}\cdots\stackrel{p^{\ell+2}_{m}}{\rightarrow} X_{m}^{\ell +1}
\stackrel{p^{\ell +1}_{m}}{\rightarrow} X^{\ell}_{m}
\stackrel{p^{\ell}_{m}}{\rightarrow} X^{\ell -1}_{m}
\stackrel{p^{\ell -1}_{m}}{\rightarrow}\cdots\stackrel{p_{m}^{1}}{\rightarrow} X_{m} 
\]
with smooth centers.
Let $\mbox{SE}(K_{X})$ be the subset of $X - \mbox{SBs}(K_{X})$ 
defined by 
\[
\mbox{SE}(K_{X}):= \{ x\in X -\mbox{SBs}(K_{X})\mid \mbox{$\Phi_{\mid m!K_{X}\mid}$ is not local isomorphism onto its image}
\]
\[
\hspace{30mm}
\mbox{ on a
 neighbourhood of $x$ for every 
$m\geq 1$}\}.
\] 
We call $SE(K_{X})$ the {\bf stable exceptional locus} 
of $K_{X}$. 
$\mbox{SE}(K_{X})$ is a divisor in 
$X - \mbox{SBs}(K_{X})$. 
Let $\hat{F}$ be the inverse image of $\mbox{SBs}(K_{X})
\cup \mbox{SE}(K_{X})$.
Then we may and do assume that $\hat{F}$ is a divisor on $\hat{X}$. 
Let 
\[
\hat{F} = \sum \hat{F}_{\alpha}
\]
be the irreducible decomposition of $\hat{F}$. 
For each $\hat{F}_{\alpha}$, there exists 
a rational fibration 
\[
f_{\alpha} : \hat{F}_{\alpha} -\cdots\rightarrow \hat{W}_{\alpha}
\]
constructed as in Theorem 4.7 or the contraction morphism 
on $\mbox{SE}(K_{X})$ induced by 
$\Phi_{\mid m!K_{X}\mid}$ for every suffficiently large $m$. 
We may assume that every $f_{\alpha}$ is a morphism. 
In fact we construct $\{\phi_{m+1}\mid \phi_{m+1} : X_{m+1}\longrightarrow 
X_{m}, m = 0,1,2,\ldots \}$  as follows. 
Let 
\[
F^{(m)} = \sum_{\alpha\in I_{m}}F^{(m)}_{\alpha}
\]
be the irreducible decomposition of $F^{(m)}$. 
Let 
\[
f^{(m)}_{\alpha} : F^{(m)}_{\alpha}\longrightarrow W^{(m)}_{\alpha}
\]
be the fibration constructed as in Theorem 4.7. 
By taking a  composition of successive blowing ups 
\[
w_{m} : \tilde{X}_{m}\longrightarrow X_{m}
\] 
with smooth centers, 
we may assume that  for the strict transform  $\hat{F}^{(m)}_{\alpha}$ 
of $F^{(m)}_{\alpha}$ in $X_{m}$, the induced rational map
\[
\hat{f}^{(m)}_{\alpha} : \hat{F}^{(m)}_{\alpha}-\cdots\rightarrow
W^{(m)}_{\alpha}
\]
is actually a morphism for every $\alpha\in I_{m}$.
We shall take 
\[
\phi_{m+1} : X_{m+1}\longrightarrow X_{m}
\]
so that it factors through $w_{m}$. 
Inductively we repeat the above procedure for all $m\geq 1$. 
Then $\hat{\pi} : \hat{X}\longrightarrow X$ has the desired property. 

We consider the equivalence relation $\sim$ generated by 
$\{ f_{\alpha}\}$, i.e. we identify all the points on a fiber of every $f_{\alpha}$. 
We set the quotient space 
\[
\hat{X}_{can} = \hat{X}/\sim
\]
and call it {\bf  the formal canonical model} of $X$. 
It is easy to see that $\hat{X}_{can}$ does not depend on 
the choice of $\{ \pi_{m}\}$. 
Let 
\[
\varpi : \hat{X}\longrightarrow \hat{X}_{can}
\]
be the natural map. 
The reason why we introduce $\hat{X}_{can}$ is 
that {\bf one may consider $(K_{X},h)$ is numerically positive 
on $\hat{X}_{can}$ as we will see in the next subsection}. 
\subsection{Concentration method on the formal canonical model}
We shall prove Proposition 5.1 in this subsection. 
The proof is similar to that of Theorem 4.1. 
The only difference is that we construct the stratification 
as in Section 5 on the formal canonical model $\hat{X}_{can}$. 
But since we have not proved $\hat{X}_{can}$ is a (projective) variety,
we cannot construct the stratification directly on $\hat{X}_{can}$.
Hence we use the fiber space structure on the stable fixed components.  
Also we use $(K_{X},h)$ ($h$ is the AZD of $K_{X}$ as before) as 
canonical divisor of $\hat{X}_{can}$. 

Let $X$ and $\tilde{X}$ be as in Proposition 5.1. 
It is sufficient to prove the case that $\tilde{X} = X$ holds.
Let $n$ be the dimension of $X$. 
Let $D$ be  a prime divisor on $X$. 
Let $h$ be the analytic Zariski decomposition of $K_{X}$ as before. 
Let 
\[
f_{D} : D -\cdots\rightarrow W
\]
be the rational fibration constructed as in Theorem 4.7.
By  successive blowing ups with smooth centers, we may assume that 
$f_{D}$ is a morphism.  
Let $x$ be a very general point on $D$, 
i.e. $x$ is outside of a union of at most countably
 many proper subvarieties of $D$ .
If we take $x$ very general we may assume that $\hat{\pi}^{-1}(x)\in \hat{X}$ is a point. 
We set 
\[
x_{can} = \varpi (\hat{\pi}^{-1}(x))
\]
and 
\[
\hat{x} = \varpi^{-1}(x_{can}).
\]
Then $\hat{x}$ is a union of at most countably many of subvarieties
in $\hat{X}$.  
We set 
\[
\mu_{0} = \mu (X,K_{X}).
\] 
We note that ${\cal I}(h^{m})$ is locally free at $x$ for every 
$m$, if we take $x\in D$ very general. 
Then since 
\[
\dim H^{0}(X,{\cal O}_{X}(mK_{X})\otimes {\cal I}(h^{m})) 
= \frac{\mu_{0}}{n!}m^{n} + o(m^{n})
\]
holds, we see that for every $\varepsilon > 0$, 
\[
H^{0}(X,{\cal O}_{X}(mK_{X})\otimes {\cal I}(h^{m})\otimes {\cal M}_{x}^{\lceil(1-\varepsilon )\sqrt[n]{\mu_{0}}m\,\rceil }) \neq 0
\]
holds for every sufficiently large $m$. 
Let $m_{0}$ be a sufficiently large positive integer and let 
\[
\sigma_{0}\in H^{0}(X,{\cal O}_{X}(m_{0}K_{X})\otimes {\cal I}(h^{m_{0}})\otimes {\cal M}_{x}^{\lceil(1-\varepsilon )\sqrt[n]{\mu_{0}}m_{0}\rceil })
\]
be a general nonzero element. 
We define the singular hermitian metric $h_{0}$ on $K_{X}$ by 
\[
h_{0} :=  \frac{1}{\mid\sigma_{0}\mid^{\frac{2}{m_{0}}}}.
\]
We set 
\[
\alpha_{0} = \inf \{ \alpha > 0\mid \hat{x}\cap \{ y\in \hat{X}\mid 
({\cal O}_{\hat{X}}(K_{\hat{X}})\otimes{\cal I}(\hat{\pi}^{*}(h^{\alpha}_{0}h^{\beta})))_{y}\subseteq 
\]
\[
\hspace{25mm} {\cal O}_{\hat{X}}(K_{\hat{X}})\otimes
{\cal I}(\hat{\pi}^{*}h^{\alpha +\beta +1})_{y}\otimes {\cal M}_{y}\}\neq \emptyset 
\hspace{5mm} \mbox{holds for every $\beta > 0$} \} .
\]
$\alpha_{0}$ is clearly finite. 
To consider  ${\cal O}_{\hat{X}}(K_{\hat{X}})\otimes{\cal I}(\hat{\pi}^{*}h^{\alpha +\beta +1})$
instead of ${\cal O}_{\hat{X}}(K_{\hat{X}})\otimes{\cal I}(\hat{\pi}^{*}h^{\alpha +\beta})$ 
reflects that we are using $\hat{\pi}^{*}(K_{X},h)$ 
instead of $\hat{\pi}^{*}K_{X}$.
Since $h_{0}$ has algebraic singularities as 
a singular hermitian metric of $K_{X}$, there exists a modification 
\[
p_{0} : Y_{0} \longrightarrow X
\]
such that the current  
$(\alpha_{0}p_{0}^{*}\Theta_{h_{0}})_{sing} (=
\alpha_{0}p_{0}^{*}\Theta_{h_{0}})$ is a divisor with normal crossings 
$B = \sum b_{i}B_{i}$.  
Then  if we define the numbers  $\{ c_{i}\}$ by 
\[
K_{Y_{0}} - p_{0}^{*}K_{X} - B = \sum_{i}c_{i}B_{i},
\] 
\[
\min \{\nu (p_{0}^{*}\Theta_{h},B_{i}) + c_{i}\mid \hat{x}\cap \hat{B}_{i}\neq \emptyset \} = -1
\]
holds, where  we have assumed that $\hat{\pi}$ factors through 
$p_{0}$ (this is clearly possible)  and $\hat{B}_{i}$ denotes the strict transform of $B_{i}$ in $\hat{X}$.  
We note that if we replace $p_{0}$ by another $p_{0}^{\prime}$ 
which factors through $p_{0}$, 
then by Corollary 2.1 the prime divisors which attain the above minimum 
are exactly the strict transforms of the ones 
associated with $p_{0}$.

By the above assumption there exists a morphism 
\[
q_{0} : \hat{X}\longrightarrow Y_{0}. 
\]
Since $p_{0}^{*}K_{X}$ is big, by Kodaira's lemma, there exists an effective 
{\bf  Q}-divisor $E_{0}$ on $Y_{0}$ such that 
$p_{0}^{*}K_{X} - E_{0}$ is ample.
Let $h_{E_{0}}$ be a $C^{\infty}$-hermitian metric on $p_{0}^{*}K_{X}-E_{0}$ 
(this is a {\bf  Q}-divisor on $Y_{0}$, but the hermitian metric is well defined) with strictly positive curvature. 
We may and do consider $h_{E_{0}}$ a singular hermitian metric on $p_{0}^{*}K_{X}$. 
If we perturb $h_{0}$ as 
\[
h_{0}:= (\frac{1}{\mid\sigma_{0}\mid^{\frac{2}{m_{0}}}})^{1-\delta_{0}}\cdot h_{E_{0}}^{\delta_{0}},
\]
where $\delta_{0}$ is a sufficiently small positive number, perturbing also 
$E_{0}$, if necessary,
we may assume that there exists a unique irreducible divisor $D_{1}=B_{i_{0}}$ 
which belongs to  
$\{ B_{i}\mid \hat{x}\cap \hat{B}_{i}\neq \emptyset \}$ such that 
\[
\nu (p_{0}^{*}\Theta_{h},D_{1})+c_{i_{0}} = -1
\]
holds.  
We set  
\[
Z_{1} : = \varpi (\hat{D}_{1}),
\]
where $\hat{D}_{1}$ is the strict transform of 
$D_{1}$ in $\hat{X}$.
We define the nonnegative integer  $n_{1}$ by 
\[
n_{1} := \dim \varpi (\hat{D}_{1}).
\]   
We note that $n_{1}$ is nothing but the dimension of the base space 
of the numerically trivial fiber space structure on $D_{1}$ 
associated with $p_{0}^{*}(K_{X},h)$. 
If $n_{1}$ is $0$, then $p_{0}^{*}(K_{X},h)$ 
is numerically trivial on $D_{1}$ (cf. Section 3.6 and 
Lemma 4.1, also \cite[Theorem 4.1]{tu4}). 
In this case we stop this process. 

Suppose that $n_{1} > 0$ holds.
We set 
\[
A_{1} := r_{1}(p_{0}^{*}K_{X}-E_{0})\mid_{D_{1}},
\]
where $r_{1}$ is a sufficiently large positive integer such that 
$r_{1}(p_{0}^{*}K_{X}-E_{0})$ is Cartier. 
We set 
\[
\mu_{1} := (n-1)!\cdot\overline{\lim}_{m\rightarrow \infty}m^{-(n-1)}
\dim H^{0}(D_{1}, {\cal O}_{D_{1}}(m(A_{1}+ p_{0}^{*}(\ell_{1}K_{X}))\otimes {\cal I}(p_{0}^{*}h^{\ell_{1}m})),
\] 
where $\ell_{1}$ is a sufficiently large positive integer which will be 
specified later. 
Let $y_{1}\in D_{1}\cap q_{0}(\varpi^{-1}(x_{can}))$ be a  point.
And we set  $x_{1} = p_{0}(y_{1}) \in X_{1}$. 
Then as before 
\[
H^{0}(D_{1},{\cal O}_{D_{1}}(m(A_{1}+ p_{0}^{*}(\ell_{1}K_{X})))\otimes {\cal I}(p_{0}^{*}h^{m})\otimes 
{\cal M}_{y_{1}}^{\lceil(1-\varepsilon )\sqrt[n-1]{\mu_{1}}m
\rceil }) \neq 0
\]
holds for every sufficiently large $m$. 
Let $m_{1} >> r_{1}$ be a sufficiently large positive integer and 
let 
\[
\sigma_{1}^{\prime}\in H^{0}(D_{1},{\cal O}_{D_{1}}(m_{1}(A_{1}+ p_{0}^{*}(\ell_{1}K_{X}))
\otimes {\cal I}(p_{0}^{*}h^{\ell_{1}m_{1}})\otimes {\cal M}_{y_{1}}^{\lceil(1-\varepsilon )\sqrt[n-1]{\mu_{1}}m_{1}\rceil })
\]
be a  general nonzero element.  
We note that since $D_{1}$ is smooth,   
\[
{\cal O}_{D_{1}}(p_{0}^{*}(m_{1}\ell_{1}K_{X}))\otimes {\cal I}(p_{0}^{*}h^{m_{1}\ell_{1}})
\]
is torsion free, since it is a subsheaf of a locally free sheaf on a smooth variety.
Then as in Lemma 4.4, we see that 
the restriction map 
\[
H^{0}(Y_{0},{\cal O}_{Y_{0}}(m(A_{1}+ p_{0}^{*}(\ell_{1}K_{X})))\otimes {\cal I}(p_{0}^{*}h^{m\ell_{1}}))\rightarrow 
\]
\[
\hspace{30mm} 
H^{0}(D_{1},{\cal O}_{D_{1}}(m(A_{1}+ p_{0}^{*}(\ell_{1}K_{X})))\otimes {\cal I}(p_{0}^{*}h^{m\ell_{1}})) 
\]
is surjective for every $m > 0$, if we take $r_{1}$ sufficiently large. 
Then 
$\sigma_{1}^{\prime}$
extends to an element $\sigma_{1}$ of 
\[
H^{0}(Y_{0},{\cal O}_{Y_{0}}(m_{1}(A_{1}+ p_{0}^{*}(\ell_{1}K_{X})))\otimes {\cal I}(p_{0}^{*}h^{\ell_{1}m_{1}})). 
\]
We define the singular hermitian metric $h_{1}$ of $(r_{1}+\ell_{1})K_{X}$ by 
\[
h_{1} = \frac{1}{\mid \sigma_{1}\mid^{\frac{2}{m_{1}}}}
\]
(Originally $h_{1}$ is considered to be a singular hermitian metric on $p_{0}^{*}(r_{1}+\ell_{1})K_{X}$, but one may consider $h_{1}$ as a singular hermitian metric on $(r_{1}+\ell_{1})K_{X}$).
Let $\varepsilon_{0}$ be a sufficiently small positive number. 
We define the positive number $\alpha_{1}$ by 
\[
\alpha_{1} := \inf\{ \alpha \mid \hat{x}\cap \{ y\in \hat{X}\mid 
({\cal O}_{\hat{X}}(K_{\hat{X}})\otimes{\cal I}(\hat{\pi}^{*}(h_{0}^{\alpha_{0}-\varepsilon_{0}}\cdot h_{1}^{\alpha}\cdot h^{\beta})))_{y}\subseteq 
\]
\[
\hspace{25mm}
{\cal O}_{X}(K_{\hat{X}})\otimes{\cal I}(\hat{\pi}^{*}h^{\alpha_{0}-\varepsilon_{0}+\alpha+\beta+1})_{y}\otimes{\cal M}_{y}\}
\neq\emptyset 
\hspace{5mm} \mbox{holds for every $\beta > 0$}\} \}.
\]
Then as in Section 5 we have the estimate :
\[
\alpha_{1}\leq \frac{n-1}{\sqrt[n_{1}]{\mu_{1}}} + O (\varepsilon_{0}).
\]
Taking $\ell_{1}$ to be sufficiently large, we may assume that 
\[
\mu_{1} >>  (\frac{(n-1)r_{1}\cdot\mbox{mult}_{D_{1}}E_{0}}{\varepsilon_{0}})^{n-1}
\]
holds.

We take a modification 
\[
f_{1} : Y_{1}\longrightarrow Y_{0}
\]
such that the singular part of $f_{1}^{*}p_{0}^{*}((\alpha_{0}-\varepsilon_{0})\Theta_{h_{0}}+\alpha_{1}\Theta_{h_{1}}))$ 
is a divisor with normal crossings in $Y_{1}$.
We may assume that $\hat{\pi} : \hat{X}\longrightarrow X$ 
factors through $p_{1}:= f_{1}\circ p_{0}$. 
Then by the procedure as before, we define a divisor $D_{2}$ as before in $Y_{1}$ and the subset $Z_{2}$ in $\hat{X}_{can}$ by 
\[
Z_{2} = \varpi (\hat{D}_{2}),
\]
where $\hat{D}_{2}$ is the strict transform of $D_{2}$ 
in $\hat{X}$. 
We note that since we have taken $\ell_{1}$ so that 
\[
\mu_{1} >>  (\frac{(n-1)r_{1}\cdot\mbox{mult}_{D_{1}}E_{0}}{\varepsilon_{0}})^{n-1}
\]
holds,
 by  the estimate of $\alpha_{1}$,
we have that 
\[
\alpha_{1}r_{1}\cdot\mbox{mult}_{D_{1}} E_{0} << \varepsilon_{0}
\]
holds. 
Hence the singularity  of $h_{1}^{\alpha_{1}}$ coming from  
$A_{1} = r_{1}(p_{0}^{*}K_{X} - E_{0})$ (roughly speaking the singularity is equal to  $\alpha_{1}r_{1}E_{0}$) is enough small so that   
 $Z_{2}$ is a proper subset of $Z_{1}$.
Inductively we define 
a sequence of modifications 
\[
X \leftarrow Y_{0} \leftarrow Y_{1}\leftarrow\cdots \leftarrow Y_{r},
\]
irreducible smooth divisors
\[
D_{i}\subset Y_{i-1} (i=1,\ldots ,r+1),
\]
points 
\[
y_{i} \in D_{i}\cap q_{i-1}(\varpi^{-1}(x_{can}))\,\,
(i=1,\ldots ,r),
\]
where 
\[
q_{i} : \hat{X}\longrightarrow Y_{i}
\]
are the natural morphisms,
singular hermitian metrics 
\[
h_{0},\ldots ,h_{r},
\]
small positive numbers 
\[
\varepsilon_{0},\ldots ,\varepsilon_{r-1},
\]
positive integers
\[
\ell_{1},\ldots ,\ell_{r},
\]
positive numbers 
\[
\alpha_{0},\ldots ,\alpha_{r},
\]
and nonnegative integers 
\[
n_{1},\ldots ,n_{r+1},
\]
positive numbers
\[
\mu_{0},\ldots \mu_{r+1}
\]
and strictly decreasing sequence of irreducible subsets 
\[
Z_{1}\supset Z_{2}\supset\cdots  \supset Z_{r+1}
\]
in $\hat{X}_{can}$. 
By the construction of this process we see that
$n_{r+1}=0$ holds. 
This means that 
\[
Z_{r+1} = x_{can}
\]
holds. 
We define the singular hermitian metric $h_{x}$ on 
 \\ 
$(\sum_{i=1}^{r-1}(\alpha_{i}-\varepsilon_{i})+(\alpha_{r}+\varepsilon_{r}))K_{X}$ 
by 
\[
h_{x} := h_{0}^{\alpha_{0}-\varepsilon_{0}}\cdot h_{1}^{\alpha_{1}-\varepsilon_{1}}\cdots h_{r}^{\alpha_{r}+\varepsilon_{r}}, 
\] 
where $\varepsilon_{r}$ is a sufficiently small positive number. 
We set for  every positive number  $\beta > \sum_{i=0}^{r} \alpha_{i}(r_{i}+\ell_{i})$ (where we have set $\ell_{0} := 1$ and $r_{0} := 0$),
\[
h_{x}(\beta ) = h^{\beta-\sum_{i=1}^{r-1}(\alpha_{i}-\varepsilon_{i})-(\alpha_{r}+\varepsilon_{r})}\cdot h_{x}. 
\]
We see  that $p_{r}^{*}(K_{X},h)$ is numerically 
trivial on $D_{r+1}$. 
Hence $\hat{D}_{r+1}$ is  contained in $\hat{x}$, where $\hat{D}_{r+1}$ is the 
strict transform of $D_{r+1}$ in $\hat{X}$.  

Let $p_{i} : Y_{i}\longrightarrow X (i=0,\ldots ,r)$ be the natural morphisms. 
\begin{lemma}
\[
H^{0}(Y_{r},{\cal O}_{Y_{r}}(K_{Y_{r}}+ mp_{r}^{*}K_{X}+D_{r+1})\otimes
{\cal I}(p_{r}^{*}h_{x}(m)))
\rightarrow 
\]
\[
\hspace{30mm} 
H^{0}(D_{r+1},{\cal O}_{D_{r+1}}(K_{D_{r+1}}+mp_{r}^{*}K_{X})\otimes
{\cal I}(p_{r}^{*}h_{x}(m)))
\]
is surjective for every positive integer $m \geq
\sum_{i=1}^{r-1}(\alpha_{i}-\varepsilon_{i})-(\alpha_{r}+\varepsilon_{r}) $.
\end{lemma}
{\bf Proof}.
The assertion follows from Theorem 2.1, since 
\[
0\rightarrow {\cal O}_{Y_{r}}(K_{Y_{r}}+ mp_{r}^{*}K_{X})\otimes
{\cal I}(p_{r}^{*}h_{x}(m))
\rightarrow 
{\cal O}_{Y_{r}}(K_{Y_{r}}+ mp_{r}^{*}K_{X}+D_{r+1})\otimes
{\cal I}(p_{r}^{*}h_{x}(m))
\rightarrow 
\]
\[
\hspace{50mm}
{\cal O}_{D_{r+1}}(K_{D_{r+1}}+mp_{r}^{*}K_{X})\otimes
{\cal I}(p_{r}^{*}h_{x}(m))\rightarrow 0
\]
is exact.
{\bf Q.E.D.} 
\subsection{Finding divisors on $D_{r+1}$}
On $D_{r+1}$, $p_{r}^{*}(K_{X},h)$ is numerically 
trivial. 
By Theorem 4.3 and Corollary 4.2 we  see that 
\[
S : = \{ x\in D_{r+1}\mid \nu_{D_{r+1}}(p_{r}^{*}\Theta_{h},x) > 0\}
\]
consists of at most {\bf countably many} prime divisors on $D_{r+1}$.
Let 
\[
S = \sum_{j\in J}E_{j}
\]
be the irredcucible decomposition of $S$.
We set 
\[
e_{j} := \nu_{D_{r+1}}(p_{r}^{*}\Theta_{h},E_{j})\hspace{15mm} (j\in J).
\]
By Theorem 4.3 and Corollary 4.2 we have the following 
lemma. 
\begin{lemma}
\[
(p_{r}^{*}K_{X} - \nu (p_{r}^{*}\Theta_{h},D_{r+1})D_{r+1})\mid_{D_{r+1}}- \sum_{j\in J}e_{j}E_{j}
\]
is numerically trivial on $D_{r+1}$. 
\end{lemma}

\subsection{Completion of the proof of Proposition 5.1.}
Let  
\[
2\pi (p_{r}^{*}\Theta_{h_{x}})_{sing} 
= \sum_{i\in I}r_{i}F_{i} 
\]
be the decomposition into irreducible component of 
the singular part of the 
current $2\pi p_{r}^{*}\Theta_{h_{x}}$. 
We may assume that $\sum_{i}F_{i}$ is a divisor with 
normal crossings. 
Changing the indices if necessary, we may assume that $F_{0} = D_{r+1}$. 
We define the {\bf R}-divisor $A_{r}^{\prime}$ on $Y_{r}$ by 
\[
A_{r}^{\prime}: = K_{Y_{r}} + \beta\cdot p_{r}^{*}K_{X} -\sum_{i\in I} r_{i}F_{i} + F_{0},
\]
where $\beta$ is the positive number defined by 
\[
\beta := \sum_{k=0}^{r-1}(\alpha_{k}-\varepsilon_{k})
+ (\alpha_{r}+\varepsilon_{r}). 
\]
We set 
\[
\nu_{i} : = \nu (p_{r}^{*}\Theta_{h},F_{i}) \hspace{5mm} (i \in I)
\]
and 
\[
K_{Y_{r}} = p_{r}^{*}K_{X} + \sum_{i\in I}a_{i}F_{i}.
\]
Then by the definition of $h_{x}$ we see that 
\[
\min_{i\in I}(-r_{i}+a_{i}+\nu_{i}) = 
= -r_{0}+a_{0}+\nu_{0} = -1
\]
and 
\[
\min_{i\neq 0}(-r_{i}+a_{i}+\nu_{i}) > -1
\]
hold.  
If we set 
\[
\theta_{i}:= -r_{i}+a_{i}+\nu_{i}
\]
for every $i\in I$, 
it is easy to verify that $A_{r}^{\prime}$ is numerically 
equivalent to the {\bf R}-divisor $A_{r}$ defined by 
\[
A_{r} = 
(\beta + 1)(p_{r}^{*}K_{X}-\sum_{i}\nu_{i}F_{i})
+ \sum_{i\neq 0}\theta_{i}F_{i}.
\]
We note that since $\theta_{i} > -1$ for every $i\neq 0$, 
we see that $\lceil\sum_{i\neq 0}\theta_{i}F_{i}\rceil$ 
is effective. 
Since $p_{r}^{*}\Theta_{h_{x}}$ is strictly positive on 
$Y_{r}$, we see that 
\[
A_{r}\mid_{F_{0}} - K_{F_{0}} 
\]
is ample.
We note that by Lemma 5.2
\[
(p_{r}^{*}K_{X} -\nu_{0}F_{0})\mid_{F_{0}}
-\sum_{j\in J}e_{j}E_{j}
\]
is numerically trivial on $F_{0} (= D_{r+1})$.
Let $J_{0}$ be the subset of the indices $J$ such that 
\[
\sum_{i\neq 0} F_{i}\mid_{F_{0}} = \sum_{j\in J_{0}}E_{j}
\]
We set 
\[
E^{*} := \sum_{j \in J-J_{0}}e_{j}E_{j}.
\]
Then $E^{*}$ defines a point on the closure of the 
cone of the effective {\bf R}-divisors on $F_{0}$ in $H^{2}(F_{0},\mbox{\bf R})$.
Let $H_{0}$ be a smooth very ample divisor on $F_{0}$ such that 
\[
\sum_{i\neq 0}F_{i}\mid_{F_{0}}
+ H_{0}
\]
is a divisor with normal crossings. 
Let $\epsilon$ be a sufficiently small positive number 
such that 
\[
A_{r}\mid_{F_{0}} - K_{F_{0}}- \epsilon H_{0} 
\]
is ample. 
Since 
\[
\epsilon H + E^{*}
\]
is numerically equivalent to an effective {\bf R}-divisor 
in $H^{2}(F_{0},\mbox{\bf R})$, 
taking a suitable modification, if necessary, we may assume that 
$A_{r}\mid_{F_{0}}$ is numerically equivalent to an 
 {\bf R}-divisor $B_{r}$ with normal crossings 
such that $\lceil B_{r}\rceil$ is effective. 

We note that 
\[
{\cal O}_{Y_{r}}(K_{Y_{r}}+m\cdot p_{r}^{*}K_{X})
\otimes {\cal I}(p_{r}^{*}h_{x}(m))
\simeq {\cal O}_{Y_{r}}(\lceil A_{r}\rceil+(m-\beta )p_{r}^{*}K_{X})
\otimes {\cal I}(h_{\lceil A_{r}\rceil}\cdot p_{r}^{*}h^{m-\beta})
\]
holds for every $m \geq \beta$, where the singular hermitian metric 
 $h_{\lceil A_{r}\rceil}$  on ${\cal O}_{Y_{r}}(\lceil A_{r}\rceil )$ 
is 
defined as the singular hermitian metric 
$h_{\lceil A\rceil}$ in Theorem 4.4. 
We note that 
there exists a positive number $\varepsilon_{0}$ 
such that if $\{ m\cdot \nu (p_{r}^{*}\Theta_{h},D_{r+1})\}
 \leq \epsilon_{0}$ holds,  
\[
{\cal O}_{D_{r+1}}(\lceil A_{r}\rceil + (m-\beta )p_{r}^{*}K_{X})
\otimes {\cal I}(\mid \tau_{r+1}\mid^{2\{ m\cdot\nu (p_{r}^{*}\Theta_{h},D_{r+1})\}}\cdot h_{\lceil B_{r}\rceil}\cdot h^{m-\beta })
\subseteq 
\]
\[
\hspace{30mm}
{\cal O}_{D_{r+1}}(\lceil A_{r}\rceil+(m-\beta )p_{r}^{*}K_{X})
\otimes {\cal I}(h_{\lceil A_{r}\rceil}\cdot p_{r}^{*}h^{m-\beta})
\]
holds by the perturbation of $A_{r}$ as in the proof of 
Theorem 4.4, where $\tau_{r+1}$ is a global section of 
${\cal O}_{Y_{r}}(D_{r+1})$ with divisor $D_{r+1}$ 
and $h_{\lceil B_{r}\rceil}$ is the singular hermitian metric 
on ${\cal O}_{D_{r+1}}(\lceil B_{r}\rceil )$ defined 
as the singular hermitian metric 
$h_{\lceil A\rceil}$ in Theorem 4.4.  
We also note that 
\[
{\cal O}_{D_{r+1}}(\lceil A_{r}\rceil+(m-\beta )p_{r}^{*}K_{X})
\otimes {\cal I}(h_{\lceil A_{r}\rceil}\cdot p_{r}^{*}h^{m-\beta})
\simeq 
{\cal O}_{D_{r+1}}(K_{D_{r+1}}+mp_{r}^{*}K_{X})\otimes
{\cal I}(p_{r}^{*}h_{x}(m))
\]
holds by the definition of $A_{r}$. 
Now we apply Theorem 4.4 and Remark 4.7 to our situation by 
setting 
$X = D_{r+1}$, $A = B_{r}$  and $L= p_{r}^{*}(K_{X},h)$. 
Then we have the following lemma.
\begin{lemma} 
For some $m > \sum_{i=0}^{r}\alpha_{i}(r_{i}+\ell_{i})$,  
\[
H^{0}(D_{r+1},{\cal O}_{D_{r+1}}(K_{D_{r+1}}+mp_{r}^{*}K_{X})\otimes
{\cal I}(p_{r}^{*}h_{x}(m)))\neq 0
\]
holds.
\end{lemma}
By Lemma 5.1 and Lemma 5.3, there exists a positive integer $m(D_{r+1})$ such that 
\[
\mbox{mult}_{D_{r+1}}\mbox{Bs}\mid p_{r}^{*}(m(D_{r+1})K_{X})\mid 
= m(D_{r+1})\cdot \nu(p_{r}^{*}\Theta_{h},D_{r+1})
\]
holds. 
Let 
\[
\hat{F} = \sum_{\alpha}\hat{F}_{\alpha} 
\]
be the inverse image of $\hat{\pi}^{-1}(\mbox{SBs}(K_{X}))$
as in Section 5.2. 
And let 
\[
\varphi_{m} : \hat{X}\longrightarrow X_{m}
\]
be the natural map. 
We note that  {\bf since $\hat{\pi}^{*}(K_{X},h)$ is numerically trivial on 
$\hat{x}$},
the divisor 
\[
\hat{\pi}^{*}K_{X} - \sum_{\alpha}\nu (\hat{\pi}^{*}\Theta_{h},\hat{F}_{\alpha})\hat{F}_{\alpha}
\]
is numerically trivial on $\hat{x}$, 
i.e., 
for every $m\geq 0$, $\pi_{m}^{*}(K_{X},h)$ is 
numerically trivial on every irreducible component of 
$\varphi_{m}(\hat{x})$. 
By Theorem 3.4, we see that
for any positive integer $\ell$ and 
any $\sigma\in \Gamma (X,{\cal O}_{X}(\ell K_{X}))$, 
\[
\mbox{mult}_{\hat{F}_{\alpha}}\hat{\pi}^{*}(\sigma )
= \mbox{mult}_{F_{\alpha ,m}}\pi_{m}^{*}(\sigma )
\geq \lceil \ell\cdot\nu (\pi_{m}^{*}\Theta_{h},F_{\alpha ,m})\rceil
= \lceil \ell\cdot\nu (\hat{\pi}^{*}\Theta_{h},\hat{F}_{\alpha})\rceil
\]
hold, where $m$ is a sufficiently large positive integer 
depending on $\alpha$ and 
$F_{\alpha ,m}$ is the prime divisor whose strict transform 
in $\hat{X}$ is $\hat{F}_{\alpha}$.

Let $\sum_{j\in J}e_{j}E_{j}$ be the divisor on $D_{r+1}$ as in Lemma 5.2. 
Then the above formula implies that for every $j\in J$, $m(D_{r+1})e_{j}$ is an integer. 
Hence $J$ is a finite set and 
\[
m(D_{r+1})(p_{r}^{*}K_{X}\mid_{D_{r+1}}- \sum_{j\in J}e_{j}E_{j})
\]
is a  Cartier divisor on  $D_{r+1}$ and is linearly equivalent to $0$. 
Let us consider $X_{m(D_{r+1})}$ and let 
$D_{r+1}^{*}$ be the divisor on $X_{m(D_{r+1})}$ defined by 
\[
D_{r+1}^{*} := (\varphi_{m(D_{r+1})})_{*}\hat{D}_{r+1},
\]
where $\hat{D}_{r+1}$ is the strict transform of 
$D_{r+1}$ in $\hat{X}$.   
Since $m(D_{r+1})(p_{r}^{*}K_{X}\mid_{D_{r+1}}- \sum_{j\in J}e_{j}E_{j})$ 
is linearly equivalent to $0$, 
we see that 
\[
\mbox{Bs}\mid \pi_{m(D_{r+1})}^{*}m(D_{r+1})(K_{X},h)\mid 
\cap D_{r+1}^{*} = \emptyset 
\]
holds. 
Hence for any subvariety $V$ on $X_{m(D_{r+1})}$ such that $D_{r+1}^{*}\cap V \neq \emptyset$, 
\[
\mbox{mult}_{V}\mid \pi_{m(D_{r+1})}^{*}(m(D_{r+1}K_{X})\mid 
= m(D_{r+1})\cdot \nu (\pi_{m(D_{r+1})}^{*}\Theta_{h},V)
\]
holds. 
Let us define the analytic subset $V_{m(D_{r+1})}$ in $X_{m(D_{r+1})}$ by 
\[
V_{m(D_{r+1})}:= \varphi_{m(D_{r+1})}(\hat{x}).
\]
Then $V_{m(D_{r+1})}$ is connected and contains $D_{r+1}^{*}$.
\begin{lemma} $\pi_{m(D_{r+1})}^{*}(K_{X},h)$ is numerically trivial 
on $V_{m(D_{r+1})}$. 
\end{lemma}
{\bf Proof of Lemma 5.4}. 
Suppose the contrary. 
Then $\varphi_{m(D_{r+1})}^{*}(\pi_{m(D_{r+1})}^{*}(K_{X},h))\mid_{\hat{x}}$ 
is {\bf not} numerically trivial. 

On the other hand  by the definition $\hat{\pi}^{*}(K_{X},h)$ 
is numerically trivial on $\hat{x}$, i.e., 
for every $m\geq 0$, $\pi_{m}^{*}(K_{X},h)$ is 
numerically trivial on every irreducible component of 
$\varphi_{m}(\hat{x})$. 

This is the contradiction. \vspace{5mm} {\bf Q.E.D.} \\
By Lemma 5.4, we see that 
\[
\mbox{Bs}\mid\pi_{m(D_{r+1})}^{*}(K_{X},h)\mid\cap V_{m(D_{r+1})}
= \emptyset
\]
holds. 
Thus we see that 
\[
[m(D_{r+1})(\hat{\pi}^{*}K_{X} - \sum_{\alpha}\nu (\hat{\pi}^{*}\Theta_{h},\hat{F}_{\alpha})\hat{F}_{\alpha})]\mid_{\hat{x}}
= \lceil m(D_{r+1})(\hat{\pi}^{*}K_{X} - \sum_{\alpha}\nu (\hat{\pi}^{*}\Theta_{h},\hat{F}_{\alpha})\hat{F}_{\alpha})\rceil\mid_{\hat{x}}
\]
and  
\[
\mbox{Bs}\mid\hat{\pi}^{*}(m(D_{r+1})(K_{X},h))\mid 
\cap \,\hat{x} = \emptyset
\]
hold. 
Hence {\bf the base point freeness propagates through  
$\hat{x}$} (in particular we may take $\hat{X}$ so that $\hat{x}$ consists of 
finitely many irreducible components).
Since $x$ is a very general point on $D$, we see that there exists a positive integer $m(D)$  such that 
\[
\mbox{mult}_{D}\mbox{Bs}\mid m(D)K_{X}\mid = m(D)\cdot \nu (\Theta_{h},D)
\]
holds. 
This completes the proof of Proposition 5.1. {\bf  Q.E.D.}
\section{Completion of the proof of Theorem 1.1}
We complete the proof of Theorem 1.1 by using a topological consideration. 
We use the same notations and conventions as in Section 5.
\begin{definition}
Let $X$ be a smooth projective variety of dimension $n$ and let 
$L$ be a big line bundle on $X$. 
Let $R = \oplus_{m\geq 0}R_{m}$ be a subring of $R(X,L)$ 
such that 
\[
\lim_{m\rightarrow \infty}{m^{-n}}\dim R_{m} > 0.
\]
For every subvariety $V$ in $X$, we set  
\[
\nu (R,V):= \lim_{m\rightarrow\infty}\frac{1}{m}\mbox{mult}_{V}\mbox{Bs}\mid R_{m}\mid ,
\]
where $\mbox{Bs}\mid R_{m}\mid$ is the base scheme as a linear subsystem 
of $\mid mL\mid$. 
Suppose that for every modification 
\[
f : Y\longrightarrow X
\]
and every prime divisor $D$ on $Y$, there exists a 
positive integer $m_{D}$ depending on $D$ such that 
\[
\nu (f^{*}R,D) = \frac{1}{m_{D}}\mbox{mult}_{D}
\mbox{Bs}\mid m_{D}f^{*}R_{m_{D}}\mid 
\]
holds. 

In this case we call that $R$ is  virtually base point free on $X$. 
\end{definition}
Proposition 5.1 implies that the canonical ring $R(X,K_{X})$
of smooth projective variety of general type $X$ is virtually 
base point free.

Let $X$ be a smooth projective variety of general type and let 
$n$ denote the dimension of $X$.
By the virtual base point freeness of $R(X,K_{X})$, 
we see that $\hat{X}_{can}$ is a complex space (possibly noncompact).
In fact by the construction and Proposition 5.1, for every compact subset 
$G$ of $\hat{X}_{can}$, there exists a positive integer $m(G)$ depending on $G$ such that 
$m(G)K_{\hat{X}_{can}}$ is Cartier on $G$ and $\mid m(G)K_{\hat{X}_{can}}\mid$ is base point free 
on every compact subset of $\hat{X}_{can}$ and is numerically positive on 
$\hat{W}\cap G$ in the obvious sense. 
This implies that $\hat{X}_{can}$ is a complex space.  
Also it is easy to see that $\hat{X}_{can}$ is normal by showing 
that $\hat{X}_{can}$ is isomorphic to the normalization. 
Moreover since by the construction of $\hat{X}_{can}$ and the virtual base point freeness, $\oplus_{m\geq 0}{\cal O}_{\hat{X}_{can}}(mK_{\hat{X}_{can}})$ 
is a finitely generated ring over ${\cal O}_{\hat{X}_{can}}$ on $G$,
 $\hat{X}_{can}$ has only canonical singularities.

Let $\hat{W}$ be the subspace of $\hat{X}_{can}$ defined by 
\[
\hat{W} : = \hat{\varpi}(\hat{F}). 
\]
Then by Theorem 4.6 $\mbox{codim}\,\hat{W} \geq 2$ holds.
We only need to consider the case : $\dim X \geq 3$. 
Now we consider the exact sequence : 
\[
H^{2}(\hat{X}_{can},\mbox{\bf Z}) \rightarrow H^{2}(\hat{W},\mbox{\bf Z})\rightarrow H^{3}(\hat{X}_{can},\hat{W},\mbox{\bf Z}).
\]
We note that since $\mbox{codim}\,\hat{W} \geq 2$, i.e., 
$\hat{\varpi}$ contracts all the irreducible components 
of $\hat{F}$ in $\hat{X}$,  
$\dim H^{2}(\hat{X}_{can},\mbox{\bf C})$ is finite.
Hence we see that 
\[
\mbox{rank}\, H^{3}(\hat{X}_{can},\hat{W},\mbox{\bf Z}) = \infty
\]
holds, if 
\[
\mbox{rank} \,\dim H^{2}(\hat{W},\mbox{\bf Z}) = \infty
\]
holds.
We note that 
\[
 H^{3}(\hat{X}_{can},\hat{W},\mbox{\bf Z})\simeq H^{3}(X,S,\mbox{\bf Z})
\]
holds, where $S$ denotes the union of the stable base locus $\mbox{SBs}(K_{X})$ and the stable exceptional locus
$\mbox{SE}(K_{X})$. 
This means that 
$H^{3}(\hat{X}_{can},\hat{W},\mbox{\bf C})$
is an finitely generated abelian group. 
This implies that
\[
\mbox{rank}\, H^{2}(\hat{W},\mbox{\bf Z}) < \infty
\]
holds (in the case of $\dim X = 3$, this immediately implies that 
$\hat{W}$ consists of finitely many irreducible components).
By the universal coefficients theorem we see that
\[
0\rightarrow \mbox{Ext}(H_{1}(\hat{X}_{can},\mbox{\bf Z}),\mbox{\bf Z})
\rightarrow H^{2}(\hat{X}_{can},\mbox{\bf Z})
\rightarrow \mbox{Hom}(H_{2}(\hat{X}_{can},\mbox{\bf Z}),\mbox{\bf Z})
\rightarrow 0
\]
is  exact.   
Since $H_{1}(\hat{X}_{can},\mbox{\bf Z})$ is finitely generated (because
$\mbox{codim}\hat{W}\geq 2$ holds), we see that 
the torsion part of  $H^{2}(\hat{X}_{can},\mbox{\bf Z})$ is finite.

Since $H^{2}(\hat{X}_{can},\mbox{\bf Z})$ is finitely generated 
and $\mbox{rank}\,H^{2}(\hat{W},\mbox{\bf Z})$ is finite, 
considering the maps:
\[
H^{2}(\hat{X}_{can},\mbox{\bf Z})
\rightarrow \mbox{Hom}(H_{2}(\hat{X}_{can},\mbox{\bf Z}),\mbox{\bf Z})
\rightarrow \mbox{Hom}(H_{2}(\hat{W},\mbox{\bf Z}),\mbox{\bf Z}),
\]
we see that the images  of some positive multiple of $c_{1}(K_{\hat{X}_{can}})
\in H^{2}(\hat{X}_{can},\mbox{\bf R})$
under the maps : 
\[
H^{2}(\hat{X}_{can},\mbox{\bf R})
\rightarrow \mbox{Hom}(H_{2}(\hat{W},\mbox{\bf R}),\mbox{\bf R}),
\]
and
\[
H^{2}(\hat{X}_{can},\mbox{\bf R})
\rightarrow \mbox{Hom}(H_{2}(\hat{X}_{can},\mbox{\bf R}),\mbox{\bf R}),
\]
are the images of elements of 
 $\mbox{Hom}(H_{2}(\hat{W},\mbox{\bf Z}),\mbox{\bf Z})$ 
and $\mbox{Hom}(H_{2}(\hat{X}_{can},\mbox{\bf Z}),\mbox{\bf Z})$ 
respectively.
This implies  that some positive multiple of 
$c_{1}(K_{\hat{X}_{can}})\in H^{2}(\hat{X}_{can},\mbox{\bf R})$ is integral
(i.e. it is in the image of the natural morsphim 
$H^{2}(\hat{X}_{can},\mbox{\bf Z})\rightarrow H^{2}(\hat{X}_{can},\mbox{\bf R})$) 
in $H^{2}(\hat{X}_{can},\mbox{\bf R})$. 
Hence some positive multiple of $K_{\hat{X}_{can}}$ is a line bundle on $\hat{X}_{can}$.

Let $r$ be a positive integer such that $rK_{\hat{X}_{can}}$ is a line bundle. 

\begin{definition}
Let $X$ be a normal complex space.  
We define the $L^{2}$-dualizing sheaf $K_{X}^{(2)}$ by
\[
K_{X}^{(2)}(U) = \{\eta \in \Gamma (U,{\cal O}_{X}(K_{X}))\mid 
\eta\wedge\bar{\eta} \in L^{1}_{loc}(U)\}.
\]
\end{definition}
The following lemma is clear by the definition of 
canonical singularities. 
\begin{lemma}
Let $X$ be a normal complex space with only canonical singularities.
Then the canonical sheaf  $K_{X} := i_{*}K_{X_{reg}}$ of $X$ is isomorphic to 
$K_{X}^{(2)}$, where $i : X_{reg}\longrightarrow X$ is the canonical injection. 
\end{lemma}
\begin{lemma} Let $Z$ be a closed $n$-dimensional subvariety of the unit open polydisk $\Delta^{N}$ with only canonical singularities 
and let $\varphi$ be a plusisubharmonic function on 
$Z \times \Delta$, where $\Delta$ is an open unit disk in ${\bf C}$.
Let $D$ be a {\bf Q}-Cartier divisor on $Z$ such that 
$K_{Z}+D$ is Cartier. 
Let $h_{D}$ be a $C^{\infty}$-hermitian metric on the 
{\bf Q}-line bundle ${\cal O}_{Z}(D)$. 
Let $t$ be the standard coordinate on $\Delta$.
Let $p_{1} : Z\times \Delta\longrightarrow Z$ be the first projection.
Then there exists a positive constant $C_{Z}$  depending only on $Z$ such that 
for every  $f\in \Gamma (Z,{\cal O}_{Z}(K_{Z}+D))$ such that 
\[
(\sqrt{-1})^{\frac{n(n-1)}{2}}\int_{Z}e^{-\varphi}\cdot h_{D}\cdot f\wedge\bar{f} < \infty
\]
there exists a holomorphic section  $F 
\in \Gamma (Z\times \Delta,{\cal O}_{Z\times\Delta}(K_{Z\times\Delta}+p_{1}^{*}D))$ such that 
\begin{enumerate}
\item $F\mid_{Z} = dt\wedge f$,
\item $(\sqrt{-1})^{\frac{n(n+1)}{2}}\int_{Z\times\Delta}e^{-\varphi}\cdot h_{D}\cdot F\wedge\bar{F}
\leq C_{Z}(\sqrt{-1})^{\frac{n(n-1)}{2}}\int_{Z}e^{-\varphi}\cdot h_{D}\cdot f\wedge\bar{f}$
\end{enumerate}
\end{lemma}
This lemma is an immediate consequence of the $L^{2}$-extension theorem 
 (\cite[p. 200, Theorem]{o-t}). 
Since $\hat{X}_{can} - \hat{W}$ is biholomorphic to $X - S$, it admits
a complete K\"{a}hler metric. 
Hence we can apply $L^{2}$-estimates for $\bar{\partial}$-operator on $\hat{X}_{can} -\hat{W}$.

We note that for every compact positive dimensional subvariety $V$ in $\hat{X}_{can}$, 
$(rK_{\hat{X}_{can}})^{\dim V}\cdot V : = \mu (V,rK_{\hat{X}_{can}}) \geq 1$.
Let $m_{0}$ be a positive integer such that 
$\mbox{Supp}\,\mbox{Bs}\mid m_{0}rK_{\hat{X}_{can}}\mid \subseteq \hat{W}$ holds. 
Such $m_{0}$ exists by Proposition 5.1. 
Let $\tau_{0}\ldots ,\tau_{N}$ be a basis of 
$\Gamma (\hat{X}_{can},{\cal O}_{\hat{X}_{can}}(m_{0}rK_{\hat{X}_{can}})$.
We define the singular hermitian metric on $rK_{\hat{X}_{can}}$ by 
\[
h_{0} := \frac{1}{(\sum_{i=0}^{N}\mid \tau_{i}\mid^{2})^{1/m_{0}}}.
\]
Let $x_{0}$ be an arbitrary point in   $\mbox{Supp}\,\mbox{Bs}\mid m_{0}rK_{\hat{X}_{can}}\mid$.
Then by Lemma 6.2, we see that for a local generator $\sigma$ of 
$rK_{\hat{X}_{can}}$ on a neighbourhood $U$ of $x_{0}$, 
for every $\alpha \geq m_{0}n$
the singular volume form 
\[
h_{0}^{\alpha}\mid\sigma\mid^{2\alpha +\frac{2}{r}}
\]
is not locally integrable on $U\cap (\hat{X}_{can}-\hat{W})$.
In fact let $x(t) (t\in\Delta )$ be a local holomorphic curve on $\hat{X}_{can}$ such that $x(0) = x_{0}$ and $x(t) \in \hat{X}_{can} -\hat{W}(t\in\Delta^{*})$.Then the limit $\lim_{t\rightarrow 0}n\cdot x(t)$ in the Douady space of $\hat{X}_{can}$ is contained in $n\cdot x(0)$. 
Hence by Lemma 6.2, the assertion follows. 

By Lemma 6.1 and Lemma 6.2 instead of Lemma 4.6, 
using the parallel argument as in \cite{a-s} or Section 4 
we conclude that 
$\mid m(rK_{\hat{X}_{can}})\mid$ is 
free at $x_{0}$ for every 
$m \geq m_{0}n+ n(n-1)/2 +2$ (we note that every strata constructed as in Section 4 except $\hat{X}_{can}$  is contained in $\hat{W}$, hence it is compact).
Since $m_{0}$ is independent of the choice of $x_{0}$, we see that 
$\mid m(rK_{\hat{X}_{can}})\mid$ is 
free on $\hat{X}_{can}$ for every 
$m \geq m_{0}(n(n+1)/2 +2)$.
Since $K_{\hat{X}_{can}}$ is numerically positive, we see that $\hat{W}$ consists of finitely many irreducible 
components and $\hat{X}_{can}$ is a projective variety 
(with only canonical singularities).
This implies that $\hat{X}_{can}$ is the canonical model
of $X$. 
Hence $R(X,K_{X})$ is finitely generated. 
This completes the proof of Theorem 1.1. 
 
Author's address\\
Hajime Tsuji\\
Department of Mathematics\\
Tokyo Institute of Technology\\
2-12-1 Ohokayama, Megro 152-8551\\
Japan \\
e-mail address: tsuji@math.titech.ac.jp
\end{document}